\documentclass[pdflatex,sn-apa]{sn-jnl}

\usepackage{graphicx}
\usepackage{amsmath,amssymb,amsfonts}
\usepackage{bm}
\usepackage{enumitem}
\usepackage{mathtools}
\usepackage{xcolor}
\usepackage{longtable,booktabs,array}
\usepackage{calc} 
\usepackage{chngcntr} 
\usepackage{multirow}
\let\orcidlogo\relax
\usepackage{orcidlink}
\usepackage{etoolbox}

\makeatletter
\patchcmd{\@maketitle}
{Corresponding author(s). E-mail(s):}
{Corresponding author. E-mail:}
{}{}
\makeatother
\hypersetup{colorlinks=true,citecolor=blue,linkcolor=blue,urlcolor=blue}

\theoremstyle{thmstyleone}

\newtheorem{theorem}{Theorem}[section]
\newtheorem{lemma}[theorem]{Lemma}
\newtheorem{corollary}[theorem]{Corollary}
\newtheorem{proposition}{Proposition}[section]
\theoremstyle{thmstyletwo}

\newtheorem{remark}{Remark}
\theoremstyle{thmstylethree}

\numberwithin{equation}{section}

\DeclareMathOperator*{\argmin}{argmin}

\newcommand{\abs}[1]{\left| #1 \right|}
\newcommand{\norm}[1]{\left\| #1 \right\|}
\newcommand{\mbf}[1]{\mathbf{#1}}
\newcommand{\rmf}[1]{\mathrm{#1}}
\newcommand{\tp}{\mathrm{tp}}
\newcommand{\BR}{\mathrm{BR}}
\newcommand{\WBR}{\mathrm{WBR}}
\newcommand{\PH}{\mathrm{PH}}

\newcommand{\meanin}{n^{-1}\sum_{i=1}^n}
\newcommand{\sumin}{\sum_{i=1}^n}
\newcommand{\sumjd}{\sum_{j=1}^d}
\newcommand{\sumkneqj}{\sum_{k\neq j}}
\newcommand{\sumellL}{\sum_{\ell=1}^L}
\newcommand{\suminl}{\sum_{\ell=1}^L}
\newcommand{\suminlp}{\sum_{\ell'=1}^L}
\newcommand{\prodjd}{\prod_{j=1}^d}

\newcommand{\E}{\normalfont \mathbb{E}}
\newcommand{\Var}{\normalfont \text{Var}}
\newcommand{\Cov}{\normalfont \text{Cov}}

\newcommand{\hatm}{\hat{m}}
\newcommand{\hatp}{\hat{p}}
\newcommand{\hatr}{\hat{r}}

\newcommand{\hatI}{\hat{I}}

\newcommand{\tildem}{\tilde{m}}

\newcommand{\bmb}{\bm{b}}
\newcommand{\bmf}{\bm{f}}
\newcommand{\bmg}{\bm{g}}
\newcommand{\bmh}{\bm{h}}
\newcommand{\bmm}{\bm{m}}
\newcommand{\bmu}{\bm{u}}
\newcommand{\bmw}{\bm{w}}
\newcommand{\bmx}{\bm{x}}

\newcommand{\bmN}{\bm{N}}
\newcommand{\bmX}{\bm{X}}

\newcommand{\bmeta}{\bm{\eta}}
\newcommand{\bmgamma}{\bm{\gamma}}
\newcommand{\bmkappa}{\bm{\kappa}}
\newcommand{\bmmu}{\bm{\mu}}

\newcommand{\bbG}{\mathbb{G}}
\newcommand{\bbI}{\mathbb{I}}
\newcommand{\bbR}{\mathbb{R}}
\newcommand{\bbW}{\mathbb{W}}

\newcommand{\cE}{\mathcal{E}}
\newcommand{\cF}{\mathcal{F}}
\newcommand{\cG}{\mathcal{G}}
\newcommand{\cI}{\mathcal{I}}
\newcommand{\cM}{\mathcal{M}}
\newcommand{\cN}{\mathcal{N}}
\newcommand{\cX}{\mathcal{X}}

\newcommand{\fraku}{\mathfrak{u}}

\newcommand{\frakU}{\mathfrak{U}}

\newcommand{\hatbmm}{\hat{\bm{m}}}
\newcommand{\hatbmM}{\hat{\mathbf{M}}}
\newcommand{\tildebmm}{\tilde{\bm{m}}}

\newcommand{\hatalpha}{\hat{\alpha}}
\newcommand{\hatbeta}{\hat{\beta}}

\renewcommand{\d}{{\rm d}}

\newcommand{\F}{{\rm F}}

\newcommand{\bmvec}{{\rm\bf{vec}}}

\newcommand{\Prob}{\normalfont \mathbb{P}}

\allowdisplaybreaks
\begin{document}

\title[Simultaneous inference for additive regression]{Simultaneous inference for additive regression}

\author{\fnm{Donghyun} \sur{Kim}\textsuperscript{1}}\email{dong1193@snu.ac.kr}

\author*{\fnm{Changwon} \sur{Choi}\textsuperscript{1*}\orcidlink{0009-0001-8053-8550}}\email{cwchoi94@gmail.com}

\affil[1]{\orgdiv{Department of Statistics}, \orgname{Seoul National University},
	\orgaddress{\city{Seoul}, \country{South Korea}}}

\abstract{
This paper develops methods for simultaneous inference on component functions in additive models.
We first propose bias-corrected smooth backfitting estimators inspired by least squares linear regression and establish their asymptotic properties.
We then show that the distributions of the uniform errors of these estimators converge to a Gumbel-type distribution.
Based on this distributional approximation, we construct simultaneous confidence bands and develop a practical residual-based procedure for testing the additive components.
Simulation studies and empirical data analyses demonstrate that the proposed estimators reduce bias and improve estimation accuracy compared with several existing estimators.
Moreover, the empirical coverage probabilities of the proposed simultaneous confidence bands are found to be comparable to those of their oracle counterparts.
}

\keywords{Simultaneous inference, simultaneous confidence band, additive models,
	smooth backfitting, bias reduction}

\maketitle

\section{Introduction}\label{sec:intro}

Nonparametric regression provides a flexible framework for capturing nonlinear structures and heterogeneous effects that commonly arise in real-world data, making it appealing in applications where the underlying functional form is unknown.
Among the most widely used kernel-based methods are the Nadaraya–Watson and local linear estimation, which approximate the regression function locally by a constant and by a linear function, respectively.
Despite their simplicity and favorable theoretical properties, such estimators typically suffer from the curse of dimensionality, leading to slow convergence rates and high variance in multivariate settings.
Moreover, these estimators exhibit a non-negligible asymptotic bias when an optimal bandwidth is used \citep{tsybakov2009introduction}, which poses challenges for statistical inference.

In view of these limitations, bias reduction has become an important research topic.
In kernel-based regression, several approaches have been developed to reduce or account for smoothing bias.
One such approach, undersmoothing, employs a bandwidth smaller than the optimal bandwidth to ensure that the bias is asymptotically negligible relative to the stochastic error \citep{chen2002confidence,fan2000simultaneous,zhang2010simultaneous}.
In contrast, explicit bias correction estimates the leading bias term and subtracts it from the original estimator, typically by employing a higher-order local polynomial or an additional smoothing parameter
\citep{xia1998bias,cheng2019nonparametric}.
Bootstrap calibration offers an alternative approach to inference in the presence of smoothing bias, without requiring undersmoothing or explicit bias correction \citep{hall2013simple}.
These approaches rely critically on data-driven bandwidth selection or calibration, which can substantially influence bias correction and statistical inference and often requires considerable empirical tuning.
Alternative methods aim to reduce bias directly through structural modifications to the estimator. Examples include techniques that exploit information from nearby design points \citep{choi1998bias,choi2000data},
and approaches based on model averaging via double smoothing \citep{he2009double}.
Although these methods can achieve nearly unbiased estimation, the complexity of their asymptotic variance structures can complicate variance estimation and potentially undermine the validity of subsequent statistical inference.
Recently, \cite{lee2020nonparametric} proposed a parametric-help approach that reduces bias by projecting the unknown regression function onto a parametric family of function spaces. Despite these advances,
such bias reduction methods remain subject to the curse of dimensionality in multivariate settings.

Bias reduction is also of central interest in sieve-based nonparametric regression, where two distinct sources of bias arise: approximation bias
resulting from the use of a finite number of basis functions and, under penalization, smoothing bias induced by the penalty term.
The asymptotic effects of the basis dimension and penalty parameter on these two sources of bias have been extensively studied \citep{li2008asymptotics,claeskens2009asymptotic}.
For smoothing splines, \cite{shang2013local} established valid pointwise and global inference via undersmoothing, rendering the regularization bias asymptotically negligible.
These studies provide important tools for understanding and addressing bias in spline-based regression.
Nevertheless, their focus on fully nonparametric regression leaves the fundamental challenge of dimensionality unresolved,
as their estimation accuracy and convergence rates can deteriorate substantially in high-dimensional settings.

Additive models provide a powerful framework for addressing the curse of dimensionality.
To estimate additive components, the smooth backfitting (SBF) technique \citep{mammen1999existence} is widely used because of its desirable theoretical properties
and coherent framework for statistical inference \citep{mammen2006simple,jeon2020additive}.
The SBF methodology has been adopted for estimating a variety of structured nonparametric models, including varying coefficient models \citep{park2015varying} and partially linear additive models \citep{yu2011semi}, highlighting its versatility and broad applicability in nonparametric and semiparametric regression. More recently, SBF has been extended to high-dimensional settings and non-Euclidean data structures
Bias reduction methods have also been developed for additive regression models. Building on the parametric-help approach of \cite{lee2020nonparametric}, \cite{hong2023additive}
proposed a parametrically assisted estimator for additive models, which was subsequently extended to varying-coefficient models by \cite{moon2025varying}.
However, the extent of bias reduction achieved by this framework depends critically on the choice of the assisting function.
Existing inferential procedures for additive models include pointwise methods \citep{gregory2021statistical,lee2024efficient} and an approach based on $L_2$-norm \citep{mammen2022backfitting}.
To the best of our knowledge, however, simultaneous inference for additive components has received little attention.

In this paper, we develop simultaneous inferential procedures for additive models.
For this, we first introduce bias-corrected estimators that do not require auxiliary functions, drawing on the (weighted) least squares approach \citep{cheng2018bias}.
The proposed estimators completely eliminate the leading asymptotic bias term, yielding faster pointwise, $L_2$, and uniform convergence rates
than the original local linear SBF estimator.
Using these bias-corrected estimators, we construct simultaneous confidence bands for the additive component functions.
Pointwise inference concerns an additive component at a fixed covariate value, whereas inference for the entire component function over an interval requires a simultaneous confidence band.
Moreover, the non-negligible asymptotic bias of existing estimators can result in improperly calibrated confidence intervals and hypothesis tests.
These considerations motivate the construction of simultaneous confidence bands and hypothesis tests based on the limiting distribution of the uniform estimation error of the proposed bias-corrected estimators.

Specifically, for each component function $m_j:[0,1]\to\bbR$ of the underlying additive model,
we consider testing the null hypothesis $H_0 : m_j(x_j) = m_{0j}(x_j) \text{ for all } x_j\in I$, where $m_{0j}:[0,1]\to\bbR$ is a known function and $I \subset (0,1)$ is a closed interval.
We establish that the distribution of the uniform estimation error for each additive component over $I \subset (0,1)$ converges to a Gumbel-type distribution, in a spirit similar to that of \cite{hall2013simple}.
Using this distributional approximation, we develop a practical procedure for conducting the simultaneous test.
We employ residual-based local linear smoothing to estimate the conditional variance of the bias-corrected SBF estimator, thereby enabling a feasible and data-driven testing procedure.

Numerical studies confirm the theoretical advantages of the proposed bias-corrected estimators of $m_j$.
Specifically, the proposed method outperforms not only standard local linear smoothing but also the parametric-help approach with naively chosen assisting functions.
Furthermore, the proposed simultaneous test achieves performance comparable to that of an oracle procedure in which all additive components other than the target component $m_j$ are known,
thereby reducing the problem to a standard univariate nonparametric test for $m_j$ \citep{cheng2018bias}.

The remainder of this paper is organized as follows.
Section~\ref{sec:methods} outlines the basic idea and properties of local linear smoothing for additive models,
introduces the least-squares-based (LS-based) bias reduction method, and establishes the asymptotic properties of the proposed bias-corrected estimators.
In Section~\ref{sec:inference}, we derive the uniform distributional approximation and develop the simultaneous test and confidence band. Simulation studies and real-data applications are presented in
Sections~\ref{sec:Numerical_Study} and \ref{sec:real_data}, respectively.
Section~\ref{sec:conclusion} concludes with brief remarks. Proofs of all theorems are provided in the Supplementary Material.

\section{Bias reduction for local linear SBF}\label{sec:methods}

\subsection{Local linear SBF}\label{subsec:LL_SBF}

Let $Y$ be a response variable and $\bmX=(X_1,\dots,X_d)^\top$ be a covariate vector taking values in $[0,1]^d$.
We observe i.i.d. copies $(\bmX_i,Y_i)$ of $(\bmX,Y)$.
Throughout, we assume that the distribution of $\bmX$ is absolutely continuous with respect to the $d$-dimensional Lebesgue measure, and hence $\bmX$ has a joint density function on $[0,1]^d$.
We consider the additive model:
\begin{equation}\label{eq:add_model}
	Y = m_0 + \sumjd m_j(X_j) + \varepsilon,
\end{equation}
where $m_0\in\bbR$ is an unknown constant, $m_j:[0,1]\to\bbR$ for $1\le j\le d$ are unknown univariate component functions, and $\varepsilon$ is a random error satisfying $\E(\varepsilon | \bmX) = 0$.
To ensure the identifiability, we assume that $\E[m_j(X_j)] = 0$ for $1\le j\le d$.

As a starting point, we consider the estimators of the additive components in \eqref{eq:add_model} proposed by \cite{jeon2022locally}.
For a baseline kernel $K:\bbR\to[0,\infty)$ and a bandwidth $h>0$, we use the normalized kernel
$K_h:[0,1]\times[0,1]\to[0,\infty)$ defined by
$K_h(u,v):= c_h(v)\cdot h^{-1} K((u-v)/h)$ for $u,v\in[0,1]$,
where $c_h(v)$ is the normalizing constant to satisfy
$\int_0^1 K_h(u,v)\d u=1$ for all $v\in[0,1]$.
Let $h_1,\dots,h_d>0$ denote the componentwise bandwidths.
For $1\le j\le d$, define the scaled derivative component as $m_{1,j}(x_j) := h_j m_j'(x_j)$ and $\bmm_j^\tp:=(m_j,m_{1,j})^\top:[0,1]\to\bbR^2$.
The details for the local linear SBF estimators $\hatbmm_j^\tp = (\hatm_j,\hatm_{1,j})^\top$, including their construction, identifiability constraints and iterative backfitting algorithm, are provided in Section~S.1 of the Supplementary Material.

We establish a uniform asymptotic expansion for the local linear SBF estimator, which serves as the foundation for our subsequent analysis.
In particular, the deterministic bias term in this expansion motivates the bias reduction method proposed in Section~\ref{subsec:BR_SBF}.
For this, we impose the following assumptions.
Let $p$, $p_{jk}$ and $p_j$ denote the density functions of $\bmX$, $(X_j,X_k)$ and $X_j$, respectively.
For positive sequences $a_n,b_n>0$, we write $a_n\lesssim b_n$ if there exist some constant $0<C<\infty$ such that $a_n \le Cb_n$ for all $n\ge1$.
We write $a_n\asymp b_n$ if $a_n\lesssim b_n$ and $b_n\lesssim a_n$.

\begin{enumerate}[label=(A\arabic*), leftmargin=*, align=left]
	
	\item
	The joint density $p$ is bounded away from zero and infinity on $[0,1]^d$.
	For each $1\le j\le d$, the derivative $p_j'$ exists and is continuous on $[0,1]$.
	Also, for each $1\le j\neq k\le d$, the partial derivatives $\frac{\partial}{\partial x_j}p_{jk}(x_j,x_k), \frac{\partial}{\partial x_k}p_{jk}(x_j,x_k)$ exist and are continuous on $[0,1]^2$.
	
	\item
	The baseline kernel $K$ is symmetric, bounded, Lipschitz continuous, supported on $[-1,1]$ and satisfies $\int K(u) \d u=1$.
	
	\item
	The random error $\varepsilon$ satisfies $\E[|\varepsilon|^\alpha]<\infty$ for some $\alpha>5/2$.
	Also, its conditional second moment function $\E(\varepsilon^2 | \bmX=\cdot)$ is bounded on $[0,1]^d$.
	
	\item
	For $1\le j\le d$, the true additive component $m_j$ has a continuous second derivative on $[0,1]$.
	
	\item
	There exist constants $0<c_j<(\alpha-2)/\alpha$ with $c_j + c_k < 1$ for all $1\le j\neq k\le d$, where $\alpha$ is the constant in (A3), such that $h_j=o(1)$, $n^{c_j}h_j\gtrsim1$ and $h_k^2\lesssim h_j$ for all $1\le j\neq k\le d$.
	
\end{enumerate}

Assumptions (A1)--(A5) are standard in the theoretical analysis of SBF estimators \citep{jeon2020additive,jeon2022locally,lee2024efficient}.
Note that (A5) accommodates the optimal bandwidths $h_j\asymp n^{-1/5}$.

For $1\le j\neq k\le d$, let $\hatbmM_{jj}:[0,1]\to\bbR^{2\times2}$ and $\hatbmM_{jk}:[0,1]^2\to\bbR^{2\times2}$ be defined by
\begin{align*}
	\hatbmM_{jj}(x_j) :=& \, \meanin \bmvec\bigg(\frac{X_{ij}-x_j}{h_j}\bigg) \bmvec\bigg(\frac{X_{ij}-x_j}{h_j}\bigg)^\top	K_{h_j}(x_j,X_{ij}), \\
	\hatbmM_{jk}(x_j,x_k) :=& \, \meanin \bmvec\bigg(\frac{X_{ij}-x_j}{h_j}\bigg)\bmvec\bigg(\frac{X_{ik}-x_k}{h_k}\bigg)^\top	K_{h_j}(x_j,X_{ij})	K_{h_k}(x_k,X_{ik}),
\end{align*}
where $\bmvec(u) := (1,u)^\top$ for $u\in\bbR$.
Also, let
\begin{equation}\label{eq:bmN_bmgama_def}
	\begin{split}
		\bmN_j(x_j) :=& \, \int_0^1 \bmvec\bigg(\frac{v-x_j}{h_j}\bigg) \bmvec\bigg(\frac{v-x_j}{h_j}\bigg)^\top	K_{h_j}(x_j,v)\, \d v, \\
		\bmgamma_j(x_j) :=& \, \int_0^1 \bigg(\frac{v-x_j}{h_j}\bigg)^2 \bmvec\bigg(\frac{v-x_j}{h_j}\bigg)	K_{h_j}(x_j,v)\,\d v.
	\end{split}	
\end{equation}
For $x_j \in [2h_j, 1-2h_j]$, the normalized kernel reduces to the ordinary kernel, so the associated kernel moments coincide with the complete kernel moments.
Consequently, $\bmN_j(x_j) = \bmN$ and $\bmgamma_j(x_j) = \bmgamma$, where $\bmN\in\bbR^{2\times 2}$ is a non-negative definite matrix and $\bmgamma \in\bbR^2$.
Define
\begin{equation*}
	\tildebmm_j^{A,\tp}(x_j) := \meanin \hatbmM_{jj}(x_j)^{-1} \bmvec\bigg(\frac{X_{ij}-x_j}{h_j}\bigg) K_{h_j}(x_j,X_{ij}) \cdot \varepsilon_i, \quad x_j\in[0,1].
\end{equation*}
The following proposition establishes a uniform asymptotic expansion of the local linear SBF estimator and is a slight modification of Theorem 2 in \cite{jeon2022locally}.

\begin{proposition}\label{prop:bias_expansion}
	Assume (A1)--(A5).
	Then for each $1\le j\le d$,
	\begin{align*}
		&\sup_{x_j\in[0,1]} \Big| \hatbmm_j^{\tp}(x_j) - \bmm_j^{\tp}(x_j) - \tildebmm_j^{A,\tp}(x_j) - \frac{1}{2}h_j^2\bmN_j(x_j)^{-1}\bmgamma_j(x_j)m_j''(x_j) \Big| \\
		&\hspace{1cm} = O_p\bigg( \sqrt{\frac{\log n}{n}}\bigg) + o_p\bigg(\sum_{k=1}^{d}h_k^2\bigg).
	\end{align*}
\end{proposition}

The stochastic term $\hatbmm_j^{A, \tp}(x_j)$ in Proposition~\ref{prop:bias_expansion} is of order $(nh_j)^{-1/2}$ in both the pointwise and $L_2$ senses,
and of order $(\log n)^{1/2}(nh_j)^{-1/2}$ uniformly over $x_j \in [0,1]$ \citep{jeon2022locally}.
Thus, if $h_j \asymp h_k$ for $1\le j\neq k\le d$, we obtain the following convergence rates:
\begin{itemize}	
	\item[(i)] Pointwise convergence: $| \hatm_j(x_j) - m_j(x_j) | = O_p\big( h_j^2 + (nh_j)^{-1/2} \big)$ for $x_j \in (0,1)$;
	
	\item[(ii)] $L_2$ convergence: $\int_0^1 | \hatm_j(x_j) - m_j(x_j) |^2 \d x_j = O_p\big(h_j^4 + (nh_j)^{-1} \big)$;
	
	\item[(iii)] Uniform convergence: $\sup_{x_j\in[0,1]} | \hatm_j(x_j) - m_j(x_j)| =  O_p\big(h_j^2 + \sqrt{\frac{\log n}{nh_j}} \big)$.
\end{itemize}

\subsection{Bias reduction}\label{subsec:BR_SBF}

To develop valid simultaneous inferential procedures for the additive components, we introduce bias reduction methods for the local linear SBF estimators
that eliminate the deterministic bias terms of size $h_j^2$ in Proposition~\ref{prop:bias_expansion}.
To specify the dependence of the estimators on the bandwidths, we let $\hatbmm_j^{\tp}(\cdot;\bmh)$ denote the local linear SBF estimators based on the bandwidth vector $\bmh = (h_1,\dots,h_d)$.
We then define the corresponding estimator of $m_j$ as $\hatm_j(x_j;\bmh) := \fraku_0^\top \cdot \hatbmm_j^{\tp}(x_j;\bmh)$ for all $x_j\in[0,1]$, where $\fraku_0 = (1,0)^\top$.

We employ the LS-based bias reduction approach proposed by \cite{cheng2018bias}.
According to Proposition~\ref{prop:bias_expansion}, for $x_j\in(0,1)$,
\begin{equation}\label{hatmj_expression}
	\begin{aligned}
		\hatm_j(x_j;\bmh)
		&\approx m_j(x_j) + \frac{1}{2}\fraku_0^\top \bmN_j(x_j)^{-1} \bmgamma_j(x_j) m_j''(x_j) \cdot h_j^2
		+ \fraku_0^\top \cdot \tildebmm_j^{A,\tp}(x_j).
	\end{aligned}
\end{equation}
Choose an integer $L \ge 2$ and a set of auxiliary bandwidth vectors $\{\bmb_{\ell}=(b_{1,\ell}, \ldots, b_{d,\ell}):1\le \ell\le L\}$.
For each $j=1,\ldots,d$, we fit the following \textit{Hilbertian} linear regression model,
treating $b_{j,\ell}^2$ and $\hatm_j(\,\cdot\,;\bmb_{\ell})$, $\ell=1,\ldots,L$, as the $L$ covariates and $L^2([0,1])$-valued responses, respectively.
\begin{equation}\label{eq:linear_regression}
	\hatm_j(\,\cdot\,;\bmb_\ell) = \alpha_j + b_{j,\ell}^2 \cdot \beta_j  + \zeta_{j,\ell}, \quad \ell=1,\dots,L,
\end{equation}
where $\alpha_j, \beta_j \in L^2([0,1])$ and
$\zeta_{j,\ell}$ is considered as a mean-zero noise variable taking values in $L^2([0,1])$.
The least squares estimators of $\alpha_j$ and $\beta_j$ are then given by
\begin{align}\label{alpha-beta-est}\begin{split}
	\hatalpha_{j,L} :=& \,
	\suminl \Bigg\{\frac{\suminlp b_{j,\ell'}^4 - \big(\suminlp b_{j,\ell'}^2 \big) b_{j,\ell}^2}{L\cdot \suminlp b_{j,\ell'}^4 - \big(\suminlp b_{j,\ell'}^2 \big)^2} \Bigg\} \cdot \hatm_j(\,\cdot\,;\bmb_\ell), \\
	\hatbeta_{j,L}	:=& \,
	\suminl \Bigg\{\frac{L \cdot b_{j,\ell}^2 - \big( \suminlp b_{j,\ell'}^2 \big)}{L \cdot \suminlp b_{j,\ell'}^4 - \big(\suminlp b_{j,\ell'}^2 \big)^2} \Bigg\} \cdot \hatm_j(\,\cdot\,;\bmb_\ell),
\end{split}
\end{align}
where, again, the sum ``$\sum$'' and scalar multiplications of $\hatm_j(\,\cdot\,;\bmb_\ell)$ are understood in the sense of addition and scalar multiplications in the Hilbert space $L^2([0,1])$.

Noting that $\alpha_j$ are parameters for the respective $m_j$, and $\beta_j$ for the respective factors of $h_j^2$ in the expansion \eqref{hatmj_expression},
we define the following bias-corrected SBF estimators of $m_j$:
\begin{equation}\label{eq:BR_def}
\hatm_j^\BR := \hatm_j(\,\cdot\,;\bmh) - h_j^2\cdot \hatbeta_{j,L}, \quad 1 \le j \le d.
\end{equation}
We note that the bandwidth vector $\bmh\equiv (h_1, \ldots, h_d)$ used in \eqref{eq:BR_def} is allowed to be different from $\bmb_{\ell}, \,1\le \ell\le L$, used in the construction of $\hatbeta_{j,L}$.

Alternatively, one may consider a weighted least squares estimator for $\beta_j$ in the Hilbertian linear model \eqref{eq:linear_regression}.
This approach can improve efficiency by accounting for the heteroscedasticity in the stochastic component of $\hatm_j(x_j;\bmb_\ell)$ across the bandwidth vectors $\bmb_\ell$.
Considering that the variance of $\hatm_j(\,\cdot\,;\bmb_\ell)$ is of order $n^{-1}b_{j,\ell}^{-1}$,
we use weights proportional to $b_{j,\ell}$ for the $\ell$th covariate-response pair $(b_{j,\ell},\,\hatm_j(\,\cdot\,;\bmb_\ell))$.
The resulting weighted least squares estimator is given by
\begin{align*}
	\hatbeta_{j,L}^W :=& \,
	\suminl	\Bigg\{	\frac{\big( \suminlp b_{j,\ell'} \big) b_{j,\ell}^3 - \big(\suminlp b_{j,\ell'}^3 \big) b_{j,\ell} }{\big(\suminlp b_{j,\ell'} \big) \big(\suminlp b_{j,\ell'}^5\big) - \big(\suminlp b_{j,\ell'}^3\big)^2}\Bigg\} \cdot \hatm_j(\,\cdot\,;\bmb_\ell).
\end{align*}
The corresponding weighted versions of $\hatm_j^\BR$ are
\begin{equation}\label{eq:WBR_def}
\hatm_{j}^\WBR := \hatm_j(\,\cdot\,;\bmh) - h_j^2 \cdot \hatbeta_{j,L}^W, \quad 1 \le j \le d.
\end{equation}

\begin{remark}\label{alt-alpha}
The estimator $\hatalpha_{j,L}$ in \eqref{alpha-beta-est} can be used as an alternative bias-corrected estimator of $m_j$.
These estimators have a similar degree of bias reduction to that of $\hatm_j^\BR$, and their finite-sample performance is also comparable.
For more details, we refer to Section~S.2 of the Supplementary Material.
\end{remark}

\subsection{Asymptotic properties of bias-corrected estimators}\label{subsec:theory_BR_SBF}

We now establish the asymptotic properties of the two bias-corrected SBF estimators defined in \eqref{eq:BR_def} and \eqref{eq:WBR_def}.
Let
\[
b_j := \max_{1 \le \ell \le L}b_{j,\ell}, \quad b_{-j} := \max_{1 \le \ell \le L}\max\{b_{k,\ell}: k \neq j\}.
\]
Also, define $h_{-j} := \max_{k \neq j} h_k$ and
\begin{equation*}
	r_{nj} := h_j^2\frac{b_{-j}^2}{b_{j}^2} + h_{-j}^2, \quad 1\le j\le d.
\end{equation*}
We make an assumption on $b_j$ and $b_{-j}$ analogous to Assumption (A5) in Section~\ref{subsec:LL_SBF}.

\begin{enumerate}[label=(A\arabic*), start=6, leftmargin=*, align=left]
	
	\item
	For $1\le \ell\le L$ and $1\le j\neq k\le d$, it holds that $b_{j,\ell} / b_j \to c_{j,\ell}>0$ and $b_{k,\ell}\asymp b_{-j}$.
	For $1\le j\le d$, there exist constants $0<c_j', c_{-j}' < (\alpha-2)/\alpha$ with $c_j'+c_{-j}'<1$, where $\alpha$ is the constant in (A3), such that $b_j=o((\log n)^{-1})$, $n^{c_j'}b_j\gtrsim 1$, $n^{c_{-j}'}b_{-j}\gtrsim 1$.
\end{enumerate}

Let $I_j :=  [2b_j, 1-2b_j]$ denote the interior of the domain of $m_j$.

\begin{theorem}\label{thm:mhat_asymp}
	Assume that (A1)--(A6) hold and that $h_j = o(b_j)$.
	Then, for $\#=\BR$ and $\#=\WBR$, the following convergence rates hold:
	\begin{itemize}
		\item[(i)] Pointwise convergence:
		\begin{align*}
			| \hatm_j^{\#}(x_j) - m_j(x_j) | &= O_p\big(n^{-1/2}h_j^{-1/2}\big) + o_p\big(h_j^2 + r_{nj}\big), \quad x_j\in I_j, \\
			| \hatm_j^{\#}(x_j) - m_j(x_j) | &= O_p\big(n^{-1/2}h_j^{-1/2} + h_j^2\big) + o_p\big(r_{nj}\big), \quad x_j\in [0,1]\setminus I_j.
		\end{align*}
		
		\item[(ii)] $L_2$ convergence:
		\begin{align*}
			\int_0^1 | \hatm_j^{\#}(x_j) - m_j(x_j) |^2 \d x_j &= O_p\big(n^{-1}h_j^{-1}\big) + o_p\big(h_j^4 + r_{nj}^2 \big).
		\end{align*}
		
		\item[(iii)] Uniform convergence:
		\begin{align*}
			\sup_{x_j \in I_j} | \hatm_j^{\#}(x_j) - m_j(x_j) | &= O_p\bigg(\sqrt{\frac{\log n}{n h_j}} \bigg) + o_p\big(h_j^2 + r_{nj} \big), \\
			\sup_{x_j\in[0,1]} | \hatm_j^{\#}(x_j) - m_j(x_j) | &= O_p\bigg(\sqrt{\frac{\log n}{n h_j}} + h_j^2 \bigg) + o_p\big(r_{nj} \big).
		\end{align*}	
	\end{itemize}
\end{theorem}

Theorem~\ref{thm:mhat_asymp} shows that the bias of $\hatm_j(\,\cdot\,;\bmh)$, which is of order $O_p(h_j^2)$, is reduced to $o_p(h_j^2)$ for $\hatm_{j}^\BR$ and $\hatm_{j}^\WBR$,
while the convergence rate of the stochastic part is preserved.
The bias reduction achieved by $\hatm_{j}^\BR$ and $\hatm_{j}^\WBR$ facilitates the statistical inferential procedures developed in Section~\ref{sec:inference}.
We remark that Theorem~\ref{thm:mhat_asymp} holds for any fixed $L\ge 2$; in particular, no assumption that $L \to \infty$ is required.

To examine the contribution of the estimation error in $\hat{\beta}_{j,L}$ to the bias-corrected estimators, let $w_{j,\ell}$ denote the weights in the definition of
$\hat{\beta}_{j,L}$ in \eqref{alpha-beta-est}, i.e., $\hat{\beta}_{j,L}=\sum_{\ell=1}^L w_{j,\ell}\cdot\hatm_j(\,\cdot\,;\bmb_\ell)$.
We focus on $\hatm_{j}^\BR$, and neglect the $\sqrt{\log n}$ factor in the case of the uniform error.
Since $\sum_{\ell=1}^L w_{j,\ell}=0$ and $\sum_{\ell=1}^L w_{j,\ell}\,b_{j,\ell}^2=1$, we get
\begin{align*}
\hat{\beta}_{j,L}=\suminl w_{j,\ell} \cdot (\alpha_j + b_{j,\ell}^2 \cdot \beta_j + \zeta_{j,\ell})=\beta_j+\suminl w_{j,\ell}\cdot\zeta_{j,\ell}.
\end{align*}
Thus, the contribution of $\hat{\beta}_{j,L}$ to the error of $\hatm_{j}^\BR$ enters through the stochastic terms $\zeta_{j,\ell}$, which are of order $n^{-1/2}b_j^{-1/2}$.
Since $w_{j,\ell} \asymp b_j^{-2}$, the overall contribution of $\hat{\beta}_{j,L}-\beta_j$ is of order $n^{-1/2}b_j^{-5/2}h_j^2$,
which negligible relative to $n^{-1/2}h_j^{-1/2}$ under the assumption $h_j=o(b_j)$ imposed in Theorem~\ref{thm:mhat_asymp}.


\begin{remark}
We note that the interior $I_j$ in Theorem~\ref{thm:mhat_asymp} for the pointwise and uniform convergence rates is $[2b_j, 1-2b_j]$, instead of $[2h_j, 1-2h_j]$.
The main reason is that the linear model \eqref{eq:linear_regression} is not valid under the sup-metric, that is, in $L^\infty([0,1])$,
because $\bmN_j(x_j)$ and $\bmgamma_j(x_j)$ defined in \eqref{eq:bmN_bmgama_def}, and hence $\beta_j(x_j)$, depend on $b_{j,\ell}$ (the square-root of the covariate values) for $x_j \notin I_j$.
To obtain a version of Theorem~\ref{thm:mhat_asymp} with $I_j=[2h_j, 1-2h_j]$, the linear model \eqref{eq:linear_regression} may be re-formulated in a pointwise manner.
Specifically, for each $x_j \in [0,1]$, one fits
\begin{equation*}
\hatm_j(x_j;\bmb_\ell) = \alpha_j(x_j) + \beta_j(x_j) Z_j(x_j;b_{j,\ell}) + \zeta_{j,\ell}(x_j), \quad \ell=1,\dots,L,
\end{equation*}
now treating $Z_j(x_j;b_{j,\ell}) := b_{j,\ell}^2 \,\fraku_0^\top \, \bmN_j(x_j;b_{j,\ell})^{-1}\bmgamma_j(x_j;b_{j,\ell})$ and
$\hatm_j(x_j;\bmb_\ell)$ as the real-valued covariates and responses, respectively.
Although the underlying principle is the same, establishing its theoretical properties requires more involved derivations. Moreover, its practical implementation is considerably more complex
than that of the $L^2([0,1])$ version. Therefore, we do not explore this option in this paper.
\end{remark}

\begin{remark}
We compare our method with the parametric-help (PH) method \citep{hong2023additive}.
The effectiveness of the PH method depends critically on the choice of the assisting functions,
and it reduces only the constant factor of the bias of $\hatm_j(\,\cdot\,;\bmh)$ without changing its order.
In contrast, our method reduces the order of the bias without requiring such assisting functions.
This feature is beneficial not only for improving estimation accuracy but also for conducting statistical inference on the additive components.
\end{remark}

\section{Simultaneous inference}\label{sec:inference}

Classical kernel-based approaches to uniform inference typically rely on explicit bias correction or bootstrap calibration.
For example, \cite{hall2013simple} uses bootstrap calibration to achieve the desired coverage, while \cite{cheng2019nonparametric} employs debiased local polynomial estimation.
However, both approaches are developed only for univariate regression.

In the additive regression setting, we establish a distributional approximation for the supremum of the estimation errors of the bias-corrected local linear SBF estimators defined in \eqref{eq:BR_def} and \eqref{eq:WBR_def}.
These estimators remove the leading interior smoothing bias, thereby reducing the deterministic component of the estimation error that would otherwise complicate the uniform approximation.
Based on this approximation, we construct simultaneous confidence bands for each component function of the additive model.

\subsection{Uniform distributional approximation}\label{subsec:uniform_CLT}

We first derive the asymptotic distributions of the uniform estimation errors of $\hatm_j^\BR$ and $\hatm_j^\WBR$ defined in \eqref{eq:BR_def} and \eqref{eq:WBR_def}, respectively, over a closed interval $I \subset(0,1)$.
Recall that $p_j$ denotes the marginal density function of $X_j$.
Let $\nu_0 := \int_{-\infty}^{\infty} K^2(u)du$ and $\sigma_j^2(x_j):=\E(\varepsilon^2| X_j=x_j)$.
Define
\begin{equation*}
	A_{nj}(I):=(2\log (|I|/h_j))^{1/2}, \quad B_{nj}(I) := A_{nj}(I) + \frac{\log(\lambda_K/(4\pi^2))}{2A_{nj}(I)},
\end{equation*}
where $\lambda_K	:=	\nu_0^{-1} \cdot \int \{K'(u)\}^2\d u$.
For $\#\in\{\BR,\WBR\}$, let
\begin{equation}\label{sup-error}
\Delta_{nj}^\#(I) := \, \nu_0^{-1/2} \cdot \sup_{x_j\in I} \Big|\{nh_j\sigma_j(x_j)^{-2}p_j(x_j)\}^{1/2} \{\hat m_j^{\#}(x_j)-m_j(x_j)\}\Big|.
\end{equation}

Below, we present the asymptotic distributions of $\Delta_{nj}^\BR(I)$ and $\Delta_{nj}^\WBR(I)$.
For this we impose the following additional assumptions.

\begin{enumerate}[label=(B\arabic*), start=1, leftmargin=*, align=left]
\item The density $p$ satisfies (A1), and is continuously differentiable on $[0,1]^d$.
\item For some $\alpha>4$, it holds that $\E[|\varepsilon|^{\alpha}]<\infty$.
Also, the conditional variance function $\sigma^2:=\E(\varepsilon^2|\bmX=\cdot)$ is bounded away from zero
and infinity on $[0,1]^d$ and has continuous partial derivatives on $[0,1]^d$.
\item For some $0<\epsilon<1/5$ and $(1+2\epsilon)/7<c_h<1/2$, it holds that $h_j,h_{-j}\asymp n^{-c_h}$ and $b_j,b_{-j} \asymp n^{-c_h+\epsilon}$.
\item For $1 \le j \le d$, the true additive component $m_j$ has a Lipschitz continuous second derivative on $[0,1]$.
\end{enumerate}

The condition (B2) ensures $1/2 < (\alpha-2)/\alpha$, which, together with (B3), yields 
\[
c_h < \frac{1}{2}<\frac{\alpha-2}{\alpha}.
\]
Thus, setting $c_j=c_h$ for $1 \le j \le d$ satisfies the bandwidth rate restrictions in (A5).
Similarly, setting $c_j'=c_h-\epsilon$ for $1 \le j \le d$ satisfies those in (A6).
The condition (B3) allows the optimal bandwidth rate $n^{-1/5}$ for $h_j$ and $h_{-j}$ since $(1+2\epsilon)/7<1/5$.
The lower bound for $c_h$ ensures that $\sqrt{nh_j}A_{nj}(I) \cdot h_j^2 b_j\to0$, while the upper bound for $c_h$ ensures that $nh_jh_k\to\infty$ for $j\neq k$, which is required for joint estimation of the $d$ component functions.

\begin{theorem}\label{thm:unif_CLT}
	Assume that (A2), (A6) and (B1)--(B4) hold and that $h_j=o(b_j)$.
	In (A6), the constant $\alpha$ is taken to be the one appearing in (B2).
	Then, for any closed interval $I\subset(0,1)$ and $\#\in\{\BR,\WBR\}$,
	\begin{equation*}
		\lim_{n \to \infty}	\Prob\Big(A_{nj}(I)\big\{\Delta_{nj}^{\#}(I)-B_{nj}(I) \big\} \le u \Big) =\exp(-2e^{-u}), \quad u\in\bbR.
	\end{equation*}
\end{theorem}

Theorem~\ref{thm:unif_CLT} establishes that the $A_{nj}(I)$-scaled uniform estimation errors $\Delta_{nj}^\BR(I)$ and $\Delta_{nj}^\WBR(I)$ converge in distribution to a Gumbel-type distribution.
The proof relies on the following approximation strategy.
Conditional on the covariates, the leading stochastic terms of the bias-corrected SBF estimators $\hatm_j^\BR$ and $\hatm_j^\WBR$ are approximated by kernel-weighted sums of conditionally independent errors.
These weighted sums are then approximated by appropriate Gaussian processes.
The two approximations yield the Gumbel-type limit established in Theorem~\ref{thm:unif_CLT}.

\subsection{Uniform testing}\label{subsec:inference}

Building on Theorem~\ref{thm:unif_CLT}, we introduce a simultaneous testing procedure for $m_j$ in the additive model \eqref{eq:add_model}.
We describe the method using $\hatm_j^\BR$. The same method and results apply to its weighted counterpart $\hatm_j^\WBR$.
Let $m_{0j}:[0,1]\to\bbR$ denote a prespecified function and $I\subset(0,1)$ be a closed interval.
We consider the following hypothesis:
\begin{equation}\label{eq:hypothesis}
	H_0 : m_j(x_j) = m_{0j}(x_j)  \text{ for all } x_j \in I \quad \text{vs} \quad  H_1 : \text{not } H_0.
\end{equation}

To test the hypothesis \eqref{eq:hypothesis}, we use the asymptotic distribution of $\Delta_{nj}^\BR(I)$ established in Theorem~\ref{thm:unif_CLT}, which requires estimating the unknown quantities $p_j$ and $\sigma_j$.
For the marginal density $p_j$, we employ the kernel density estimator $\hatp_j$ defined by
\begin{equation*}
	\hatp_j(x_j) := \meanin K_{h_j}(x_j,X_{ij}), \quad x_j\in[0,1].
\end{equation*}
To estimate the conditional variance function $\sigma_j^2:=\E(\varepsilon^2|X_j=\cdot)$, we apply a local linear approach.
For an estimator $\hatm_j$ of $m_j$, let
\[
\hatr_{ij} := \{Y_i - (\bar{Y} + \sumjd \hatm_j(X_{ij}))\}^2, \quad 1\le i\le n.
\]
For each $x_j\in(0,1)$, we estimate $\sigma_j^2(x_j)$ by $\hat\sigma_j^2(x_j) := \widehat{\alpha}_j$, where
\begin{equation*}
(\widehat{\alpha}_j, \widehat{\beta}_j) := \argmin_{\alpha, \beta \in \bbR} \meanin K_{h_j}(x_j, X_{ij}) \big\{ \hatr_{ij} - \alpha - \beta (X_{ij} - x_j) \big\}^2.
\end{equation*}
We take $\hatm_j$ to be the local linear SBF estimator, $\hatm_j(\,\cdot\,;\bmh)$.
Alternatively, one may use the bias-corrected estimators introduced in Section~\ref{subsec:BR_SBF}.

We let $\widehat\Delta_{nj}^{\#}$ be defined as $\Delta_{nj}^{\#}$ in \eqref{sup-error} with $p_j$ and $\sigma_j^2$ replaced by $\hatp_j$ and $\hat\sigma_j^2$, respectively.
For $\#\in\{\BR,\WBR\}$, define the test statistic $T^{\#}_{nj}(I)$ by
\begin{equation}\label{eq:test_statistic}
	T^{\#}_{nj}(I) := A_{nj}(I)\big\{\widehat\Delta_{nj}^{\#}(I)-B_{nj}(I) \big\}.
\end{equation}
Since $A_{nj}(I)\Delta_{nj}^{\#}(I) \lesssim \log n$ and
\begin{align*}
\Big|A_{nj}(I)\big\{\widehat\Delta_{nj}^{\#}(I)-\Delta_{nj}^{\#}(I)\big\}\Big| \le A_{nj}(I)\Delta_{nj}^{\#}(I)\cdot \sup_{x_j \in I}\bigg|\sqrt{\frac{\hatp_j(x_j)\sigma(x_j)^2}{p_j(x_j)\hat\sigma(x_j)^2}}-1\bigg|,
\end{align*}
we obtain the following corollary.


\begin{corollary}\label{cor:unif_CLT}
	Assume that the conditions of Theorem~\ref{thm:unif_CLT} hold.
	Then, for any closed interval $I\subset(0,1)$ and $\#\in\{\BR,\WBR\}$,
	\begin{equation*}
		\lim_{n \to \infty}	\Prob\Big(T^{\#}_{nj}(I) \le u \Big) =\exp(-2e^{-u}), \quad u\in\bbR.
	\end{equation*}
\end{corollary}

Based on this approximation, given a significance level $\alpha\in(0,1)$, we reject the null hypothesis $H_0$ in \eqref{eq:hypothesis} if
\begin{equation}\label{eq:reject_test}
	T^{\#}_{nj}(I) > -\log\bigg(\frac{1}{2}\log \frac{1}{1-\alpha}\bigg) =: C_\alpha.
\end{equation}
Accordingly, the p-value is given by
\begin{equation*}
	\text{p-value} := 1-\exp\big(-2\exp(-T^{\#}_{nj}(I))\big).
\end{equation*}
Furthermore, Corollary~\ref{cor:unif_CLT} enables us to construct a simultaneous confidence band (SCB) for $m_j$ over $I \subset (0,1)$ satisfying
\begin{align*}
	\Prob\Big( m_j(x_j) \in \text{SCB}(x_j;\alpha,I) \text{ for all } x_j \in I \Big) \to 1-\alpha.
\end{align*}
Here, the SCB is defined by
\begin{equation}\label{asymp_scb}
	\text{SCB}(x_j;\alpha,I) := \big[\hatm_j^\BR(x_j) - \Delta_{1-\alpha}(x_j;I), \hatm_j^\BR(x_j) + \Delta_{1-\alpha}(x_j;I) \big], \quad x_j \in I,
\end{equation}
where
\begin{equation*}
	\Delta_{1-\alpha}(x_j;I) := \bigg(\frac{\nu_0 \hat{\sigma}_j^2(x_j)}{nh_j\hatp_j(x_j)}\bigg)^{1/2} \cdot \bigg( \frac{C_\alpha}{A_{nj}(I)} + B_{nj}(I)\bigg), \quad x_j\in I.
\end{equation*}

We compare our proposed simultaneous testing procedure with existing methods for additive regression.
First, \cite{lee2024efficient} introduced a debiased SBF estimator but considered only pointwise inference, and their method excludes the optimal bandwidth.
\cite{gregory2021statistical} proposed an additive estimator that is asymptotically equivalent to the oracle estimator; however, its non-negligible bias complicates statistical inference.
Furthermore, the $L_2$-based test proposed by \cite{mammen2022backfitting} applies only to testing $H_0: m_j \equiv 0$, for which the asymptotic bias vanishes.
In contrast, our method accommodates the optimal bandwidth and allows testing the null hypothesis $H_0$ in \eqref{eq:hypothesis} for any known function $m_{0j}:[0,1]\to\bbR$.
Thus, our approach provides a more flexible framework for inference on additive component functions.

\section{Simulation Studies}\label{sec:Numerical_Study}

\subsection{Estimation}\label{subsec:bias_reduction_sim}

Here, we assess the finite-sample performance of the proposed bias reduction methods.
We first introduce the data generating process. We choose $d=2$.
The covariates are $X_j := \Phi(Z_j)$ for $j=1,2$, where $\Phi:\bbR\to[0,1]$ denotes the cumulative distribution function of the standard normal distribution, and $(Z_1,Z_2)^\top$ is a bivariate normal random vector with $\E(Z_j)=0$ and $\Cov(Z_j,Z_k) = 0.5^{|j-k|}$.
By this construction, each $X_j$ is marginally uniformly distributed on $[0,1]$ while retaining dependence between them.
For some selected parameters $\rho_1,\rho_2\in\bbR$, we consider $m(x_1,x_2;\rho_1,\rho_2)=m_1(x_1;\rho_1)+m_2(x_2;\rho_2)$, where
\begin{align*}
m_1(x_1;\rho_1)&=\sin(2\pi x_1) +  \rho_1\exp(-16(4x_1-2)^2), \\
m_2(x_2;\rho_2)&= \cos(2\pi x_2)+\rho_2 \cdot \big[ 0.3\exp(-16(4x_2-1)^2) + 0.7\exp(-16(4x_2-3)^2) \big].
\end{align*}
The response variable $Y$ is then generated from the following two-dimensional additive model:
\begin{equation}\label{eq:sim_model}
	Y = m(X_1,X_2;\rho_1,\rho_2) + \sigma(X_1,X_2) \, \varepsilon,
\end{equation}
where $\sigma(x_1, x_2) = 0.2(1+\sin(2\pi x_1) + x_2^2)$ and the random error $\varepsilon\sim N(0,1)$ is independent of $(X_1,X_2)$.
We choose $(\rho_1,\rho_2)=(0,3), \, (2,0)$ and $(2.3)$.
These values are chosen so that the two components in each $m_j(\cdot,\rho_j)$ make a balanced contribution.

For comparison, we compute the local linear SBF estimator $\hatm_j^{\rm LL} := \fraku_0^\top \cdot \hatbmm_j^\tp(\cdot;\bmh)$, where $\hatbmm_j^\tp$ solves the SBF equation (S.1) subject to the constraint (S.2).
We also compute the parametric-help (PH) SBF estimator \citep{hong2023additive}.
Recall that the PH method requires the choice of the assisting functions $g_j$.
We considered a PH strategy that constructs the estimator $\hatm_{\rmf{naive}}^{\rmf{PH}}$ using the periodic components as assisting functions:
\begin{align}\label{assist-ftn}
g_1^{\rm naive}(x_1) := \sin(2\pi x_1), \quad g_2^{\rm naive}(x_2) := \cos(2 \pi x_2).
\end{align}

To select the bandwidths for the proposed bias-reduced estimation, we first obtain the optimal bandwidth $h_{j, \mathrm{opt}}$ for each $j$ using the coordinate-wise bandwidth selection (CBS) algorithm \citep{jeon2020additive}.
Specifically, given an initial set $(h_1^{(0)},\dots,h_d^{(0)})$, we iteratively update the bandwidths using the most recent optimal values, i.e.,
\begin{equation*}
	h_j^{(\ell)} := \argmin_{g_j\in\cG_j} \text{CV}\big(h_1^{(\ell)},\ldots,h_{j-1}^{(\ell)},g_j,h_{j+1}^{(\ell-1)},\ldots,h_d^{(\ell-1)}\big), \quad 1\le j\le d, \ \ell\ge1,
\end{equation*}
where $\cG_j:= \{ h_{j,\min} + 0.005 \cdot s :  s=1,\dots,10 \}$ and $h_{j,\min}$ is the minimal bandwidth chosen to satisfy $\inf_{x_j\in[0,1]}\lambda_{\min}(\hatbmM_{jj}(x_j))>0$ \citep{jeon2022locally}.
We use the 5-fold cross-validation prediction error defined by
\begin{equation}\label{CV-crit}
	\text{CV}(g_1,\dots,g_d) := \frac{1}{5} \sum_{k=1}^{5} \frac{1}{|I_k|} \sum_{i\in I_k} \big(Y_i - \widehat{m}^{(-k)}(\bmX_i;\bmg) \big)^2,
\end{equation}
where $I_1,\ldots,I_5$ form a partition of the index set $\{1,\ldots,n\}$ into five disjoint folds, and $\widehat m^{(-k)}(\cdot;\bmg)$ denotes the local linear SBF estimator computed with a bandwidth vector $\bmg=(g_1,\ldots,g_d)$ using all observations except those in the $k$-th fold.
For the auxiliary bandwidths, we set $b_{j,\ell} = (1 + C_{\ell} / 20)h_{j, \mathrm{opt}}$ for $\ell=1,\ldots,13$, where $C_\ell = \ell - 7$.
We then compute the proposed bias-corrected estimators using the auxiliary bandwidth vectors $\bmb_{\ell} = (b_{1,\ell}, \dots, b_{d,\ell})^{\top}$ for $\ell = 1, \ldots, 13$
along with the optimal bandwidth vector $\bmh = (h_{1,\mathrm{opt}}, \ldots, h_{d,\mathrm{opt}})$.

To select the bandwidths for the PH estimator, we also employ the CBS scheme with the 5-fold cross-validation in \eqref{CV-crit}, now replacing $\widehat{m}^{(-k)}(\bmX_i;\bmg)$
with $\widehat{m}_{\rm{PH}}^{(-k)}(\bmX_i;\bmg)$, 
where $\widehat{m}_{\rm{PH}}^{(-k)}(\cdot;\bmg)$ denotes the PH estimator computed using all observations except those in the $k$-th fold.
Throughout this paper, we use the Epanechnikov kernel $K(u) = (3/4)(1-u^2)\cdot\bbI(|u|<1)$ as the baseline kernel function.

To evaluate finite-sample performance, we compute Monte Carlo approximations of the integrated squared bias (ISB), the integrated variance (IV), and the mean integrated squared error (MISE),
for the estimators of each additive component $m_j$ as well as for the additive regression function $m$.
For $i=1,\dots,M$, let $\hatm^{(i)}$ denote the estimator of the additive regression function $m$ obtained from the $i$-th Monte Carlo replication.
We choose $M=100$. Then, the ISB, IV and MISE for $m$ are given by
\begin{align*}
	\text{ISB} & \coloneqq \int_{[0,1]^2} \bigg( M^{-1}\sum_{i=1}^M\hatm^{(i)}(\bmx) - m(\bmx) \bigg)^2 d\bmx,\\
	\text{IV} & \coloneqq \int_{[0,1]^2} M^{-1}\sum_{i=1}^M\bigg(\hatm^{(i)}(\bmx) - M^{-1}\sum_{i=1}^M \hatm^{(i)}(\bmx) \bigg)^2 d\bmx, \\
	\text{MISE} & \coloneqq \text{ISB} + \text{IV}.
\end{align*}
Similarly, for $j=1,2$, let $\hatm_j^{(i)}$ denote an estimator of $m_j$ obtained from the $i$-th Monte Carlo replication.
Then, the corresponding criteria for $m_j$ are given by
\begin{align*}
	\text{ISB}_j & \coloneqq \int_0^1 \bigg( M^{-1}\sum_{i=1}^M\hatm_j^{(i)}(x_j) - m_j(x_j) \bigg)^2 \d x_j, \\
	\text{IV}_j & \coloneqq \int_0^1 M^{-1}\sum_{i=1}^M\bigg(\hatm_j^{(i)}(x_j) - M^{-1}\sum_{i=1}^M \hatm_j^{(i)}(x_j) \bigg)^2 \d x_j, \\	
	\text{MISE}_j & \coloneqq \text{ISB}_j + \text{IV}_j.
\end{align*}

\begingroup
\setlength{\aboverulesep}{0pt}
\setlength{\belowrulesep}{0pt}
\renewcommand{\arraystretch}{0.8}
\begin{table}[h!]
	\centering
	\caption{The values of ISB, IV and MISE, multiplied by $10^3$, for the local linear SBF estimator, the PH estimator, and the proposed bias-corrected estimators.}\label{tab:MISE_overall}%
	\small
	\renewcommand{\arraystretch}{0.9}
	\setlength{\tabcolsep}{5pt}
	\begin{tabular}{@{}ccccccc@{}}
		\toprule
		\rule{0pt}{3.2ex}$(\rho_1,\rho_2)$ & $n$ &
		& $\hatm^{\rmf{LL}}$ & $\hatm_{\rmf{naive}}^{\rmf{PH}}$ & $\hatm^{\BR}$ & $\hatm^{\rmf{WBR}}$ \\[2pt]
		\midrule
		
		\multirow{6}{*}{(0,3)}
		& \multirow{3}{*}{400}
		& ISB    & 1.89 & 1.04 & 0.68 & 0.67 \\
		& & IV   & 14.00 & 14.60 & 14.80 & 14.70 \\
		& & MISE & 15.89 & 15.64 & 15.48 & 15.36 \\
		
		\cmidrule(lr){2-7}
		
		& \multirow{3}{*}{1000}
		& ISB    & 1.65 & 0.44 & 0.19 & 0.19 \\
		& & IV   & 5.36 & 6.17 & 6.20 & 6.17 \\
		& & MISE & 7.02 & 6.61 & 6.39 & 6.36 \\
		
		\cmidrule(lr){1-7}
		
		\multirow{6}{*}{(2,0)}
		& \multirow{3}{*}{400}
		& ISB    & 3.95 & 1.30 & 0.76 & 0.74 \\
		& & IV   & 13.95 & 14.76 & 14.92 & 14.79 \\
		& & MISE & 17.90 & 16.06 & 15.69 & 15.53 \\
		
		\cmidrule(lr){2-7}
		
		& \multirow{3}{*}{1000}
		& ISB    & 3.72 & 0.47 & 0.29 & 0.29 \\
		& & IV   & 4.46 & 4.77 & 4.76 & 4.74 \\
		& & MISE & 8.18 & 5.24 & 5.04 & 5.03 \\
		
		\cmidrule(lr){1-7}
		
		\multirow{6}{*}{(2,3)}
		& \multirow{3}{*}{400}
		& ISB    & 5.71 & 1.47 & 0.87 & 0.85 \\
		& & IV   & 15.07 & 15.23 & 15.70 & 15.58 \\
		& & MISE & 20.78 & 16.70 & 16.57 & 16.43 \\
		
		\cmidrule(lr){2-7}
		
		& \multirow{3}{*}{1000}
		& ISB    & 5.29 & 0.48 & 0.35 & 0.36 \\
		& & IV   & 5.34 & 6.39 & 6.33 & 6.31 \\
		& & MISE & 10.63 & 6.87 & 6.68 & 6.67 \\
		
		\bottomrule
	\end{tabular}
\end{table}

\begin{table}[h!]
	\centering
	\caption{The values of ISB$_1$, IV$_1$ and MISE$_1$, multiplied by $10^3$, for the local linear SBF estimator, the PH estimator, and the proposed bias-corrected estimators.}\label{tab:MISE_first_comp}%
	\small
	\renewcommand{\arraystretch}{0.9}
	\setlength{\tabcolsep}{5pt}
	\begin{tabular}{@{}ccccccc@{}}
		\toprule
		\rule{0pt}{3.2ex}$(\rho_1,\rho_2)$ & $n$ &
		& $\hatm_1^{\rmf{LL}}$ & $\hatm_{1, \rmf{naive}}^{\rmf{PH}}$ & $\hatm_1^{\BR}$ & $\hatm_1^{\rmf{WBR}}$ \\[2pt]
		\midrule
		
		\multirow{6}{*}{(0,3)}
		& \multirow{3}{*}{400}
		& ISB$_1$    & 1.33 & 0.18 & 0.31 & 0.30 \\
		& & IV$_1$   & 3.62 & 3.32 & 4.37 & 4.32 \\
		& & MISE$_1$ & 4.94 & 3.50 & 4.68 & 4.62 \\
		
		\cmidrule(lr){2-7}
		
		& \multirow{3}{*}{1000}
		& ISB$_1$    & 0.49 & 0.09 & 0.10 & 0.10 \\
		& & IV$_1$   & 1.72 & 1.53 & 2.08 & 2.07 \\
		& & MISE$_1$ & 2.21 & 1.62 & 2.18 & 2.16 \\
		
		\cmidrule(lr){1-7}
		
		\multirow{6}{*}{(2,0)}
		& \multirow{3}{*}{400}
		& ISB$_1$    & 3.68 & 0.83 & 0.57 & 0.55 \\
		& & IV$_1$   & 5.75 & 5.80 & 6.11 & 6.04 \\
		& & MISE$_1$ & 9.44 & 6.63 & 6.68 & 6.59 \\
		
		\cmidrule(lr){2-7}
		
		& \multirow{3}{*}{1000}
		& ISB$_1$    & 3.54 & 0.40 & 0.25 & 0.25 \\
		& & IV$_1$   & 2.91 & 3.10 & 3.15 & 3.15 \\
		& & MISE$_1$ & 6.45 & 3.50 & 3.39 & 3.41 \\
		
		\cmidrule(lr){1-7}
		
		\multirow{6}{*}{(2,3)}
		& \multirow{3}{*}{400}
		& ISB$_1$      & 3.80 & 0.69 & 0.59 & 0.58 \\
		& & IV$_1$     & 6.83 & 7.45 & 7.32 & 7.25 \\
		& & MISE$_1$   & 10.62 & 8.14 & 7.91 & 7.83 \\
		
		\cmidrule(lr){2-7}
		
		& \multirow{3}{*}{1000}
		& ISB$_1$    & 3.60 & 0.33 & 0.23 & 0.23 \\
		& & IV$_1$   & 2.99 & 3.26 & 3.33 & 3.31 \\
		& & MISE$_1$ & 6.59 & 3.59 & 3.56 & 3.54 \\
		\bottomrule
	\end{tabular}
\end{table}

\begin{table}[h!]
	\centering
	\caption{The values of ISB$_2$, IV$_2$ and MISE$_2$, multiplied by $10^3$, for the local linear SBF estimator, the PH estimator, and the proposed bias-corrected estimators.}\label{tab:MISE_second_comp}%
	\small
	\renewcommand{\arraystretch}{0.9}
	\setlength{\tabcolsep}{5pt}
	\begin{tabular}{@{}ccccccc@{}}
		\toprule
		\rule{0pt}{3.2ex}$(\rho_1,\rho_2)$ & $n$ &
		& $\hatm_2^{\rmf{LL}}$ & $\hatm_{2, \rmf{naive}}^{\rmf{PH}}$ & $\hatm_2^{\BR}$ & $\hatm_2^{\rmf{WBR}}$ \\[2pt]
		\midrule
		
		\multirow{6}{*}{(0,3)}
		& \multirow{3}{*}{400}
		& ISB$_2$    & 1.66 & 1.47 & 0.28 & 0.28 \\
		& & IV$_2$   & 6.16 & 5.66 & 6.79 & 6.75 \\
		& & MISE$_2$ & 7.81 & 7.13 & 7.07 & 7.03 \\
		
		\cmidrule(lr){2-7}
		
		& \multirow{3}{*}{1000}
		& ISB$_2$    & 1.48 & 1.32 & 0.10 & 0.10 \\
		& & IV$_2$   & 2.35 & 2.14 & 2.70 & 2.69 \\
		& & MISE$_2$ & 3.83 & 3.47 & 2.80 & 2.78 \\
		
		\cmidrule(lr){1-7}
		
		\multirow{6}{*}{(2,0)}
		& \multirow{3}{*}{400}
		& ISB$_2$    & 1.03 & 0.14 & 0.20 & 0.19 \\
		& & IV$_2$   & 4.80 & 4.42 & 5.28 & 5.24 \\
		& & MISE$_2$ & 5.83 & 4.57 & 5.48 & 5.43 \\
		
		\cmidrule(lr){2-7}
		
		& \multirow{3}{*}{1000}
		& ISB$_2$    & 0.50 & 0.10 & 0.09 & 0.08 \\
		& & IV$_2$   & 1.96 & 2.02 & 2.26 & 2.25 \\
		& & MISE$_2$ & 2.47 & 2.12 & 2.36 & 2.33 \\
		
		\cmidrule(lr){1-7}
		
		\multirow{6}{*}{(2,3)}
		& \multirow{3}{*}{400}
		& ISB$_2$    & 1.86 & 0.46 & 0.35 & 0.35 \\
		& & IV$_2$   & 5.93 & 6.62 & 6.66 & 6.62 \\
		& & MISE$_2$ & 7.79 & 7.08 & 7.01 & 6.97 \\
		
		\cmidrule(lr){2-7}
		
		& \multirow{3}{*}{1000}
		& ISB$_2$    & 1.70 & 0.16 & 0.13 & 0.13 \\
		& & IV$_2$   & 2.24 & 2.87 & 2.84 & 2.84 \\
		& & MISE$_2$ & 3.94 & 3.03 & 2.97 & 2.97 \\
		
		\bottomrule
	\end{tabular}
\end{table}
\endgroup

Tables~\ref{tab:MISE_overall}--\ref{tab:MISE_second_comp} report the values of the performance measures for $n\in\{400,1000\}$.
Table~\ref{tab:MISE_overall} shows that, in terms of the estimation accuracy of the additive function $m$, the proposed
$\hatm^{\BR}$ and $\hatm^{\WBR}$ achieve both smaller bias and better MISE performance than the local linear
and parametric-help SBF estimators.
The PH method exhibits better bias and MISE performance than the local linear SBF estimator.
In terms of estimating individual component functions, Tables~\ref{tab:MISE_first_comp}--\ref{tab:MISE_second_comp} demonstrate that the proposed estimators and the PH estimator yield smaller ISB values
than the local linear SBF estimator, demonstrating their effectiveness in reducing bias.

To turn to a more focused comparison of $\hatm_j^{\BR}$ and $\hatm_j^{\WBR}$ with the PH estimator $\hatm_{j, \rmf{naive}}^{\rmf{PH}}$ for $j=1,2$, we first note that the assisting function $g_j^{\rmf{naive}}$ in \eqref{assist-ftn}
coincides with the true additive component $m_j$, up to an additive constant, for $j=1$ when $(\rho_1,\rho_2)=(0,3)$ and for $j=2$ when $(\rho_1,\rho_2)=(2,0)$.
The leading biases of $\hatm_{j,\rmf{naive}}^\PH$ vanish in the respective cases, and are therefore expected to exhibit bias performance  comparable to those of $\hatm_j^{\BR}$ and $\hatm_j^{\WBR}$ in these cases.
Indeed, the $\mathrm{ISB}_1$ values of $\hatm_{1,\rmf{naive}}^\PH$ are comparable to those of $\hatm_1^{\BR}$ and $\hatm_1^{\WBR}$ when
$(\rho_1,\rho_2)=(0,3)$, while the $\mathrm{ISB}_2$ values of $\hatm_{2,\rmf{naive}}^\PH$ are comparable to those of $\hatm_2^{\BR}$ and $\hatm_2^{\WBR}$ when $(\rho_1,\rho_2)=(2,0)$.
Overall, $\hatm_j^{\BR}$ and $\hatm_j^{\WBR}$ exhibit better bias and MISE performance than $\hatm_{j,\rmf{naive}}^\PH$, except the case
where $\hatm_{j,\rmf{naive}}^\PH$ uses the respective \textit{oracle} assisting functions $g_j^{\rmf{naive}}$.

\subsection{Statistical inference}

We now conduct simulation studies to evaluate the practical performance of the proposed statistical inference introduced in Section~\ref{sec:inference}.
The simulation data are generated according to \eqref{eq:sim_model} with $(\rho_1,\rho_2)=(0,0)$, that is,
$Y = m_1(X_1) + m_2(X_2)+ \sigma(X_1,X_2) \, \varepsilon$ with $m(x_1)=\sin(2\pi x_1)$ and $m_2(x_2)=\cos(2\pi x_2)$.
We consider two conditional variance structures:
\begin{align}\label{var-str}
\sigma_1(x_1, x_2;\alpha) = C_\sigma, \quad \sigma_2(x_1, x_2;\alpha) = C_\sigma \cdot \big\{2.5 +0.4 \sin(6 x_1) - 2x_2 \big\}.
\end{align}
For the auxiliary bandwidths, we use $b_{j,\ell} = (1+ \ell / 20)h_{j,\rmf{opt}}\,n^{1/5}/2$ for $\ell=1,\dots,10$, where $h_{j,\rmf{opt}}$ are the optimal bandwidths obtained via the CBS algorithm introduced in Section~\ref{subsec:bias_reduction_sim}.
The simultaneous confidence bands (SCBs) for each additive component are constructed over the interval $I = [0.1, 0.9]$.

We compare the proposed testing procedure with an oracle approach.
Specifically, we define oracle responses as
\begin{align*}
	Y^{\rmf{oracle},1} := Y - m_2(X_2) = m_1(X_1) + \sigma(X_1,X_2) \, \varepsilon, \\
	Y^{\rmf{oracle},2} := Y - m_1(X_1) = m_2(X_2) + \sigma(X_1,X_2) \, \varepsilon.
\end{align*}
Using these oracle responses, we construct an SCB for $m_j$ over $I$ by applying the method of \cite{cheng2018bias}, denoted as SCB$_\rmf{oracle}$.

\begin{figure}[h]
	\centering
	
	\begin{minipage}{0.48\linewidth}
		\centering
		\includegraphics[width=\linewidth]{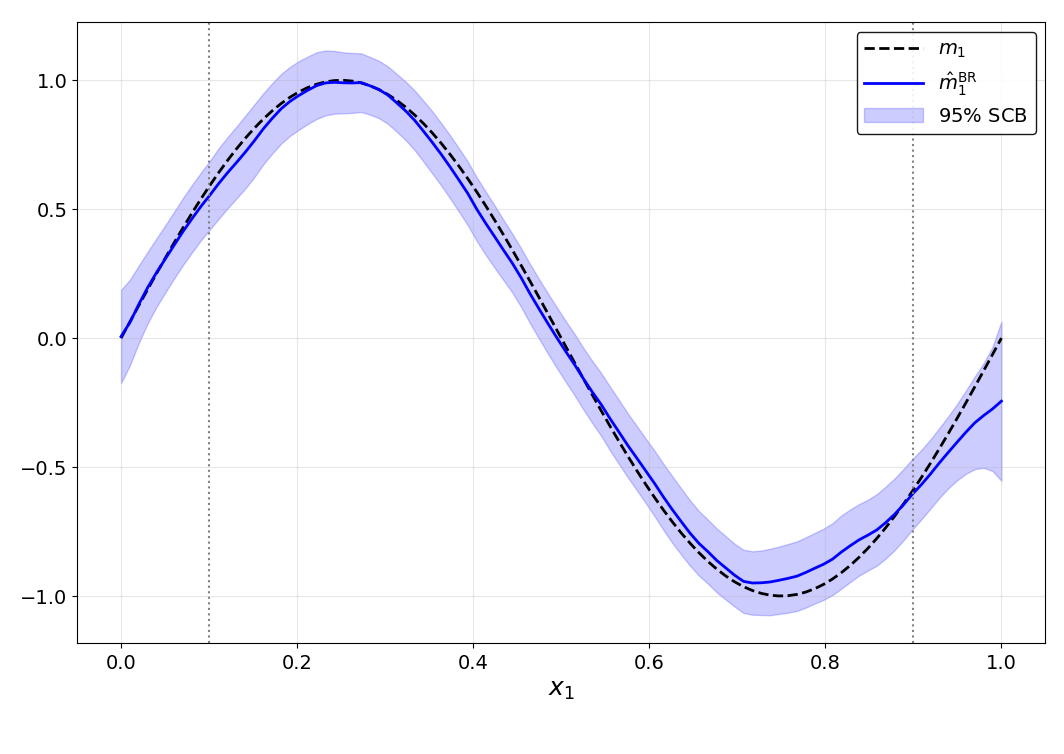}
		{(a) $\hatm_1$ under $\sigma_1$ with $C_\sigma=0.4$}
	\end{minipage}
	\hfill
	\begin{minipage}{0.48\linewidth}
		\centering
		\includegraphics[width=\linewidth]{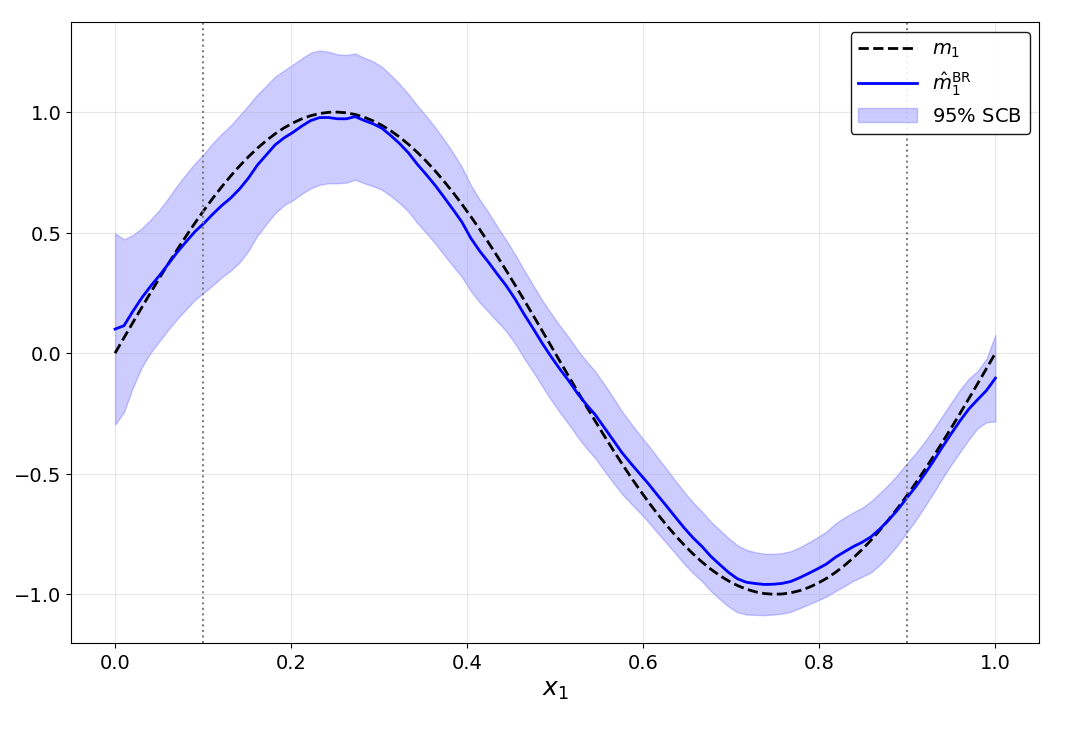}
		{(b) $\hatm_1$ under $\sigma_2$ with $C_\sigma=0.4$}
	\end{minipage}
	
	\vspace{0.5cm}
	
	\begin{minipage}{0.48\linewidth}
		\centering
		\includegraphics[width=\linewidth]{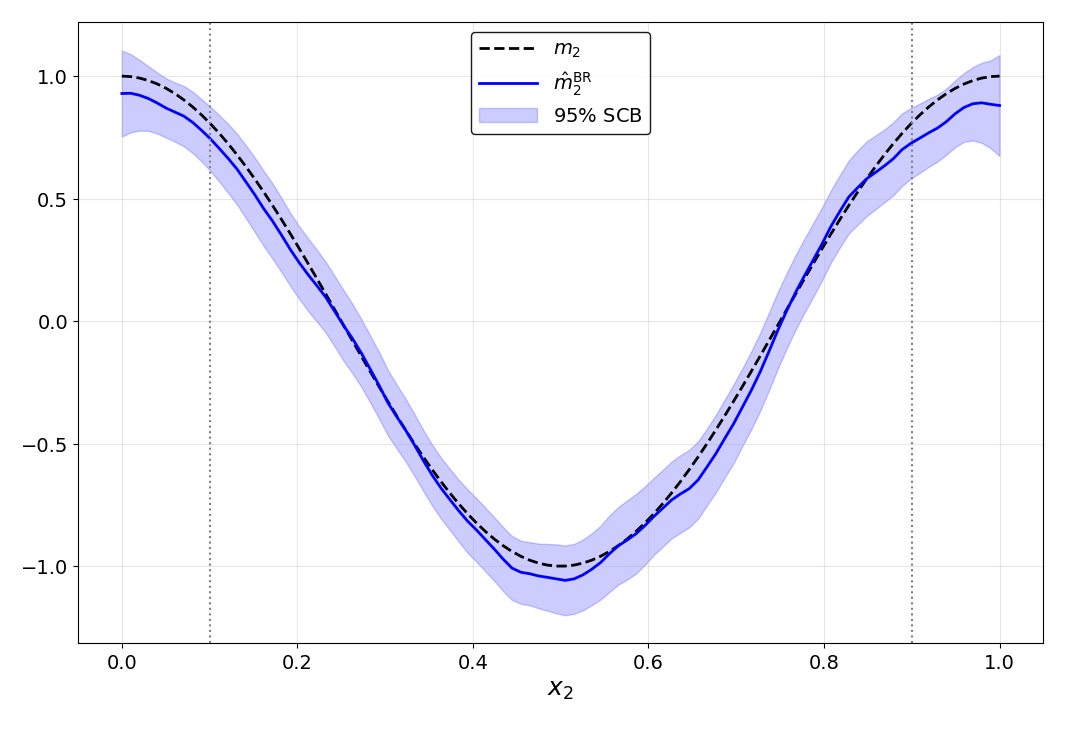}
		{(c) $\hatm_2$ under $\sigma_1$ with $C_\sigma=0.4$}
	\end{minipage}
	\hfill
	\begin{minipage}{0.48\linewidth}
		\centering
		\includegraphics[width=\linewidth]{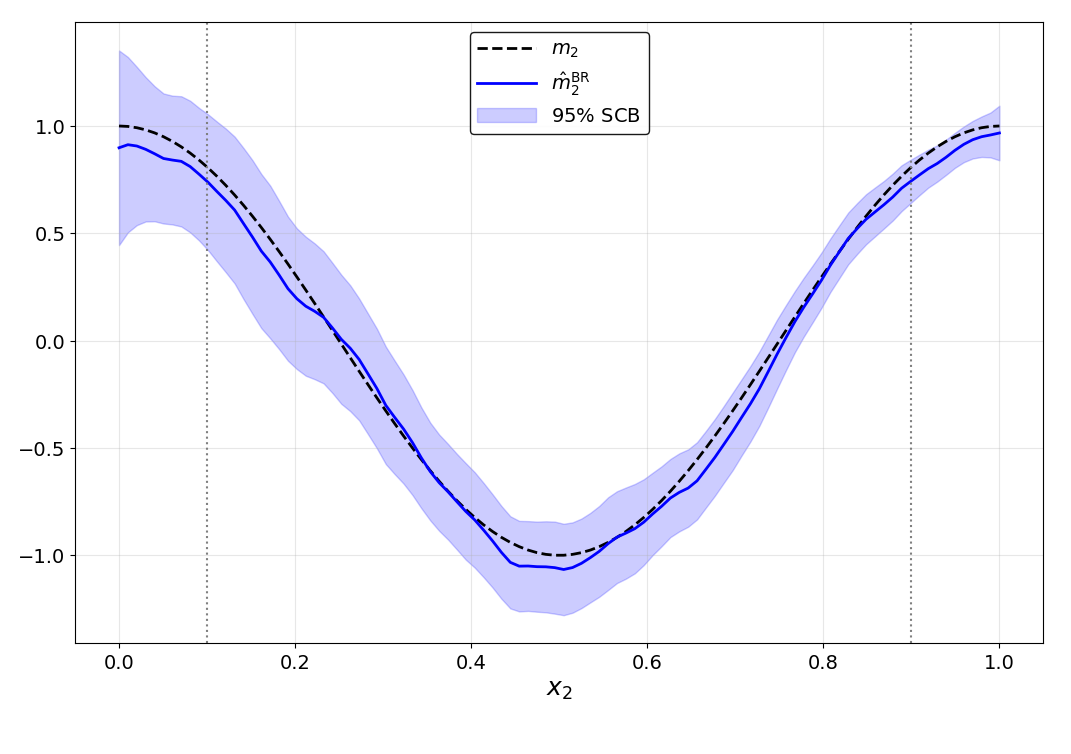}
		{(d) $\hatm_2$ under $\sigma_2$ with $C_\sigma=0.4$}
	\end{minipage}

	\caption{Simultaneous confidence bands for $m_j$ under the conditional variance structures $\sigma_1$ and $\sigma_2$ in
\eqref{var-str} with $C_\sigma=0.4$.}\label{fig:SCB_const_hetero}
\end{figure}

\begin{table}[h!]
	\centering
	\caption{Coverage rates of the 95\% simultaneous confidence bands (SCBs)
		for the component functions $m_1$ and $m_2$
		under the conditional variance structure $\sigma_1$ in \eqref{var-str} with different variance levels $C_\sigma$
		and sample sizes $n$.
	}\label{tab:SCB_sigma1}%
	\renewcommand{\arraystretch}{0.95}
	\small
	\begin{tabular}{@{}cccccc@{}}
		\toprule
		\multirow{2}{*}{$n$} & \multirow{2}{*}{$C_\sigma$} & \multicolumn{2}{c}{$m_1$} & \multicolumn{2}{c}{$m_2$} \\ \cmidrule(lr){3-4}\cmidrule(lr){5-6}
		& & SCB$_{\rmf{oracle}}$ (\%) & SCB (\%) & SCB$_{\rmf{oracle}}$ (\%) & SCB (\%) \\ \midrule
		\multirow{2}{*}{1000}
		& 0.4 & 92 & 89 & 97 & 91\\
		& 0.7 & 96 & 95 & 95 & 94\\
		\midrule
		\multirow{2}{*}{1500}
		& 0.4 & 96 & 92 & 93 & 92\\
		& 0.7 & 91 & 96 & 93 & 92\\
		\bottomrule
	\end{tabular}
\end{table}

\begin{table}[h!]
	\centering
	\caption{Coverage rates of the 95\% simultaneous confidence bands (SCBs)
		for the component functions $m_1$ and $m_2$
		under the conditional variance structure $\sigma_2$ in \eqref{var-str} with different variance levels $C_\sigma$
		and sample sizes $n$.
	}\label{tab:SCB_sigma2}%
	\renewcommand{\arraystretch}{0.95}
	\small
	\begin{tabular}{@{}cccccc@{}}
		\toprule
		\multirow{2}{*}{$n$} & \multirow{2}{*}{$C_\sigma$} & \multicolumn{2}{c}{$m_1$} & \multicolumn{2}{c}{$m_2$} \\ \cmidrule(lr){3-4}\cmidrule(lr){5-6}
		& & SCB$_{\rmf{oracle}}$ (\%) & SCB (\%) & SCB$_{\rmf{oracle}}$ (\%) & SCB (\%) \\ \midrule
		\multirow{2}{*}{1000}
		& 0.4 & 91 & 90 & 92 & 91\\
		& 0.7 & 95 & 93 & 94 & 91\\
		\midrule
		\multirow{2}{*}{1500}
		& 0.4 & 95 & 92 & 95 & 93 \\
		& 0.7 & 94 & 91 & 97 & 93 \\
		\bottomrule
	\end{tabular}
\end{table}

Figure~\ref{fig:SCB_const_hetero} presents the bias-corrected estimators $\hatm_j^{\BR}$
and their corresponding $95\%$ SCBs under two conditional variance structures $\sigma_1$ and $\sigma_2$ in \eqref{var-str} with $C_\sigma=0.4$.
Under $\sigma_1$, the SCB width remains approximately constant over the domain, whereas under $\sigma_2$, it varies with $x_2$.
Tables~\ref{tab:SCB_sigma1} and \ref{tab:SCB_sigma2} report the empirical coverage rates of the proposed and oracle SCBs,
with a nominal coverage level of $95\%$.
The SCBs constructed using $\hatm^{\BR}$ perform comparably to the oracle SCBs based on univariate estimation,
with empirical coverage rates close to the nominal level of $95\%$.

\section{Real data applications}\label{sec:real_data}

\subsection{Bike rental data}\label{subsec:bike_rental}

Predicting daily bicycle rental demand is an important task in urban mobility planning, as accurate forecasts facilitate efficient bike allocation and help prevent shortages.
Previous studies have documented that meteorological conditions play a significant role in predicting bicycle usage \citep{faghih2014land,el2017effects}.
We apply the proposed additive regression method to the Bikeshare dataset from the R package \texttt{ISLR2}, which contains observations collected from the Capital Bikeshare system in Washington, D.C., during 2011--2012.
In this analysis, the response variable is the daily bicycle rental count and the covariates are daily average temperature and daily average humidity.
After removing outliers, we use $n=941$ observations.
To improve predictive performance, we apply logarithmic transformations to both the response and the covariates.
Also, the covariates are linearly transformed to lie in the unit interval $[0,1]$.

\begin{table}[h]
	\centering
	\caption{Comparison of 5-fold cross-validated ASPE values in the analysis of the Bikeshare dataset.}\label{tab:ASPE_real_1}%
	\renewcommand{\arraystretch}{0.95}
	\small
	\begin{tabular}{@{}lcccc@{}}
		\toprule
		Estimator & $\hatm^{\rmf{LL}}$ & $\hatm^{\rmf{PH}}$	& $\hatm^{\BR}$ & $\hatm^{\rmf{WBR}}$ \\ \midrule
		ASPE & 1.7630 & 1.7628 & 1.7453 & 1.7451 \\ \bottomrule
	\end{tabular}
\end{table}

We compare our method with the local linear SBF estimator $\hatm^{\rmf{LL}}$ and the parametric-help SBF estimator $\hatm^{\rmf{PH}}$.
For the proposed bias reduction methods, the auxiliary bandwidths are chosen as $b_{j,\ell} = (1 + C_\ell / 20)\, h_{j,\rmf{opt}}$, where $C_\ell = 1,2,\ldots,7$.
Also, for the implementation of the PH estimator, we employ $g_1(x_1)=\sin(2\pi x_1)$ and $g_2(x_2)=\cos(2\pi x_2)$ as the assisting functions.
To assess predictive performance, we adopt 5-fold cross-validated average squared prediction error (ASPE) defined by
\begin{equation*}
	\rmf{ASPE} = \frac{1}{n}\sum_{k=1}^{5}\sum_{i \in I_k} \big( Y_i - \hat{Y}_i^{(-I_k)} \big)^2,
\end{equation*}
where $I_k$ denotes the index set of the $k$-th fold, and $\hat{Y}_i^{(-I_k)}$ denotes the prediction of $Y_i$ obtained from the model trained without the observations in $I_k$.
Table~\ref{tab:ASPE_real_1} reports the 5-fold cross-validated ASPE values for the competing estimators.
The proposed bias-corrected estimators $\hatm^{\BR}$ and $\hatm^{\rmf{WBR}}$ yield lower ASPE values than $\hatm^{\rmf{LL}}$ and $\hatm^{\rmf{PH}}$, suggesting that the proposed bias reduction also leads to improved predictive performance.

\begin{figure}[h]
	\centering
	
	\begin{minipage}{0.6\linewidth}
		\centering
		\includegraphics[width=\linewidth]{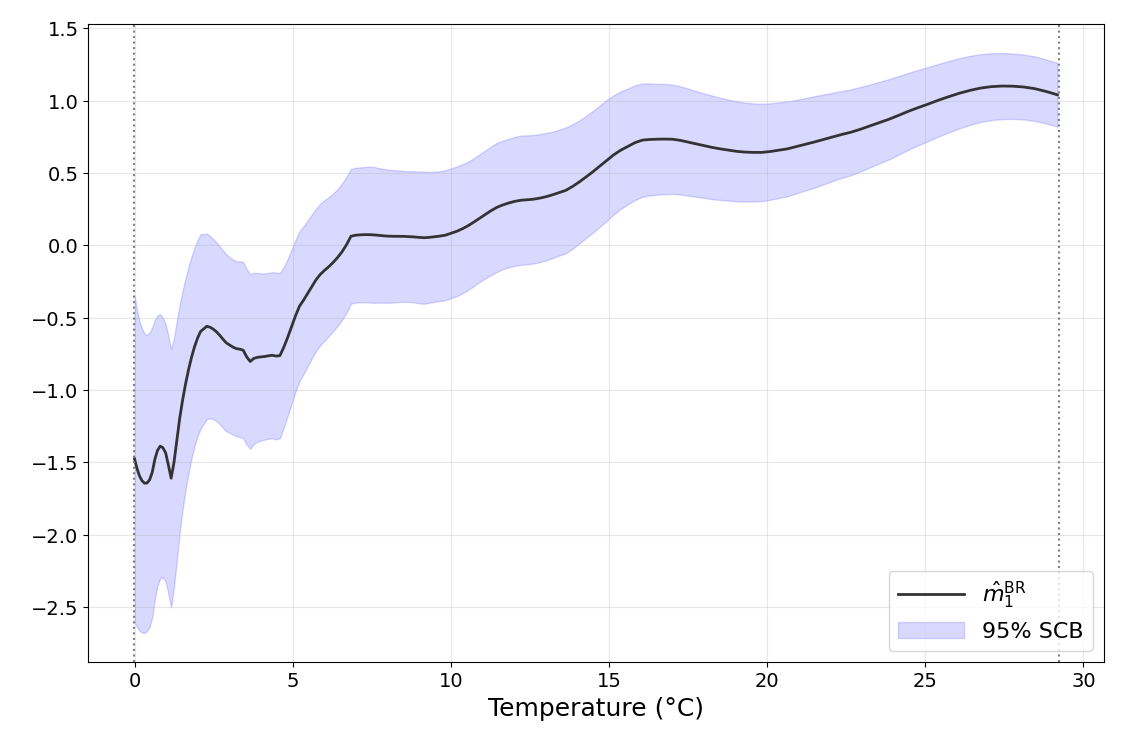}
		{(a) The effect of temperature}
	\end{minipage}
	\hfill
	\begin{minipage}{0.6\linewidth}
		\centering
		\includegraphics[width=\linewidth]{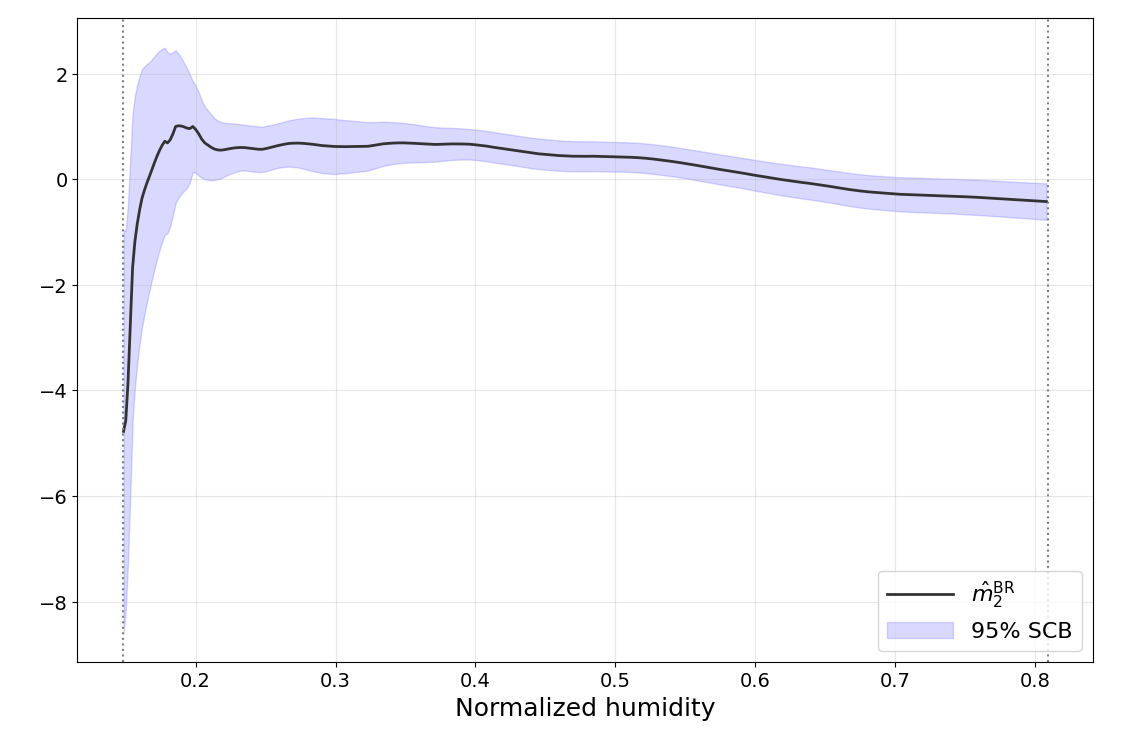}
		{(b) The effect of humidity}
	\end{minipage}
	
	\caption{Estimated additive component functions $\hatm_j^\BR$ in the analysis of the Bikeshare dataset.
	The 95\% simultaneous confidence bands (SCBs) for the panels are constructed.
	}
	\label{fig:real_1}%
\end{figure}

Figure~\ref{fig:real_1} depicts the estimated additive component functions $\hatm_j^\BR$ and their 95$\%$ simultaneous confidence bands on $I = [0.1, 0.9]$.
The top panel shows a sharp decrease in bike rental counts at low temperatures, which is consistent with our expectation that people tend to rent fewer bikes in cold weather.
For humidity, rental counts decrease at both low and high humidity levels, with a more pronounced decrease at low humidity levels.
Furthermore, we test the null hypothesis
\begin{equation}\label{eq:real_data_hypothesis}
	H_0 : m_j(x_j) = 0  \text{ for all } x_j \in I \quad \text{vs} \quad  H_1 : \text{not } H_0
\end{equation}
using the procedure proposed in Section~\ref{subsec:inference}.
The test statistics $T_{nj}(I)$ defined at \eqref{eq:test_statistic} are 14.82 and 9.24 for the temperature and humidity, respectively.
Since the test statistics exceed the critical value $C_\alpha=3.66$ for $\alpha=0.05$ defined at \eqref{eq:reject_test}, the null hypothesis \eqref{eq:real_data_hypothesis} is rejected for both components.
These results indicate that both temperature and humidity have statistically significant effects on bike rental counts at the 5\% significance level.

\subsection{Flights data}

Since excessive flight delays can lead to congestion and operational inefficiencies, the accurate prediction of flight departure delays is vital for airport scheduling and air traffic management.
Weather conditions significantly impact these delays, as adverse meteorological situations such as extreme temperatures or high wind speed can disrupt airport operations.
We analyze the flights dataset in the R package \texttt{nycflights13} augmented with the NYC Weather Dataset from \url{https://www.kaggle.com/datasets/muhammedjunayed/nyc-weather-2009-2015}.
In this analysis, we restrict attention to flights departing from John F. Kennedy International Airport in 2013.
The response variable is the log-transformed delay time, while the covariates are temperature and wind speed, which describe atmospheric conditions that possibly influence aircraft operations and airport ground handling.
After removing missing values, we use a final sample of $n=978$ observations to fit the additive regression model.

\begin{table}[h!]
	\centering
	\caption{Comparison of 5-fold cross-validated ASPE values in the analysis of the flights dataset.}\label{tab:ASPE_real_2}%
	\renewcommand{\arraystretch}{0.95}
	\small
	\begin{tabular}{@{}lcccc@{}}
		\toprule
		Estimator & $\hatm^{\rmf{LL}}$ & $\hatm^{\rmf{PH}}$ & $\hatm^{\BR}$ & $\hatm^{\rmf{WBR}}$ \\ \midrule
		ASPE & 1.5500 & 1.5488 & 1.5487 & 1.5478 \\ \bottomrule
	\end{tabular}
\end{table}

We again compare the proposed estimators with the local linear and parametric-help SBF estimators.
We use the same auxiliary bandwidths and the prespecified functions as those used in Section~\ref{subsec:bike_rental}.
Table~\ref{tab:ASPE_real_2} reports the 5-fold cross-validated ASPE values for the competing estimators.
Consistent with the analysis of the Bikeshare dataset, the ASPE values for the proposed bias-corrected estimators are lower than or comparable to those of the competing estimators.

\begin{figure}[h]
	\centering
	
	\begin{minipage}{0.65\linewidth}
		\centering
		\includegraphics[width=\linewidth]{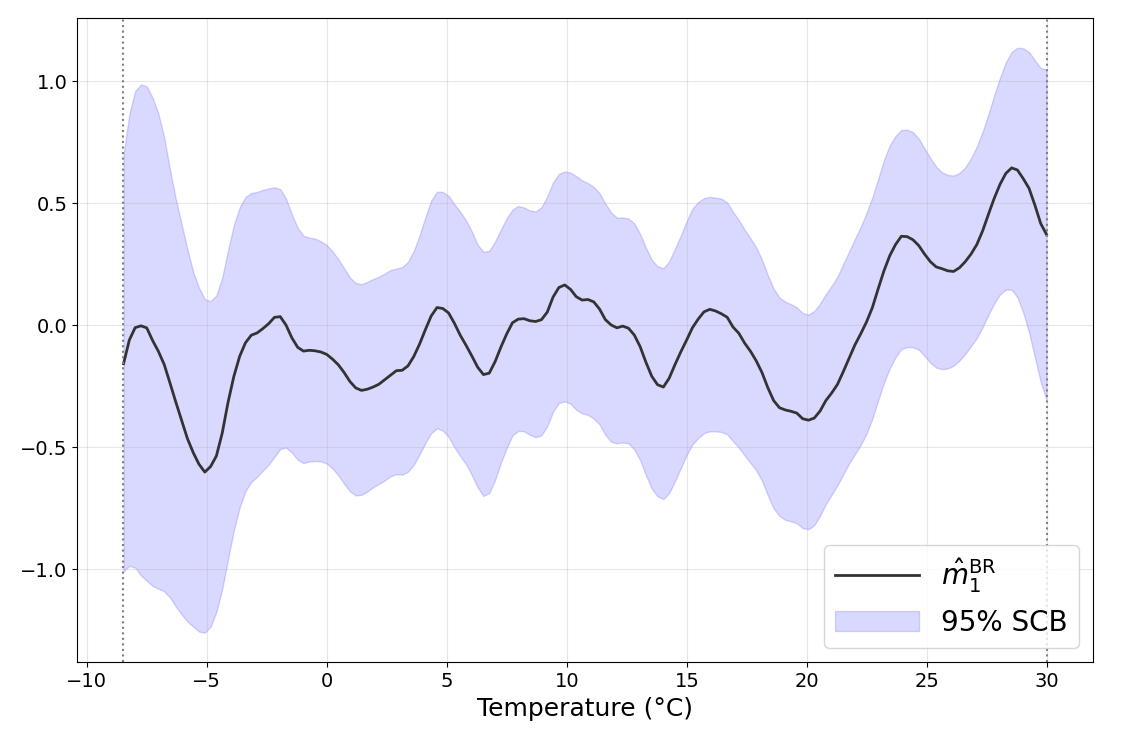}
		{(a) The effect of temperature}
	\end{minipage}
	\hfill
	\begin{minipage}{0.65\linewidth}
		\centering
		\includegraphics[width=\linewidth]{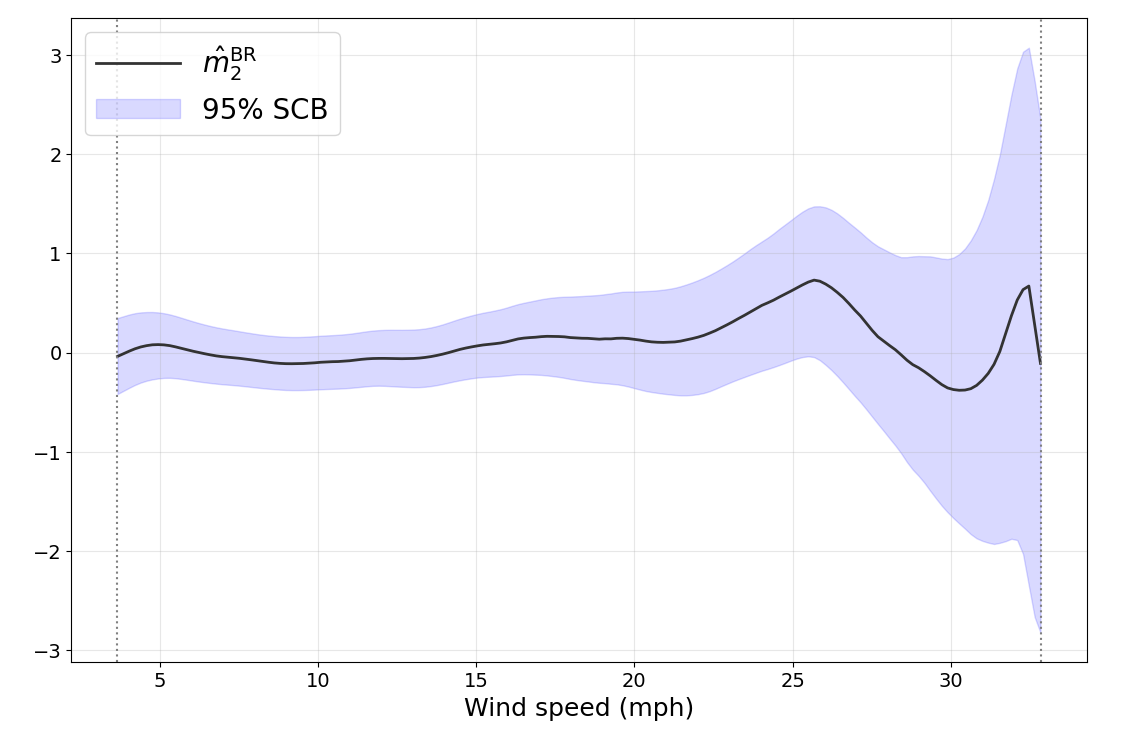}
		{(b) The effect of wind speed}
	\end{minipage}
	
	\caption{Estimated additive component functions $\hatm_j^\BR$ in the analysis of the flights dataset.
	The 95\% simultaneous confidence bands (SCBs) for the panels are constructed.
	}
	\label{fig:real_2}%
\end{figure}

Figure~\ref{fig:real_2} displays the estimated additive component functions $\hatm_j^{\BR}$, with the shaded region representing the simultaneous confidence bands over $I=[0.1,0.9]$.
The left panel illustrates that temperature has a gradual positive effect on flight delays.
Similarly, wind speed positively impacts the delay time, but it exhibits substantially larger variability at high values.
Furthermore, we again test the null hypothesis \eqref{eq:real_data_hypothesis}.
The test statistics $T_{nj}(I)$ defined at \eqref{eq:test_statistic} are 5.74 and 3.48 for temperature and wind speed, respectively.
Comparing these to the critical value $C_\alpha=3.66$ for $\alpha=0.05$, the null hypothesis is rejected for temperature but is not rejected for wind speed.
These results provide evidence that temperature has a significant effect on flight delays at the 5\% significance level, whereas the effect of wind speed is not statistically significant.

\section{Conclusion}\label{sec:conclusion}

In this paper, we develop a framework for simultaneous inference on additive component functions.
To this end, we propose a least-squares-based bias reduction method that incorporates a least-squares adjustment
into the local linear smooth backfitting framework. The resulting bias-reduced estimators eliminate the leading asymptotic bias.
Building on these estimators, we establish the limiting distributions of the uniform estimation errors
of the additive component estimators and develop a simultaneous inference procedure based on these results.
Simulation studies demonstrate the favorable finite-sample performance of the proposed method in terms of bias reduction
and coverage accuracy relative to existing approaches. Furthermore, applications to real data illustrate the practical utility
of the proposed procedure for capturing nonlinear relationships while providing reliable uncertainty quantification.
Overall, our approach provides an effective and theoretically grounded framework for inference on component functions in additive models.

Several avenues for future research remain. First, the proposed framework for simultaneous inference could be extended
to other structured nonparametric regression models, including varying coefficient models \citep{moon2025varying}
and partially linear additive models \citep{cho2024partially,choi2025high}.
Second, it would be worthwhile to develop methods for constructing simultaneous confidence bands for generalized additive models, where nonlinear link functions introduce additional technical challenges \citep{liu2013oracally}.
Finally, the proposed methodology could be extended to non-Euclidean responses in the context of Fréchet regression \citep{petersen2019frechet,chen2022uniform}.

\backmatter


\section*{Declarations}

\bmhead{Funding}
Changwon Choi’s research was supported by Basic Science Research Program through the National Research Foundation of Korea (NRF) funded by the Ministry of Education (No. RS-2024-00446677).
Byeong U. Park’s research was supported by the National Research Foundation (NRF) grant funded by Korea government (MSIT) (No. RS-2026-25606859).


\bibliography{bibliography}

\begin{thebibliography}{}
\renewcommand{\doi}[1]{\url{https://doi.org/#1}}
\bibcommenthead

\bibitem [\protect \citeauthoryear {%
Adler%
\ \BBA {} Taylor%
}{%
Adler%
\ \BBA {} Taylor%
}{%
{\protect \APACyear {2007}}%
}]{%
adler2007random}
\APACinsertmetastar {%
adler2007random}%
\begin{APACrefauthors}%
Adler, R.J.%
\BCBT {}\ \BBA {} Taylor, J.E.%
\end{APACrefauthors}%
\unskip\
\newblock
\APACrefYear{2007}.
\newblock
\APACrefbtitle {Random fields and geometry} {Random fields and geometry}.
\newblock
\APACaddressPublisher{}{Springer}.
\PrintBackRefs{\CurrentBib}

\bibitem [\protect \citeauthoryear {%
S.X.~Chen%
\ \BBA {} Qin%
}{%
S.X.~Chen%
\ \BBA {} Qin%
}{%
{\protect \APACyear {2002}}%
}]{%
chen2002confidence}
\APACinsertmetastar {%
chen2002confidence}%
\begin{APACrefauthors}%
Chen, S.X.%
\BCBT {}\ \BBA {} Qin, Y.S.%
\end{APACrefauthors}%
\unskip\
\newblock
\APACrefYearMonthDay{2002}{}{}.
\newblock
{\BBOQ}\APACrefatitle {Confidence intervals based on local linear smoother}
  {Confidence intervals based on local linear smoother}.{\BBCQ}
\newblock
\APACjournalVolNumPages{Scandinavian Journal of Statistics}{29}{1}{89--99,}
\newblock
\begin{APACrefDOI} \doi{10.1111/1467-9469.00273} \end{APACrefDOI}
\newblock

\newblock

\PrintBackRefs{\CurrentBib}

\bibitem [\protect \citeauthoryear {%
Y.~Chen%
\ \BBA {} M\"{u}ller%
}{%
Y.~Chen%
\ \BBA {} M\"{u}ller%
}{%
{\protect \APACyear {2022}}%
}]{%
chen2022uniform}
\APACinsertmetastar {%
chen2022uniform}%
\begin{APACrefauthors}%
Chen, Y.%
\BCBT {}\ \BBA {} M\"{u}ller, H\BHBI G.%
\end{APACrefauthors}%
\unskip\
\newblock
\APACrefYearMonthDay{2022}{}{}.
\newblock
{\BBOQ}\APACrefatitle {Uniform convergence of local {F}r\'{e}chet regression
  with applications to locating extrema and time warping for metric space
  valued trajectories} {Uniform convergence of local {F}r\'{e}chet regression
  with applications to locating extrema and time warping for metric space
  valued trajectories}.{\BBCQ}
\newblock
\APACjournalVolNumPages{The Annals of Statistics}{50}{3}{1573--1592,}
\newblock
\begin{APACrefDOI} \doi{10.1214/21-aos2163} \end{APACrefDOI}
\newblock

\newblock

\PrintBackRefs{\CurrentBib}

\bibitem [\protect \citeauthoryear {%
G.~Cheng%
\ \BBA {} Chen%
}{%
G.~Cheng%
\ \BBA {} Chen%
}{%
{\protect \APACyear {2019}}%
}]{%
cheng2019nonparametric}
\APACinsertmetastar {%
cheng2019nonparametric}%
\begin{APACrefauthors}%
Cheng, G.%
\BCBT {}\ \BBA {} Chen, Y\BHBI C.%
\end{APACrefauthors}%
\unskip\
\newblock
\APACrefYearMonthDay{2019}{}{}.
\newblock
{\BBOQ}\APACrefatitle {Nonparametric inference via bootstrapping the debiased
  estimator} {Nonparametric inference via bootstrapping the debiased
  estimator}.{\BBCQ}
\newblock
\APACjournalVolNumPages{Electronic Journal of Statistics}{13}{1}{2194--2256,}
\newblock
\begin{APACrefDOI} \doi{10.1214/19-EJS1575} \end{APACrefDOI}
\newblock

\newblock

\PrintBackRefs{\CurrentBib}

\bibitem [\protect \citeauthoryear {%
M\BHBI Y.~Cheng%
, Huang%
, Liu%
\BCBL {}\ \BBA {} Peng%
}{%
M\BHBI Y.~Cheng%
\ \protect \BOthers {.}}{%
{\protect \APACyear {2018}}%
}]{%
cheng2018bias}
\APACinsertmetastar {%
cheng2018bias}%
\begin{APACrefauthors}%
Cheng, M\BHBI Y.%
, Huang, T.%
, Liu, P.%
\BCBL {} Peng, H.%
\end{APACrefauthors}%
\unskip\
\newblock
\APACrefYearMonthDay{2018}{}{}.
\newblock
{\BBOQ}\APACrefatitle {Bias reduction for nonparametric and semiparametric
  regression models} {Bias reduction for nonparametric and semiparametric
  regression models}.{\BBCQ}
\newblock
\APACjournalVolNumPages{Statistica Sinica}{28}{5}{2749--2770,}
\newblock
\begin{APACrefDOI} \doi{10.5705/ss.202017.0058} \end{APACrefDOI}
\newblock

\newblock

\PrintBackRefs{\CurrentBib}

\bibitem [\protect \citeauthoryear {%
Chernozhukov%
, Chetverikov%
\BCBL {}\ \BBA {} Kato%
}{%
Chernozhukov%
\ \protect \BOthers {.}}{%
{\protect \APACyear {2014}}%
{\protect \APACexlab {{\protect \BCnt {1}}}}}]{%
chernozhukov2014anti}
\APACinsertmetastar {%
chernozhukov2014anti}%
\begin{APACrefauthors}%
Chernozhukov, V.%
, Chetverikov, D.%
\BCBL {} Kato, K.%
\end{APACrefauthors}%
\unskip\
\newblock
\APACrefYearMonthDay{2014{\protect \BCnt {1}}}{}{}.
\newblock
{\BBOQ}\APACrefatitle {Anti-concentration and honest, adaptive confidence
  bands} {Anti-concentration and honest, adaptive confidence bands}.{\BBCQ}
\newblock
\APACjournalVolNumPages{The Annals of Statistics}{42}{5}{1787--1818,}
\newblock
\begin{APACrefDOI} \doi{10.1214/14-AOS1235} \end{APACrefDOI}
\newblock

\newblock

\PrintBackRefs{\CurrentBib}

\bibitem [\protect \citeauthoryear {%
Chernozhukov%
, Chetverikov%
\BCBL {}\ \BBA {} Kato%
}{%
Chernozhukov%
\ \protect \BOthers {.}}{%
{\protect \APACyear {2014}}%
{\protect \APACexlab {{\protect \BCnt {2}}}}}]{%
chernozhukov2014gaussian}
\APACinsertmetastar {%
chernozhukov2014gaussian}%
\begin{APACrefauthors}%
Chernozhukov, V.%
, Chetverikov, D.%
\BCBL {} Kato, K.%
\end{APACrefauthors}%
\unskip\
\newblock
\APACrefYearMonthDay{2014{\protect \BCnt {2}}}{}{}.
\newblock
{\BBOQ}\APACrefatitle {Gaussian approximation of suprema of empirical
  processes} {Gaussian approximation of suprema of empirical processes}.{\BBCQ}
\newblock
\APACjournalVolNumPages{The Annals of Statistics}{42}{4}{1564--1597,}
\newblock
\begin{APACrefDOI} \doi{10.1214/14-AOS1230} \end{APACrefDOI}
\newblock

\newblock

\PrintBackRefs{\CurrentBib}

\bibitem [\protect \citeauthoryear {%
Chernozhukov%
, Chetverikov%
\BCBL {}\ \BBA {} Kato%
}{%
Chernozhukov%
\ \protect \BOthers {.}}{%
{\protect \APACyear {2015}}%
}]{%
chernozhukov2015comparison}
\APACinsertmetastar {%
chernozhukov2015comparison}%
\begin{APACrefauthors}%
Chernozhukov, V.%
, Chetverikov, D.%
\BCBL {} Kato, K.%
\end{APACrefauthors}%
\unskip\
\newblock
\APACrefYearMonthDay{2015}{}{}.
\newblock
{\BBOQ}\APACrefatitle {Comparison and anti-concentration bounds for maxima of
  {G}aussian random vectors} {Comparison and anti-concentration bounds for
  maxima of {G}aussian random vectors}.{\BBCQ}
\newblock
\APACjournalVolNumPages{Probability Theory and Related
  Fields}{162}{1-2}{47--70,}
\newblock
\begin{APACrefDOI} \doi{10.1007/s00440-014-0565-9} \end{APACrefDOI}
\newblock

\newblock

\PrintBackRefs{\CurrentBib}

\bibitem [\protect \citeauthoryear {%
Cho%
, Jeon%
, Kim%
, Yu%
\BCBL {}\ \BBA {} Park%
}{%
Cho%
\ \protect \BOthers {.}}{%
{\protect \APACyear {2024}}%
}]{%
cho2024partially}
\APACinsertmetastar {%
cho2024partially}%
\begin{APACrefauthors}%
Cho, S.%
, Jeon, J.M.%
, Kim, D.%
, Yu, K.%
\BCBL {} Park, B.U.%
\end{APACrefauthors}%
\unskip\
\newblock
\APACrefYearMonthDay{2024}{}{}.
\newblock
{\BBOQ}\APACrefatitle {Partially linear additive regression with a general
  {H}ilbertian response} {Partially linear additive regression with a general
  {H}ilbertian response}.{\BBCQ}
\newblock
\APACjournalVolNumPages{Journal of the American Statistical
  Association}{119}{546}{942--956,}
\newblock
\begin{APACrefDOI} \doi{10.1080/01621459.2022.2149407} \end{APACrefDOI}
\newblock

\newblock

\PrintBackRefs{\CurrentBib}

\bibitem [\protect \citeauthoryear {%
C.~Choi%
, Lin%
\BCBL {}\ \BBA {} Park%
}{%
C.~Choi%
\ \protect \BOthers {.}}{%
{\protect \APACyear {2025}}%
}]{%
choi2025high}
\APACinsertmetastar {%
choi2025high}%
\begin{APACrefauthors}%
Choi, C.%
, Lin, Z.%
\BCBL {} Park, B.U.%
\end{APACrefauthors}%
\unskip\
\newblock
\APACrefYearMonthDay{2025}{}{}.
\newblock
{\BBOQ}\APACrefatitle {High-dimensional partially linear additive models on
  {R}iemannian manifolds} {High-dimensional partially linear additive models on
  {R}iemannian manifolds}.{\BBCQ}
\newblock
\APACjournalVolNumPages{Bernoulli}{31}{4}{3285--3308,}
\newblock
\begin{APACrefDOI} \doi{10.3150/24-BEJ1847} \end{APACrefDOI}
\newblock

\newblock

\PrintBackRefs{\CurrentBib}

\bibitem [\protect \citeauthoryear {%
E.~Choi%
\ \BBA {} Hall%
}{%
E.~Choi%
\ \BBA {} Hall%
}{%
{\protect \APACyear {1998}}%
}]{%
choi1998bias}
\APACinsertmetastar {%
choi1998bias}%
\begin{APACrefauthors}%
Choi, E.%
\BCBT {}\ \BBA {} Hall, P.%
\end{APACrefauthors}%
\unskip\
\newblock
\APACrefYearMonthDay{1998}{}{}.
\newblock
{\BBOQ}\APACrefatitle {On bias reduction in local linear smoothing} {On bias
  reduction in local linear smoothing}.{\BBCQ}
\newblock
\APACjournalVolNumPages{Biometrika}{85}{2}{333--345,}
\newblock
\begin{APACrefDOI} \doi{10.1093/biomet/85.2.333} \end{APACrefDOI}
\newblock

\newblock

\PrintBackRefs{\CurrentBib}

\bibitem [\protect \citeauthoryear {%
E.~Choi%
, Hall%
\BCBL {}\ \BBA {} Rousson%
}{%
E.~Choi%
\ \protect \BOthers {.}}{%
{\protect \APACyear {2000}}%
}]{%
choi2000data}
\APACinsertmetastar {%
choi2000data}%
\begin{APACrefauthors}%
Choi, E.%
, Hall, P.%
\BCBL {} Rousson, V.%
\end{APACrefauthors}%
\unskip\
\newblock
\APACrefYearMonthDay{2000}{}{}.
\newblock
{\BBOQ}\APACrefatitle {Data sharpening methods for bias reduction in
  nonparametric regression} {Data sharpening methods for bias reduction in
  nonparametric regression}.{\BBCQ}
\newblock
\APACjournalVolNumPages{The Annals of Statistics}{28}{5}{1339--1355,}
\newblock
\begin{APACrefDOI} \doi{10.1214/aos/1015957396} \end{APACrefDOI}
\newblock

\newblock

\PrintBackRefs{\CurrentBib}

\bibitem [\protect \citeauthoryear {%
Claeskens%
, Krivobokova%
\BCBL {}\ \BBA {} Opsomer%
}{%
Claeskens%
\ \protect \BOthers {.}}{%
{\protect \APACyear {2009}}%
}]{%
claeskens2009asymptotic}
\APACinsertmetastar {%
claeskens2009asymptotic}%
\begin{APACrefauthors}%
Claeskens, G.%
, Krivobokova, T.%
\BCBL {} Opsomer, J.D.%
\end{APACrefauthors}%
\unskip\
\newblock
\APACrefYearMonthDay{2009}{}{}.
\newblock
{\BBOQ}\APACrefatitle {Asymptotic properties of penalized spline estimators}
  {Asymptotic properties of penalized spline estimators}.{\BBCQ}
\newblock
\APACjournalVolNumPages{Biometrika}{96}{3}{529--544,}
\newblock
\begin{APACrefDOI} \doi{10.1093/biomet/asp035} \end{APACrefDOI}
\newblock

\newblock

\PrintBackRefs{\CurrentBib}

\bibitem [\protect \citeauthoryear {%
Dudley%
}{%
Dudley%
}{%
{\protect \APACyear {2016}}%
}]{%
dudley2016suprema}
\APACinsertmetastar {%
dudley2016suprema}%
\begin{APACrefauthors}%
Dudley, R.M.%
\end{APACrefauthors}%
\unskip\
\newblock
\APACrefYearMonthDay{2016}{}{}.
\newblock
{\BBOQ}\APACrefatitle {V. {N}. {S}udakov's work on expected suprema of
  {G}aussian processes} {V. {N}. {S}udakov's work on expected suprema of
  {G}aussian processes}.{\BBCQ}
\newblock
 \APACrefbtitle {{High Dimensional Probability VII: The Carg\`{e}se Volume}}
  {{High Dimensional Probability VII: The Carg\`{e}se Volume}}\ (\BPGS\
  37--43).
\newblock
\APACaddressPublisher{Cham}{Springer International Publishing}.
\PrintBackRefs{\CurrentBib}

\bibitem [\protect \citeauthoryear {%
El-Assi%
, Salah~Mahmoud%
\BCBL {}\ \BBA {} Nurul~Habib%
}{%
El-Assi%
\ \protect \BOthers {.}}{%
{\protect \APACyear {2017}}%
}]{%
el2017effects}
\APACinsertmetastar {%
el2017effects}%
\begin{APACrefauthors}%
El-Assi, W.%
, Salah~Mahmoud, M.%
\BCBL {} Nurul~Habib, K.%
\end{APACrefauthors}%
\unskip\
\newblock
\APACrefYearMonthDay{2017}{}{}.
\newblock
{\BBOQ}\APACrefatitle {Effects of built environment and weather on bike sharing
  demand: a station level analysis of commercial bike sharing in {Toronto}}
  {Effects of built environment and weather on bike sharing demand: a station
  level analysis of commercial bike sharing in {Toronto}}.{\BBCQ}
\newblock
\APACjournalVolNumPages{Transportation}{44}{3}{589--613,}
\newblock
\begin{APACrefDOI} \doi{10.1007/s11116-015-9669-z} \end{APACrefDOI}
\newblock

\newblock

\PrintBackRefs{\CurrentBib}

\bibitem [\protect \citeauthoryear {%
Faghih-Imani%
, Eluru%
, El-Geneidy%
, Rabbat%
\BCBL {}\ \BBA {} Haq%
}{%
Faghih-Imani%
\ \protect \BOthers {.}}{%
{\protect \APACyear {2014}}%
}]{%
faghih2014land}
\APACinsertmetastar {%
faghih2014land}%
\begin{APACrefauthors}%
Faghih-Imani, A.%
, Eluru, N.%
, El-Geneidy, A.M.%
, Rabbat, M.%
\BCBL {} Haq, U.%
\end{APACrefauthors}%
\unskip\
\newblock
\APACrefYearMonthDay{2014}{}{}.
\newblock
{\BBOQ}\APACrefatitle {How land-use and urban form impact bicycle flows:
  evidence from the bicycle-sharing system ({BIXI}) in {Montreal}} {How
  land-use and urban form impact bicycle flows: evidence from the
  bicycle-sharing system ({BIXI}) in {Montreal}}.{\BBCQ}
\newblock
\APACjournalVolNumPages{Journal of Transport Geography}{41}{}{306--314,}
\newblock
\begin{APACrefDOI} \doi{10.1016/j.jtrangeo.2014.01.013} \end{APACrefDOI}
\newblock

\newblock

\PrintBackRefs{\CurrentBib}

\bibitem [\protect \citeauthoryear {%
Fan%
\ \BBA {} Zhang%
}{%
Fan%
\ \BBA {} Zhang%
}{%
{\protect \APACyear {2000}}%
}]{%
fan2000simultaneous}
\APACinsertmetastar {%
fan2000simultaneous}%
\begin{APACrefauthors}%
Fan, J.%
\BCBT {}\ \BBA {} Zhang, W.%
\end{APACrefauthors}%
\unskip\
\newblock
\APACrefYearMonthDay{2000}{}{}.
\newblock
{\BBOQ}\APACrefatitle {Simultaneous confidence bands and hypothesis testing in
  varying-coefficient models} {Simultaneous confidence bands and hypothesis
  testing in varying-coefficient models}.{\BBCQ}
\newblock
\APACjournalVolNumPages{Scandinavian Journal of Statistics}{27}{4}{715--731,}
\newblock
\begin{APACrefDOI} \doi{10.1111/1467-9469.00218} \end{APACrefDOI}
\newblock

\newblock

\PrintBackRefs{\CurrentBib}

\bibitem [\protect \citeauthoryear {%
G\"{o}tze%
\ \BBA {} Zaitsev%
}{%
G\"{o}tze%
\ \BBA {} Zaitsev%
}{%
{\protect \APACyear {2009}}%
}]{%
gotze2008bounds}
\APACinsertmetastar {%
gotze2008bounds}%
\begin{APACrefauthors}%
G\"{o}tze, F.%
\BCBT {}\ \BBA {} Zaitsev, A.Y.%
\end{APACrefauthors}%
\unskip\
\newblock
\APACrefYearMonthDay{2009}{}{}.
\newblock
{\BBOQ}\APACrefatitle {Bounds for the rate of strong approximation in the
  multidimensional invariance principle} {Bounds for the rate of strong
  approximation in the multidimensional invariance principle}.{\BBCQ}
\newblock
\APACjournalVolNumPages{Theory of Probability and Its
  Applications}{53}{1}{59--80,}
\newblock
\begin{APACrefDOI} \doi{10.1137/S0040585X9798350X} \end{APACrefDOI}
\newblock

\newblock

\PrintBackRefs{\CurrentBib}

\bibitem [\protect \citeauthoryear {%
Gregory%
, Mammen%
\BCBL {}\ \BBA {} Wahl%
}{%
Gregory%
\ \protect \BOthers {.}}{%
{\protect \APACyear {2021}}%
}]{%
gregory2021statistical}
\APACinsertmetastar {%
gregory2021statistical}%
\begin{APACrefauthors}%
Gregory, K.%
, Mammen, E.%
\BCBL {} Wahl, M.%
\end{APACrefauthors}%
\unskip\
\newblock
\APACrefYearMonthDay{2021}{}{}.
\newblock
{\BBOQ}\APACrefatitle {Statistical inference in sparse high-dimensional
  additive models} {Statistical inference in sparse high-dimensional additive
  models}.{\BBCQ}
\newblock
\APACjournalVolNumPages{The Annals of Statistics}{49}{3}{1514--1536,}
\newblock
\begin{APACrefDOI} \doi{10.1214/20-aos2011} \end{APACrefDOI}
\newblock

\newblock

\PrintBackRefs{\CurrentBib}

\bibitem [\protect \citeauthoryear {%
Hall%
\ \BBA {} Horowitz%
}{%
Hall%
\ \BBA {} Horowitz%
}{%
{\protect \APACyear {2013}}%
}]{%
hall2013simple}
\APACinsertmetastar {%
hall2013simple}%
\begin{APACrefauthors}%
Hall, P.%
\BCBT {}\ \BBA {} Horowitz, J.%
\end{APACrefauthors}%
\unskip\
\newblock
\APACrefYearMonthDay{2013}{}{}.
\newblock
{\BBOQ}\APACrefatitle {A simple bootstrap method for constructing nonparametric
  confidence bands for functions} {A simple bootstrap method for constructing
  nonparametric confidence bands for functions}.{\BBCQ}
\newblock
\APACjournalVolNumPages{The Annals of Statistics}{41}{4}{1892--1921,}
\newblock
\begin{APACrefDOI} \doi{10.1214/13-AOS1137} \end{APACrefDOI}
\newblock

\newblock

\PrintBackRefs{\CurrentBib}

\bibitem [\protect \citeauthoryear {%
He%
\ \BBA {} Huang%
}{%
He%
\ \BBA {} Huang%
}{%
{\protect \APACyear {2009}}%
}]{%
he2009double}
\APACinsertmetastar {%
he2009double}%
\begin{APACrefauthors}%
He, H.%
\BCBT {}\ \BBA {} Huang, L\BHBI S.%
\end{APACrefauthors}%
\unskip\
\newblock
\APACrefYearMonthDay{2009}{}{}.
\newblock
{\BBOQ}\APACrefatitle {Double-smoothing for bias reduction in local linear
  regression} {Double-smoothing for bias reduction in local linear
  regression}.{\BBCQ}
\newblock
\APACjournalVolNumPages{Journal of Statistical Planning and
  Inference}{139}{3}{1056--1072,}
\newblock
\begin{APACrefDOI} \doi{10.1016/j.jspi.2008.06.011} \end{APACrefDOI}
\newblock

\newblock

\PrintBackRefs{\CurrentBib}

\bibitem [\protect \citeauthoryear {%
Hong%
, Lee%
\BCBL {}\ \BBA {} Park%
}{%
Hong%
\ \protect \BOthers {.}}{%
{\protect \APACyear {2023}}%
}]{%
hong2023additive}
\APACinsertmetastar {%
hong2023additive}%
\begin{APACrefauthors}%
Hong, H.%
, Lee, Y.K.%
\BCBL {} Park, B.U.%
\end{APACrefauthors}%
\unskip\
\newblock
\APACrefYearMonthDay{2023}{}{}.
\newblock
{\BBOQ}\APACrefatitle {Additive regression with parametric help} {Additive
  regression with parametric help}.{\BBCQ}
\newblock
\APACjournalVolNumPages{Bernoulli}{29}{4}{3059--3092,}
\newblock
\begin{APACrefDOI} \doi{10.3150/22-bej1575} \end{APACrefDOI}
\newblock

\newblock

\PrintBackRefs{\CurrentBib}

\bibitem [\protect \citeauthoryear {%
Jeon%
, Lee%
, Mammen%
\BCBL {}\ \BBA {} Park%
}{%
Jeon%
\ \protect \BOthers {.}}{%
{\protect \APACyear {2022}}%
}]{%
jeon2022locally}
\APACinsertmetastar {%
jeon2022locally}%
\begin{APACrefauthors}%
Jeon, J.M.%
, Lee, Y.K.%
, Mammen, E.%
\BCBL {} Park, B.U.%
\end{APACrefauthors}%
\unskip\
\newblock
\APACrefYearMonthDay{2022}{}{}.
\newblock
{\BBOQ}\APACrefatitle {Locally polynomial {H}ilbertian additive regression}
  {Locally polynomial {H}ilbertian additive regression}.{\BBCQ}
\newblock
\APACjournalVolNumPages{Bernoulli}{28}{3}{2034--2066,}
\newblock
\begin{APACrefDOI} \doi{10.3150/21-bej1410} \end{APACrefDOI}
\newblock

\newblock

\PrintBackRefs{\CurrentBib}

\bibitem [\protect \citeauthoryear {%
Jeon%
\ \BBA {} Park%
}{%
Jeon%
\ \BBA {} Park%
}{%
{\protect \APACyear {2020}}%
}]{%
jeon2020additive}
\APACinsertmetastar {%
jeon2020additive}%
\begin{APACrefauthors}%
Jeon, J.M.%
\BCBT {}\ \BBA {} Park, B.U.%
\end{APACrefauthors}%
\unskip\
\newblock
\APACrefYearMonthDay{2020}{}{}.
\newblock
{\BBOQ}\APACrefatitle {Additive regression with {H}ilbertian responses}
  {Additive regression with {H}ilbertian responses}.{\BBCQ}
\newblock
\APACjournalVolNumPages{The Annals of Statistics}{48}{5}{2671--2697,}
\newblock
\begin{APACrefDOI} \doi{10.1214/19-AOS1902} \end{APACrefDOI}
\newblock

\newblock

\PrintBackRefs{\CurrentBib}

\bibitem [\protect \citeauthoryear {%
Leadbetter%
, Lindgren%
\BCBL {}\ \BBA {} Rootz\'{e}n%
}{%
Leadbetter%
\ \protect \BOthers {.}}{%
{\protect \APACyear {1983}}%
}]{%
leadbetter1983extremes}
\APACinsertmetastar {%
leadbetter1983extremes}%
\begin{APACrefauthors}%
Leadbetter, M.R.%
, Lindgren, G.%
\BCBL {} Rootz\'{e}n, H.%
\end{APACrefauthors}%
\unskip\
\newblock
\APACrefYear{1983}.
\newblock
\APACrefbtitle {Extremes and related properties of random sequences and
  processes} {Extremes and related properties of random sequences and
  processes}.
\newblock
\APACaddressPublisher{}{Springer}.
\PrintBackRefs{\CurrentBib}

\bibitem [\protect \citeauthoryear {%
E.R.~Lee%
, Park%
, Mammen%
\BCBL {}\ \BBA {} Park%
}{%
E.R.~Lee%
\ \protect \BOthers {.}}{%
{\protect \APACyear {2024}}%
}]{%
lee2024efficient}
\APACinsertmetastar {%
lee2024efficient}%
\begin{APACrefauthors}%
Lee, E.R.%
, Park, S.%
, Mammen, E.%
\BCBL {} Park, B.U.%
\end{APACrefauthors}%
\unskip\
\newblock
\APACrefYearMonthDay{2024}{}{}.
\newblock
{\BBOQ}\APACrefatitle {Efficient functional {L}asso kernel smoothing for
  high-dimensional additive regression} {Efficient functional {L}asso kernel
  smoothing for high-dimensional additive regression}.{\BBCQ}
\newblock
\APACjournalVolNumPages{The Annals of Statistics}{52}{4}{1741--1773,}
\newblock
\begin{APACrefDOI} \doi{10.1214/24-aos2415} \end{APACrefDOI}
\newblock

\newblock

\PrintBackRefs{\CurrentBib}

\bibitem [\protect \citeauthoryear {%
Y.K.~Lee%
, Mammen%
, Nielsen%
\BCBL {}\ \BBA {} Park%
}{%
Y.K.~Lee%
\ \protect \BOthers {.}}{%
{\protect \APACyear {2020}}%
}]{%
lee2020nonparametric}
\APACinsertmetastar {%
lee2020nonparametric}%
\begin{APACrefauthors}%
Lee, Y.K.%
, Mammen, E.%
, Nielsen, J.P.%
\BCBL {} Park, B.U.%
\end{APACrefauthors}%
\unskip\
\newblock
\APACrefYearMonthDay{2020}{}{}.
\newblock
{\BBOQ}\APACrefatitle {Nonparametric regression with parametric help}
  {Nonparametric regression with parametric help}.{\BBCQ}
\newblock
\APACjournalVolNumPages{Electronic Journal of Statistics}{14}{2}{3845--3868,}
\newblock
\begin{APACrefDOI} \doi{10.1214/20-EJS1760} \end{APACrefDOI}
\newblock

\newblock

\PrintBackRefs{\CurrentBib}

\bibitem [\protect \citeauthoryear {%
Li%
\ \BBA {} Ruppert%
}{%
Li%
\ \BBA {} Ruppert%
}{%
{\protect \APACyear {2008}}%
}]{%
li2008asymptotics}
\APACinsertmetastar {%
li2008asymptotics}%
\begin{APACrefauthors}%
Li, Y.%
\BCBT {}\ \BBA {} Ruppert, D.%
\end{APACrefauthors}%
\unskip\
\newblock
\APACrefYearMonthDay{2008}{}{}.
\newblock
{\BBOQ}\APACrefatitle {On the asymptotics of penalized splines} {On the
  asymptotics of penalized splines}.{\BBCQ}
\newblock
\APACjournalVolNumPages{Biometrika}{95}{2}{415--436,}
\newblock
\begin{APACrefDOI} \doi{10.1093/biomet/asn010} \end{APACrefDOI}
\newblock

\newblock

\PrintBackRefs{\CurrentBib}

\bibitem [\protect \citeauthoryear {%
Liu%
, Yang%
\BCBL {}\ \BBA {} H\"{a}rdle%
}{%
Liu%
\ \protect \BOthers {.}}{%
{\protect \APACyear {2013}}%
}]{%
liu2013oracally}
\APACinsertmetastar {%
liu2013oracally}%
\begin{APACrefauthors}%
Liu, R.%
, Yang, L.%
\BCBL {} H\"{a}rdle, W.K.%
\end{APACrefauthors}%
\unskip\
\newblock
\APACrefYearMonthDay{2013}{}{}.
\newblock
{\BBOQ}\APACrefatitle {Oracally efficient two-step estimation of generalized
  additive model} {Oracally efficient two-step estimation of generalized
  additive model}.{\BBCQ}
\newblock
\APACjournalVolNumPages{Journal of the American Statistical
  Association}{108}{502}{619--631,}
\newblock
\begin{APACrefDOI} \doi{10.1080/01621459.2013.763726} \end{APACrefDOI}
\newblock

\newblock

\PrintBackRefs{\CurrentBib}

\bibitem [\protect \citeauthoryear {%
Mammen%
, Linton%
\BCBL {}\ \BBA {} Nielsen%
}{%
Mammen%
\ \protect \BOthers {.}}{%
{\protect \APACyear {1999}}%
}]{%
mammen1999existence}
\APACinsertmetastar {%
mammen1999existence}%
\begin{APACrefauthors}%
Mammen, E.%
, Linton, O.%
\BCBL {} Nielsen, J.P.%
\end{APACrefauthors}%
\unskip\
\newblock
\APACrefYearMonthDay{1999}{}{}.
\newblock
{\BBOQ}\APACrefatitle {The existence and asymptotic properties of a backfitting
  projection algorithm under weak conditions} {The existence and asymptotic
  properties of a backfitting projection algorithm under weak
  conditions}.{\BBCQ}
\newblock
\APACjournalVolNumPages{The Annals of Statistics}{27}{5}{1443--1490,}
\newblock
\begin{APACrefDOI} \doi{10.1214/aos/1017939137} \end{APACrefDOI}
\newblock

\newblock

\PrintBackRefs{\CurrentBib}

\bibitem [\protect \citeauthoryear {%
Mammen%
\ \BBA {} Park%
}{%
Mammen%
\ \BBA {} Park%
}{%
{\protect \APACyear {2006}}%
}]{%
mammen2006simple}
\APACinsertmetastar {%
mammen2006simple}%
\begin{APACrefauthors}%
Mammen, E.%
\BCBT {}\ \BBA {} Park, B.U.%
\end{APACrefauthors}%
\unskip\
\newblock
\APACrefYearMonthDay{2006}{}{}.
\newblock
{\BBOQ}\APACrefatitle {A simple smooth backfitting method for additive models}
  {A simple smooth backfitting method for additive models}.{\BBCQ}
\newblock
\APACjournalVolNumPages{The Annals of Statistics}{34}{5}{2252--2271,}
\newblock
\begin{APACrefDOI} \doi{10.1214/009053606000000696} \end{APACrefDOI}
\newblock

\newblock

\PrintBackRefs{\CurrentBib}

\bibitem [\protect \citeauthoryear {%
Mammen%
\ \BBA {} Sperlich%
}{%
Mammen%
\ \BBA {} Sperlich%
}{%
{\protect \APACyear {2022}}%
}]{%
mammen2022backfitting}
\APACinsertmetastar {%
mammen2022backfitting}%
\begin{APACrefauthors}%
Mammen, E.%
\BCBT {}\ \BBA {} Sperlich, S.%
\end{APACrefauthors}%
\unskip\
\newblock
\APACrefYearMonthDay{2022}{}{}.
\newblock
{\BBOQ}\APACrefatitle {Backfitting tests in generalized structured models}
  {Backfitting tests in generalized structured models}.{\BBCQ}
\newblock
\APACjournalVolNumPages{Biometrika}{109}{1}{137--152,}
\newblock
\begin{APACrefDOI} \doi{10.1093/biomet/asaa108} \end{APACrefDOI}
\newblock

\newblock

\PrintBackRefs{\CurrentBib}

\bibitem [\protect \citeauthoryear {%
Moon%
, Park%
\BCBL {}\ \BBA {} Lee%
}{%
Moon%
\ \protect \BOthers {.}}{%
{\protect \APACyear {2025}}%
}]{%
moon2025varying}
\APACinsertmetastar {%
moon2025varying}%
\begin{APACrefauthors}%
Moon, S.H.%
, Park, B.U.%
\BCBL {} Lee, Y.K.%
\end{APACrefauthors}%
\unskip\
\newblock
\APACrefYearMonthDay{2025}{}{}.
\newblock
{\BBOQ}\APACrefatitle {Varying coefficient regression: revisit and parametric
  help} {Varying coefficient regression: revisit and parametric help}.{\BBCQ}
\newblock
\APACjournalVolNumPages{Bernoulli}{31}{4}{2569--2596,}
\newblock
\begin{APACrefDOI} \doi{10.3150/24-BEJ1817} \end{APACrefDOI}
\newblock

\newblock

\PrintBackRefs{\CurrentBib}

\bibitem [\protect \citeauthoryear {%
Park%
, Mammen%
, Lee%
\BCBL {}\ \BBA {} Lee%
}{%
Park%
\ \protect \BOthers {.}}{%
{\protect \APACyear {2015}}%
}]{%
park2015varying}
\APACinsertmetastar {%
park2015varying}%
\begin{APACrefauthors}%
Park, B.U.%
, Mammen, E.%
, Lee, Y.K.%
\BCBL {} Lee, E.R.%
\end{APACrefauthors}%
\unskip\
\newblock
\APACrefYearMonthDay{2015}{}{}.
\newblock
{\BBOQ}\APACrefatitle {Varying coefficient regression models: a review and new
  developments} {Varying coefficient regression models: a review and new
  developments}.{\BBCQ}
\newblock
\APACjournalVolNumPages{International Statistical Review}{83}{1}{36--64,}
\newblock
\begin{APACrefDOI} \doi{10.1111/insr.12029} \end{APACrefDOI}
\newblock

\newblock

\PrintBackRefs{\CurrentBib}

\bibitem [\protect \citeauthoryear {%
Petersen%
\ \BBA {} M\"{u}ller%
}{%
Petersen%
\ \BBA {} M\"{u}ller%
}{%
{\protect \APACyear {2019}}%
}]{%
petersen2019frechet}
\APACinsertmetastar {%
petersen2019frechet}%
\begin{APACrefauthors}%
Petersen, A.%
\BCBT {}\ \BBA {} M\"{u}ller, H\BHBI G.%
\end{APACrefauthors}%
\unskip\
\newblock
\APACrefYearMonthDay{2019}{}{}.
\newblock
{\BBOQ}\APACrefatitle {Fr\'{e}chet regression for random objects with
  {E}uclidean predictors} {Fr\'{e}chet regression for random objects with
  {E}uclidean predictors}.{\BBCQ}
\newblock
\APACjournalVolNumPages{The Annals of Statistics}{47}{2}{691--719,}
\newblock
\begin{APACrefDOI} \doi{10.1214/17-AOS1624} \end{APACrefDOI}
\newblock

\newblock

\PrintBackRefs{\CurrentBib}

\bibitem [\protect \citeauthoryear {%
Shang%
\ \BBA {} Cheng%
}{%
Shang%
\ \BBA {} Cheng%
}{%
{\protect \APACyear {2013}}%
}]{%
shang2013local}
\APACinsertmetastar {%
shang2013local}%
\begin{APACrefauthors}%
Shang, Z.%
\BCBT {}\ \BBA {} Cheng, G.%
\end{APACrefauthors}%
\unskip\
\newblock
\APACrefYearMonthDay{2013}{}{}.
\newblock
{\BBOQ}\APACrefatitle {Local and global asymptotic inference in smoothing
  spline models} {Local and global asymptotic inference in smoothing spline
  models}.{\BBCQ}
\newblock
\APACjournalVolNumPages{The Annals of Statistics}{41}{5}{2608--2638,}
\newblock
\begin{APACrefDOI} \doi{10.1214/13-AOS1164} \end{APACrefDOI}
\newblock

\newblock

\PrintBackRefs{\CurrentBib}

\bibitem [\protect \citeauthoryear {%
Tsybakov%
}{%
Tsybakov%
}{%
{\protect \APACyear {2009}}%
}]{%
tsybakov2009introduction}
\APACinsertmetastar {%
tsybakov2009introduction}%
\begin{APACrefauthors}%
Tsybakov, A.B.%
\end{APACrefauthors}%
\unskip\
\newblock
\APACrefYear{2009}.
\newblock
\APACrefbtitle {Introduction to nonparametric estimation} {Introduction to
  nonparametric estimation}\ (\BVOL~11).
\newblock
\APACaddressPublisher{}{Springer}.
\PrintBackRefs{\CurrentBib}

\bibitem [\protect \citeauthoryear {%
Xia%
}{%
Xia%
}{%
{\protect \APACyear {1998}}%
}]{%
xia1998bias}
\APACinsertmetastar {%
xia1998bias}%
\begin{APACrefauthors}%
Xia, Y.%
\end{APACrefauthors}%
\unskip\
\newblock
\APACrefYearMonthDay{1998}{}{}.
\newblock
{\BBOQ}\APACrefatitle {Bias-corrected confidence bands in nonparametric
  regression} {Bias-corrected confidence bands in nonparametric
  regression}.{\BBCQ}
\newblock
\APACjournalVolNumPages{Journal of the Royal Statistical Society Series B:
  Statistical Methodology}{60}{4}{797--811,}
\newblock
\begin{APACrefDOI} \doi{10.1111/1467-9868.00155} \end{APACrefDOI}
\newblock

\newblock

\PrintBackRefs{\CurrentBib}

\bibitem [\protect \citeauthoryear {%
Yu%
, Mammen%
\BCBL {}\ \BBA {} Park%
}{%
Yu%
\ \protect \BOthers {.}}{%
{\protect \APACyear {2011}}%
}]{%
yu2011semi}
\APACinsertmetastar {%
yu2011semi}%
\begin{APACrefauthors}%
Yu, K.%
, Mammen, E.%
\BCBL {} Park, B.U.%
\end{APACrefauthors}%
\unskip\
\newblock
\APACrefYearMonthDay{2011}{}{}.
\newblock
{\BBOQ}\APACrefatitle {Semi-parametric regression: efficiency gains from
  modeling the nonparametric part} {Semi-parametric regression: efficiency
  gains from modeling the nonparametric part}.{\BBCQ}
\newblock
\APACjournalVolNumPages{Bernoulli}{17}{2}{736--748,}
\newblock
\begin{APACrefDOI} \doi{10.3150/10-BEJ296} \end{APACrefDOI}
\newblock

\newblock

\PrintBackRefs{\CurrentBib}

\bibitem [\protect \citeauthoryear {%
Zhang%
\ \BBA {} Peng%
}{%
Zhang%
\ \BBA {} Peng%
}{%
{\protect \APACyear {2010}}%
}]{%
zhang2010simultaneous}
\APACinsertmetastar {%
zhang2010simultaneous}%
\begin{APACrefauthors}%
Zhang, W.%
\BCBT {}\ \BBA {} Peng, H.%
\end{APACrefauthors}%
\unskip\
\newblock
\APACrefYearMonthDay{2010}{}{}.
\newblock
{\BBOQ}\APACrefatitle {Simultaneous confidence band and hypothesis test in
  generalised varying-coefficient models} {Simultaneous confidence band and
  hypothesis test in generalised varying-coefficient models}.{\BBCQ}
\newblock
\APACjournalVolNumPages{Journal of Multivariate Analysis}{101}{7}{1656--1680,}
\newblock
\begin{APACrefDOI} \doi{10.1016/j.jmva.2010.03.003} \end{APACrefDOI}
\newblock

\newblock

\PrintBackRefs{\CurrentBib}

\end{thebibliography}

\clearpage

\counterwithout{equation}{section}

\renewcommand{\thesection}{S.\arabic{section}}
\renewcommand{\theequation}{S.\arabic{equation}}
\renewcommand{\theHsection}{supp.\arabic{section}}
\renewcommand{\theHsubsection}{supp.\arabic{section}.\arabic{subsection}}
\renewcommand{\thelemma}{S.\arabic{lemma}}
\setcounter{section}{0}
\setcounter{equation}{0}

\begin{center}
	{\Large\bf SUPPLEMENTARY MATERIAL}
\end{center}
{\bf Organization.}
Section~\ref{suppsec:LL_SBF} provides details on local linear SBF estimation and its computation.
Section~\ref{suppsec:tildem} presents the theoretical properties and numerical results for bias-reduced estimators based on the direct use of $\hat{\alpha}_{j,L}$, as introduced in Section~\ref{subsec:BR_SBF}.
Section~\ref{suppsec:proof_theorems} contains the proofs of Proposition~\ref{prop:bias_expansion}, Theorems~\ref{thm:mtilde_asymp}, \ref{thm:mhat_asymp} and \ref{thm:unif_CLT}, and Corollary~\ref{cor:unif_CLT}.
The proofs of additional theorems and lemmas are provided in Section~\ref{suppsec:proof_lemmas}.

\section{Local linear SBF estimation and computation}\label{suppsec:LL_SBF}

This section provides the construction and computational details of the local linear SBF estimator introduced in Section~\ref{subsec:LL_SBF}.

Here, we do not distinguish between a univariate function $f_j:[0,1]\to\bbR$ and its extension to $[0,1]^d$ defined by $\bmx\mapsto f_j(x_j)$.
With this slight abuse of notation, let $m_+ := \E(Y|\bmX=\cdot) = m_0 + m_1 + \cdots + m_d : [0,1]^d \to \bbR$ denote the true additive mean function.
The local linear smoothing considers the first-order local approximation of $m_+$ around each point $\bmx\in[0,1]^d$:
\begin{equation*}
	m_+(\bmu) \simeq m_+(\bmx) + \sumjd (u_j - x_j)\frac{\partial m_j}{\partial x_j}(x_j) \quad \text{for} \quad \bmu\simeq \bmx.
\end{equation*}
Let $\bar{Y}:=\meanin Y_i$ and $\bmm^{\F,\tp} = (\bar{Y}+m_1 +\cdots + m_d, m_{1,1}, \dots, m_{1,d})^{\top} : [0,1]^d \to \bbR^{d+1}$ denote a tuple of true additive component functions.
We then define the estimator $\hatbmm^{\F,\tp} := (\bar{Y}+\hatm_1+\cdots+\hatm_d,\hatm_{1,1},\dots,\hatm_{1,d})^\tp:[0,1]^d\to\bbR^{d+1}$ as the minimizer of
\begin{equation*}
	\mathcal{L}(\bmf^{\F,\tp}) = \int_{[0,1]^d}	\Bigg[	\meanin \bigg\{	Y_i - f_+(\bmx) - \sumjd \bigg(\frac{X_{ij}-x_j}{h_j}\bigg) f_{1,j}(\bmx) \bigg\}^2 \cdot \prodjd K_{h_j}(x_j,X_{ij}) \Bigg]\,\d\bmx
\end{equation*}
over the space of tuples $\bmf^{\F,\tp} = (f_+, f_{1,1},\dots, f_{1,d})^{\top}$ where $f_+(\bmx) = f_1(x_1) + \cdots + f_d(x_d)$ for some univariate functions $f_j: [0,1] \to \bbR$.

We define the marginal local linear estimator of $\E(Y| X_j=x_j)$ as
\begin{equation*}
	\tildebmm_j^{\tp}(x_j) := \hatbmM_{jj}(x_j)^{-1} \cdot \meanin  \bmvec\bigg(\frac{X_{ij}-x_j}{h_j}\bigg)K_{h_j}(x_j,X_{ij})Y_i.
\end{equation*}
If $\hatbmM_{jj}(x_j)$ is invertible for all $1\le j\le d$ and $x_j\in[0,1]$, the estimator $\hatbmm^{\F,\tp}$ satisfies the following system of integral equations: for $1\le j\le d$,
\begin{equation}\label{eq:SBF_integral_eq}
	\hatbmm_j^{\tp}(x_j) = \tildebmm_j^{\tp}(x_j) -\bar{Y} \cdot \fraku_0 - \sum_{k\neq j} \int_0^1 \hatbmM_{jj}(x_j)^{-1}	\hatbmM_{jk}(x_j,x_k) \hatbmm_k^{\tp}(x_k) \d x_k,
\end{equation}
where $\hatbmm_j^\tp:=(\hatm_j,\hatm_{1,j})^\top:[0,1]\to\bbR^2$ for $1\le j\le d$ and $\fraku_0:=(1,0)^\top\in\bbR^2$.

According to \cite{jeon2022locally}, under mild conditions, there exists a unique tuple $\hatbmm^{\F,\tp} = (\bar{Y}+\hatm_1 +\cdots + \hatm_d, \hatm_{1,1}, \dots, \hatm_{1,d})^{\top}$ satisfying the SBF equations \eqref{eq:SBF_integral_eq} with probability tending to one.
However, the individual components $\hatm_j$ are not uniquely determined without additional constraints.
Indeed, for arbitrary constants $c_1,\ldots,c_d \in \bbR$ such that $\sumjd c_j=0$, the two additive decompositions
\begin{align*}
	(\hatm_1,\hatm_2,\ldots,\hatm_d) \quad\text{and}\quad (\hatm_1+c_1,\hatm_2+c_2,\ldots,\hatm_d+c_d)
\end{align*}
constitute the same full tuple $\hatbmm^{\F,\tp}$ and satisfy the SBF equations at \eqref{eq:SBF_integral_eq}.
Note that the estimators $\hatm_{1,j}$ are uniquely determined without any constraints.
To ensure the identifiability of $\hatm_j$, we impose the following constraints:
\begin{equation}\label{eq:SBF_constraint}
	\fraku_0^\top \cdot \int_0^1 \hatbmM_{jj}(x_j) \hatbmm_j^{\tp}(x_j) \d x_j=0 \quad \text{for} \quad 1\le j\le d.
\end{equation}

To solve the SBF equations \eqref{eq:SBF_integral_eq}, we adopt an iterative backfitting algorithm.
For $1\le j\le d$, let $\hatbmm_j^{[0],\tp} = (\hatm_j^{[0]},\hatm_{1,j}^{[0]})^\top:[0,1]\to\bbR^2$ denote an initial tuple of additive components such that $\fraku_0^\top \cdot \int_0^1 \hatbmM_{jj}(x_j)\hatbmm_j^{[0],\tp}(x_j) \d x_j=0$.
At each iteration step, the $j$-th additive component tuple is updated using the most recently available estimates of the remaining components.
Specifically, for $\ell\ge1$, the estimators $\hatbmm_j^{[\ell],\tp}:=(\hatm_j^{[\ell]},\hatm_{1,j}^{[\ell]})^\top$ at the $\ell$-th iteration are updated sequentially across $1\le j\le d$ as follows:
\begin{align*}
	\hatbmm_j^{[\ell],\tp}(x_j) &= \tildebmm_j^{\tp}(x_j) - \bar{Y} \cdot \fraku_0 - \sum_{k \le j-1} \int_0^1 \hatbmM_{jj}(x_j)^{-1} \hatbmM_{jk}(x_j,x_k) \hatbmm_k^{[\ell],\tp}(x_k)\, \d x_k \\
	&\hspace{2cm} - \sum_{k \ge j+1} \int_0^1 \hatbmM_{jj}(x_j)^{-1} \hatbmM_{jk}(x_j,x_k) \hatbmm_k^{[\ell-1],\tp}(x_k)\, \d x_k.
\end{align*}
By the choice of the initial additive tuples, all subsequent updates $\hatbmm_j^{[\ell],\tp}$ inherently satisfy $\fraku_0^\top \cdot \int_0^1 \hatbmM_{jj}(x_j)\hatbmm_j^{[\ell],\tp}(x_j) \d x_j=0$ for all $\ell\ge1$ and $1\le j\le d$.
As the initial tuples, we may use $\hatbmm_j^{[0],\tp} = \tildebmm_j^{\tp} - \bar{Y} \cdot \fraku_0$ that satisfies the constraint $\fraku_0^\top \cdot \int_0^1 \hatbmM_{jj}(x_j)\hatbmm_j^{[0],\tp}(x_j) \d x_j=0$.
This initialization is used throughout the numerical studies in Sections~\ref{sec:Numerical_Study} and \ref{sec:real_data}.

For $\ell\ge0$, let $\hatbmm^{[\ell],\F,\tp}:= (\bar{Y}+\hatm_1^{[\ell]}+\cdots+\hatm_d^{[\ell]},\hatm_{1,1}^{[\ell]},\cdots,\hatm_{1,d}^{[\ell]})^\top$.
To analyze the convergence of the SBF algorithm and to formulate a stopping rule, we define a semi-norm for the full tuples $\bmf^{\F,\tp} = (f_+, f_{1,1}, \dots, f_{1,d})^{\top}$ with $f_+(\bmx) = f_1(x_1) + \cdots + f_d(x_d)$ for some $f_j:[0,1]\to\bbR$.
For $\bmx\in[0,1]^d$, let $\hatp(\bmx) := \meanin \prodjd K_{h_j}(x_j,X_{ij})$ and $\hatp_j(x_j) := \meanin K_{h_j}(x_j,X_{ij})$.
Using these kernel density estimators, define the semi-norm $\|\cdot\|_{2,n}$ by
\begin{align*}
	\|\bmf^{\F,\tp}\|_{2,n}^2 := \int_{[0,1]^d} \abs{f_+(\bmx)}^2 \hatp(\bmx)\d\bmx + \sumjd \int_0^1 \abs{f_{1,j}(x_j)}^2 \hatp_j(x_j)\d x_j.
\end{align*}
By Theorem 1 of \cite{jeon2022locally}, under mild regularity conditions, the following holds with probability tending to one: there exist constants $0<C<\infty$ and $0<\gamma<1$ such that, if the initial estimator satisfies $\|\hatbmm^{[0],\F,\tp}\|_{2,n} \le C_0$ for some constant $0<C_0<\infty$, then $\| \hatbmm^{[\ell],\F,\tp} - \hatbmm^{\F,\tp} \|_{2,n} \le CC_0\gamma^\ell$ for all $\ell\ge1$.
Hence, the iterative SBF update converges to the unique SBF solution $\hatbmm^{\F,\tp}$ satisfying \eqref{eq:SBF_integral_eq}.
In practice, given a small $\epsilon>0$ (e.g., $\epsilon=10^{-5}$), we terminate the SBF algorithm when $\| \hatbmm^{[\ell],\F,\tp} - \hatbmm^{[\ell-1],\F,\tp} \|_{2,n} < \epsilon$.

\section{Bias-reduced estimators based on $\hat{\alpha}_{j,L}$}\label{suppsec:tildem}
In this section, we present theoretical properties and numerical results for bias-reduced estimators based on the direct use of $\hat{\alpha}_{j,L}$.
This bias-reduction method was introduced in Section~\ref{subsec:BR_SBF}.
Analogously to the weighted least-squares estimator $\hatbeta_{j,L}$, we can consider the weighted estimator $\hatalpha_{j,L}^{\rm W}$, with weights proportional to $b_{j,\ell}$.
The resulting alternative bias-reduced estimators are given by
\begin{align*}
	\tildem_j^{\rm BR} :=& \, \hatalpha_{j,L} = \,
	\suminl \Bigg\{\frac{\suminlp b_{j,\ell'}^4 - \big(\suminlp b_{j,\ell'}^2 \big) b_{j,\ell}^2}{L\cdot \suminlp b_{j,\ell'}^4 - \big(\suminlp b_{j,\ell'}^2 \big)^2} \Bigg\} \cdot \hatm_j(\cdot;\bmb_\ell), \\
	\tildem_j^{\rm WBR}	:=& \, \hatalpha_{j,L}^{\rm W} = \,
	\suminl	\Bigg\{	\frac{ \big(\suminlp b_{j,\ell'}^5\big)b_{j,\ell} - \big(\suminlp b_{j,\ell'}^3 \big) b_{j,\ell}^3}{ \big(\suminlp b_{j,\ell'} \big) \big( \suminlp b_{j,\ell'}^5\big) - \big(\suminlp b_{j,\ell'}^3 \big)^2} \Bigg\} \cdot \hatm_j(\cdot;\bmb_\ell).
\end{align*}

\subsection{Theoretical properties of $\tildem_j^\BR$ and $\tildem_j^\WBR$}
In this subsection, we establish the convergence rates for $\tildem_j^\BR$ and $\tildem_j^\WBR$.
\begin{theorem}\label{thm:mtilde_asymp}
	Assume (A1)--(A6).
	Then, for $\#=\BR$ and $\#=\WBR$, the following convergence rates hold:
	\begin{itemize}
		\item[(i)] Pointwise convergence:
		\begin{align*}
			| \tildem_j^{\#}(x_j) - m_j(x_j) | &= O_p\big(n^{-1/2}b_j^{-1/2}\big) +  o_p\big(b_j^2 + b_{-j}^2\big), \quad x_j\in I_j, \\
			| \tildem_j^{\#}(x_j) - m_j(x_j) | &= O_p\big(n^{-1/2}b_j^{-1/2} + b_j^2\big) + o_p\big(b_{-j}^2\big), \quad x_j\in [0,1]\setminus I_j.
		\end{align*}
		
		\item[(ii)] $L_2$ convergence:
		\begin{align*}
			\int_{I_j} | \tildem_j^{\#}(x_j) - m_j(x_j) |^2 \d x_j &= O_p\big(n^{-1}b_j^{-1}\big) + o_p\bigr(b_j^4 + b_{-j}^4\bigr), \\
			\int_0^1 | \tildem_j^{\#}(x_j) - m_j(x_j) |^2 \d x_j &= O_p\big(n^{-1}b_j^{-1}\big) + o_p\bigr(b_j^4 + b_{-j}^4\bigr).
		\end{align*}
		
		\item[(iii)] Uniform convergence:
		\begin{align*}
			\sup_{x_j \in I_j} | \tildem_j^{\#}(x_j) - m_j(x_j) | &= O_p\bigg(\sqrt{\frac{\log n}{n b_j}} \bigg) + o_p\big(b_j^2 + b_{-j}^2\big), \\
			\sup_{x_j\in[0,1]} | \tildem_j^{\#}(x_j) - m_j(x_j) | &= O_p\bigg(\sqrt{\frac{\log n}{n b_j}} + b_j^2 \bigg)+ o_p\big(b_{-j}^2\big).
		\end{align*}	
	\end{itemize}
\end{theorem}

Theorem~\ref{thm:mtilde_asymp} shows that, as is the case for $\hatm_j^\BR$ and $\hatm_j^\WBR$, the leading bias term of order $O_p(b_j^2)$ is reduced to $o_p(b_j^2)$.
Theorem~\ref{thm:mtilde_asymp} also shows that $\tildem_j^\BR$ and $\tildem_j^\WBR$ exhibit boundary effects for the same reason as $\hatm_j^\BR$ and $\hatm_j^\WBR$ in Section~\ref{subsec:theory_BR_SBF}.
It is worth noting that the convergence rates of $\tildem_j^\BR$ and $\tildem_j^\WBR$ are governed by the bandwidth vectors $\{\bmb_{\ell}:1\le\ell\le L\}$ rather than by $\bmh$, unlike $\hatm_j^\BR$ and $\hatm_j^\WBR$.
Indeed, both estimators directly estimate $m_j$ as weighted sums of $\hatm_j(\cdot;\bmb_{\ell})$.
Since $\tildem_j^\BR$ and $\tildem_j^\WBR$ estimate $m_j$ directly as weighted sums of $\hatm_j(\cdot;\bmb_{\ell})$, their convergence rates and interior regions depend only on $\{\bmb_{\ell}:1\le\ell\le L\}$.
In contrast, the convergence rates and interior regions of $\hatm_j^\BR$ and $\hatm_j^\WBR$ depend on both $\bmh$ and the bandwidth vectors used in their construction.
The proof of Theorem~\ref{thm:mtilde_asymp} is given in \ref{subsec:proof_mtilde}.

\subsection{Numerical results for $\tildem_j^\BR$ and $\tildem_j^\WBR$}
\begingroup
\setlength{\aboverulesep}{0pt}
\setlength{\belowrulesep}{0pt}
\renewcommand{\arraystretch}{0.8}
\begin{table}[h!]
	\centering
	\caption{The values of ISB, IV and MISE, multiplied by $10^3$, for the local linear SBF estimator, the PH estimator, and the proposed bias-reduced estimators.}\label{tab:MISE_overall_supp}%
	\small
	\renewcommand{\arraystretch}{1.0}
	\setlength{\tabcolsep}{5pt}
	\begin{tabular}{@{}ccccccccc@{}}
		\toprule
		\rule{0pt}{3.2ex}$(\rho_1,\rho_2)$ & $n$ &
		& $\hatm^{\rmf{LL}}$ & $\hatm_{\rmf{naive}}^{\rmf{PH}}$ & $\tildem^{\BR}$ & $\tildem^{\rmf{WBR}}$ & $\hatm^{\BR}$ & $\hatm^{\rmf{WBR}}$ \\[2pt]
		\midrule
		
		\multirow{6}{*}{(0,3)}
		& \multirow{3}{*}{400}
		& ISB    & 1.89 & 1.04 & 0.73 & 0.70 & 0.68 & 0.67 \\
		& & IV   & 14.00 & 14.60 & 15.16 & 14.93 & 14.80 & 14.70 \\
		& & MISE & 15.89 & 15.64 & 15.89 & 15.64 & 15.48 & 15.36 \\
		
		\cmidrule(lr){2-9}
		
		& \multirow{3}{*}{1000}
		& ISB    & 1.65 & 0.44 & 0.20 & 0.19 & 0.19 & 0.19 \\
		& & IV   & 5.36 & 6.17 & 6.28 & 6.23 & 6.20 & 6.17 \\
		& & MISE & 7.02 & 6.61 & 6.49 & 6.42 & 6.39 & 6.36 \\
		
		\cmidrule(lr){1-9}
		
		\multirow{6}{*}{(2,0)}
		& \multirow{3}{*}{400}
		& ISB    & 3.95 & 1.30 & 0.82 & 0.78 & 0.76 & 0.74 \\
		& & IV   & 13.95 & 14.76 & 15.34 & 15.09 & 14.92 & 14.79 \\
		& & MISE & 17.90 & 16.06 & 16.16 & 15.87 & 15.69 & 15.53 \\
		
		\cmidrule(lr){2-9}
		
		& \multirow{3}{*}{1000}
		& ISB    & 3.72 & 0.47 & 0.28 & 0.28 & 0.29 & 0.29 \\
		& & IV   & 4.46 & 4.77 & 4.85 & 4.80 & 4.76 & 4.74 \\
		& & MISE & 8.18 & 5.24 & 5.13 & 5.08 & 5.04 & 5.03 \\
		
		\cmidrule(lr){1-9}
		
		\multirow{6}{*}{(2,3)}
		& \multirow{3}{*}{400}
		& ISB    & 5.71 & 1.47 & 0.95 & 0.90 & 0.87 & 0.85 \\
		& & IV   & 15.07 & 15.23 & 16.08 & 15.84 & 15.70 & 15.58 \\
		& & MISE & 20.78 & 16.70 & 17.03 & 16.75 & 16.57 & 16.43 \\
		
		\cmidrule(lr){2-9}
		
		& \multirow{3}{*}{1000}
		& ISB    & 5.29 & 0.48 & 0.35 & 0.35 & 0.35 & 0.36 \\
		& & IV   & 5.34 & 6.39 & 6.42 & 6.36 & 6.33 & 6.31 \\
		& & MISE & 10.63 & 6.87 & 6.77 & 6.71 & 6.68 & 6.67 \\
		\bottomrule
	\end{tabular}
\end{table}

\begin{table}[h!]
	\centering
	\caption{The values of ISB$_1$, IV$_1$ and MISE$_1$, multiplied by $10^3$, for the local linear SBF estimator, the PH estimator, and the proposed bias-reduced estimators.}\label{tab:MISE_first_comp_supp}%
	\small
	\renewcommand{\arraystretch}{1.0}
	\setlength{\tabcolsep}{5pt}
	\begin{tabular}{@{}ccccccccc@{}}
		\toprule
		\rule{0pt}{3.2ex}$(\rho_1,\rho_2)$ & $n$ &
		& $\hatm_1^{\rmf{LL}}$ & $\hatm_{1, \rmf{naive}}^{\rmf{PH}}$ & $\tildem_1^{\BR}$ & $\tildem_1^{\rmf{WBR}}$ & $\hatm_1^{\BR}$ & $\hatm_1^{\rmf{WBR}}$ \\[2pt]
		\midrule
		
		\multirow{6}{*}{(0,3)}
		& \multirow{3}{*}{400}
		& ISB$_1$    & 1.33 & 0.18 & 0.34 & 0.33 & 0.31 & 0.30 \\
		& & IV$_1$   & 3.62 & 3.32 & 4.52 & 4.45 & 4.37 & 4.32 \\
		& & MISE$_1$ & 4.94 & 3.50 & 4.86 & 4.77 & 4.68 & 4.62 \\
		
		\cmidrule(lr){2-9}
		
		& \multirow{3}{*}{1000}
		& ISB$_1$    & 0.49 & 0.09 & 0.11 & 0.11 & 0.10 & 0.10 \\
		& & IV$_1$   & 1.72 & 1.53 & 2.10 & 2.10 & 2.08 & 2.07 \\
		& & MISE$_1$ & 2.21 & 1.62 & 2.21 & 2.21 & 2.18 & 2.16 \\
		
		\cmidrule(lr){1-9}
		
		\multirow{6}{*}{(2,0)}
		& \multirow{3}{*}{400}
		& ISB$_1$    & 3.68 & 0.83 & 0.62 & 0.60 & 0.57 & 0.55 \\
		& & IV$_1$   & 5.75 & 5.80 & 6.28 & 6.23 & 6.11 & 6.04 \\
		& & MISE$_1$ & 9.44 & 6.63 & 6.90 & 6.84 & 6.68 & 6.59 \\
		
		\cmidrule(lr){2-9}
		
		& \multirow{3}{*}{1000}
		& ISB$_1$    & 3.54 & 0.40 & 0.27 & 0.25 & 0.25 & 0.25 \\
		& & IV$_1$   & 2.91 & 3.10 & 3.20 & 3.17 & 3.15 & 3.15 \\
		& & MISE$_1$ & 6.45 & 3.50 & 3.47 & 3.42 & 3.39 & 3.41 \\
		
		\cmidrule(lr){1-9}
		
		\multirow{6}{*}{(2,3)}
		& \multirow{3}{*}{400}
		& ISB$_1$      & 3.80 & 0.69 & 0.66 & 0.62 & 0.59 & 0.58 \\
		& & IV$_1$     & 6.83 & 7.45 & 7.52 & 7.44 & 7.32 & 7.25 \\
		& & MISE$_1$   & 10.62 & 8.14 & 8.18 & 8.05 & 7.91 & 7.83 \\
		
		\cmidrule(lr){2-9}
		
		& \multirow{3}{*}{1000}
		& ISB$_1$    & 3.60 & 0.33 & 0.23 & 0.23 & 0.23 & 0.23 \\
		& & IV$_1$   & 2.99 & 3.26 & 3.35 & 3.35 & 3.33 & 3.31 \\
		& & MISE$_1$ & 6.59 & 3.59 & 3.58 & 3.58 & 3.56 & 3.54 \\
		\bottomrule
	\end{tabular}
\end{table}

\begin{table}[h!]
	\centering
	\caption{The values of ISB$_2$, IV$_2$ and MISE$_2$, multiplied by $10^3$, for the local linear SBF estimator, the PH estimator, and the proposed bias-reduced estimators.}\label{tab:MISE_second_comp_supp}%
	\small
	\renewcommand{\arraystretch}{1.0}
	\setlength{\tabcolsep}{5pt}
	\begin{tabular}{@{}ccccccccc@{}}
		\toprule
		\rule{0pt}{3.2ex}$(\rho_1,\rho_2)$ & $n$ &
		& $\hatm_2^{\rmf{LL}}$ & $\hatm_{2, \rmf{naive}}^{\rmf{PH}}$ & $\tildem_2^{\BR}$ & $\tildem_2^{\rmf{WBR}}$ & $\hatm_2^{\BR}$ & $\hatm_2^{\rmf{WBR}}$ \\[2pt]
		\midrule
		
		\multirow{6}{*}{(0,3)}
		& \multirow{3}{*}{400}
		& ISB$_2$    & 1.66 & 1.47 & 0.29 & 0.29 & 0.28 & 0.28 \\
		& & IV$_2$   & 6.16 & 5.66 & 6.92 & 6.86 & 6.79 & 6.75 \\
		& & MISE$_2$ & 7.81 & 7.13 & 7.21 & 7.14 & 7.07 & 7.03 \\
		
		\cmidrule(lr){2-9}
		
		& \multirow{3}{*}{1000}
		& ISB$_2$    & 1.48 & 1.32 & 0.11 & 0.10 & 0.10 & 0.10 \\
		& & IV$_2$   & 2.35 & 2.14 & 2.74 & 2.72 & 2.70 & 2.69 \\
		& & MISE$_2$ & 3.83 & 3.47 & 2.85 & 2.82 & 2.80 & 2.78 \\
		
		\cmidrule(lr){1-9}
		
		\multirow{6}{*}{(2,0)}
		& \multirow{3}{*}{400}
		& ISB$_2$    & 1.03 & 0.14 & 0.22 & 0.21 & 0.20 & 0.19 \\
		& & IV$_2$   & 4.80 & 4.42 & 5.40 & 5.35 & 5.28 & 5.24 \\
		& & MISE$_2$ & 5.83 & 4.57 & 5.62 & 5.57 & 5.48 & 5.43 \\
		
		\cmidrule(lr){2-9}
		
		& \multirow{3}{*}{1000}
		& ISB$_2$    & 0.50 & 0.10 & 0.10 & 0.10 & 0.09 & 0.08 \\
		& & IV$_2$   & 1.96 & 2.02 & 2.31 & 2.28 & 2.26 & 2.25 \\
		& & MISE$_2$ & 2.47 & 2.12 & 2.41 & 2.38 & 2.36 & 2.33 \\
		
		\cmidrule(lr){1-9}
		
		\multirow{6}{*}{(2,3)}
		& \multirow{3}{*}{400}
		& ISB$_2$    & 1.86 & 0.46 & 0.37 & 0.36 & 0.35 & 0.35 \\
		& & IV$_2$   & 5.93 & 6.62 & 6.80 & 6.74 & 6.66 & 6.62 \\
		& & MISE$_2$ & 7.79 & 7.08 & 7.17 & 7.09 & 7.01 & 6.97 \\
		
		\cmidrule(lr){2-9}
		
		& \multirow{3}{*}{1000}
		& ISB$_2$    & 1.70 & 0.16 & 0.14 & 0.14 & 0.13 & 0.13 \\
		& & IV$_2$   & 2.24 & 2.87 & 2.87 & 2.86 & 2.84 & 2.84 \\
		& & MISE$_2$ & 3.94 & 3.03 & 3.01 & 3.01 & 2.97 & 2.97 \\
		
		\bottomrule
	\end{tabular}
\end{table}
\endgroup

Tables~\ref{tab:MISE_overall_supp}--\ref{tab:MISE_second_comp_supp} report the ISB, IV, and MISE values for the simulation design described in Section~\ref{sec:Numerical_Study}, including results for $\tildem_j^\BR$ and $\tildem_j^\WBR$.
NOtably, these tables extend the corresponding simulation results reported in Tables~\ref{tab:MISE_overall}--\ref{tab:MISE_second_comp}.
Tables~\ref{tab:ASPE_real_1_supp} and \ref{tab:ASPE_real_2_supp} present the 5-fold cross-validated ASPE values for the estimators applied to the bikeshare and flights data, respectively, introduced in Section~\ref{sec:real_data}.

\begin{table}[h]
	\centering
	\caption{Comparison of 5-fold cross-validated ASPE values in the analysis of the Bikeshare dataset.}\label{tab:ASPE_real_1_supp}%
	\renewcommand{\arraystretch}{0.95}
	\small
	\begin{tabular}{@{}lcccccc@{}}
		\toprule
		Estimator & $\hatm^{\rmf{LL}}$ & $\hatm^{\rmf{PH}}$	& $\tildem^{\BR}$ & $\tildem^{\rmf{WBR}}$ & $\hatm^{\BR}$ & $\hatm^{\rmf{WBR}}$ \\ \midrule
		ASPE & 1.7630 & 1.7628 & 1.7459 & 1.7457 & 1.7453 & 1.7451 \\ \bottomrule
	\end{tabular}
\end{table}

\begin{table}[h!]
	\centering
	\caption{Comparison of 5-fold cross-validated ASPE values in the analysis of the flights dataset.}\label{tab:ASPE_real_2_supp}%
	\renewcommand{\arraystretch}{0.95}
	\small
	\begin{tabular}{@{}lcccccc@{}}
		\toprule
		Estimator & $\hatm^{\rmf{LL}}$ & $\hatm^{\rmf{PH}}$ & $\tildem^{\BR}$ & $\tildem^{\rmf{WBR}}$ & $\hatm^{\BR}$ & $\hatm^{\rmf{WBR}}$ \\ \midrule
		ASPE & 1.5500 & 1.5488 & 1.5483 & 1.5488 & 1.5487 & 1.5478 \\ \bottomrule
	\end{tabular}
\end{table}

\clearpage
\section{Proof of main results}\label{suppsec:proof_theorems}

\subsection{Proof of Proposition~\ref{prop:bias_expansion}}

We first introduce some notation to facilitate the proof of Proposition~\ref{prop:bias_expansion}.
Let $\cM$ be the space of $(d+1)$-tuples of maps $\bmf^{\F,\tp} = (\bmf^0, \bmf^{1,1}, \dots, \bmf^{1,d})$, where $\bmf^0, \bmf^{1,j}:[0,1]^d \to \bbR$ satisfy $\int_{[0,1]^d}\bmf^0(\bmx)^2 \d\bmx < \infty$ and $\int_{[0,1]^d} \bmf^{1,j}(\bmx)^2\d\bmx < \infty$ for all $1 \le j \le d$.
Let $\mathcal{M_{\text{add}}}$ be a subspace of $\cM$ consisting of functions that satisfy $\bmf^0(\bmx) = f_1(x_1) + \cdots +f_d(x_d)$ and $\bmf^{1,j}(\bmx) = f_{1,j}(x_j)$ for some univariate functions $f_j$ and $f_{1,j}$.
For each $1\le j\le d$, we define the subspace $\cM_j \subset \cM_{\text{add}}$ by
\begin{align*}
	\cM_j := \Big\{&\bmf^{\F,\tp}\in \cM: \bmf^{1, k}=\mbf{0} \text{ for } 1 \le k \neq j \le d  \text{ with }\\
	&\bmf^0(\bmx) = f_j(x_j) \text{ and } \bmf^{1, j}(\bmx) = f_{1,j}(x_j) \text{ for some univariate maps } f_j \text{ and } f_{1,j} \Big\}.
\end{align*}
For $\bmu \in \bbR^d$, let $\mathfrak{vec}(\bmu) = (1, u_1, u_2, \dots, u_d)^{\top} \in \bbR^{d+1}$.
Let $\mathfrak{N}^\F$ be the $(d+1)\times(d+1)$ positive definite matrix given by
\begin{equation*}
	\mathfrak{N}^\F := \int_{[-1,1]^d}\mathfrak{vec}(\bmu)\cdot \mathfrak{vec}(\bmu)^{\top} \prod_{j=1}^dK(u_j) \d\bmu.
\end{equation*}
We define $\|\cdot\|_2: \cM \to [0, \infty)$ by $\|\bmf^{\F,\tp}\|_2^2 = \langle \bmf^{\F,\tp}, \bmf^{\F,\tp} \rangle$, where
\begin{equation*}
	\langle \bmf^{\F,\tp}, \bmg^{\F,\tp} \rangle :=
	\int_{[0,1]^{d}} \Big\langle \bmf^{\F,\tp}(\bmx), \mathfrak{N}^\F \cdot \bmg^{\F,\tp}(\bmx) \Big\rangle_{\bbR^{d+1}}
	p(\bmx) \d\bmx.
\end{equation*}
Under mild conditions, $\|\cdot \|_2$ forms a norm on $\cM$ \citep{jeon2022locally}.

Define $\hat{\mathfrak{M}}^\F:[0,1]^d \to \bbR^{(d+1)\times(d+1)}$ by
\begin{equation*}
	\hat{\mathfrak{M}}^{\F}(\bmx) := \meanin \mathfrak{vec}\big((\bmX_i-\bmx) / \bmh\big)  \cdot \mathfrak{vec}\big((\bmX_i-\bmx) / \bmh\big)^{\top}  \cdot    \prod_{\ell=1}^d K_{h_\ell}(x_{\ell}, X_{i\ell}),
\end{equation*}
where $\bmu/\bmh = (u_1/h_1, \dots, u_d/h_d)^\top$.
Let $\hat{\pi}_j : \cM \to \cM_j$ be the projection operator defined by
\begin{equation*}
	\hat{\pi}_j(\bmf^{\F,\tp}) := \frakU_j^{\top} \cdot \hatbmM_{jj}(x_j)^{-1} \cdot	\int_{[0,1]^{d-1}}\frakU_j \cdot \hat{\mathfrak{M}}^{\F}(\bmx) \cdot \bmf^{\F,\tp}(\bmx) \d\bmx_{-j},
\end{equation*}
where $\frakU_j$ denotes the $2 \times (d+1)$ matrix that transform the full tuple into the $j$-th sub-tuple, given by
\begin{equation*}
	\frakU_j  :=
	\begin{pmatrix}
		1 & \mbf{0}_{j-1}^{\top} & 0 & \mbf{0}_{d-j}^{\top} \\
		0 & \mbf{0}_{j-1}^{\top} & 1 & \mbf{0}_{d-j}^{\top}
	\end{pmatrix}.
\end{equation*}
In particular, for $\bmf_k^{\F,\tp} = (\bmf^k, \mbf{0}_{k-1}^{\top}, \bmf^{1,k}, \mbf{0}_{d-k}^{\top})^{\top} \in \cM_k$, let $\bmf_k^{\tp}(x_k) = (f_k, f_{1,k})^{\top}$ with $f_k, f_{1,k}$ satisfying $\bmf^k(\bmx) = f_k(x_k)$ and $\bmf^{1,k}(\bmx) = f_{1,k}(x_k)$.
Then, it holds that
\begin{equation}\label{eq:projection_restricted}
	\hat{\pi}_j \big(\frakU_k^\top \cdot \bmf_k^{\tp} \big)(x_j) = \frakU_j^{\top} \cdot \int_0^1\hatbmM_{jj}(x_j)^{-1} \hatbmM_{jk}(x_j, x_k) \bmf_k^{\tp}(x_k) \d x_k.
\end{equation}

For $1 \le j \le d$, define
\begin{equation*}
	\hat{\Delta}_j^{\tp}(x_j) := \hatbmm_j^{\tp}(x_j) - \bmm_j^{\tp}(x_j) - \tildebmm_j^{A, \tp}(x_j) - \frac{h_j^2}{2}\bmN_j(x_j)^{-1}\bm{\gamma}_j(x_j)m_j''(x_j), \quad x_j\in[0,1].
\end{equation*}

Lemma~\ref{lem:delta_j_k} establishes the asymptotic properties of $\hat\Delta_j^\tp$.
The proof of Lemma~\ref{lem:delta_j_k} is given in Section~\ref{subsec:proof_delta_j_k}.

\begin{lemma}\label{lem:delta_j_k}
	Assume (A1)--(A5).
	Then, for $1\le j\le d$,
	\begin{align*}
		& \sup_{x_j\in[0,1]}\bigg| \hat{\Delta}_j^{\tp}(x_j) + \sumkneqj \int_0^1 \hatbmM_{jj}(x_j)^{-1}\hatbmM_{jk}(x_j, x_k)\hat{\Delta}_k^{\tp}(x_k) \d x_k \bigg| \\
		&\hspace{1cm} = o_p\bigg(\sum_{m=1}^d h_m^2\bigg) + O_p\bigg(\sqrt{\frac{\log n}{n}} \bigg).
	\end{align*}
	If we further assume that each $m_j''$ is Lipschitz continuous on $[0,1]$, then for $1\le j\le d$,
	\begin{align*}
		\sup_{x_j\in[0,1]}\bigg| \hat{\Delta}_j^{\tp}(x_j) + \sumkneqj \int_0^1 \hatbmM_{jj}(x_j)^{-1}\hatbmM_{jk}(x_j, x_k)\hat{\Delta}_k^{\tp}(x_k) \d x_k \bigg| = O_p\bigg(\sum_{m=1}^d h_m^3 + \sqrt{\frac{\log n}{n}} \bigg).
	\end{align*}
\end{lemma}

\begin{proof}[Proof of Proposition~{\upshape\ref{prop:bias_expansion}}]
	Let $\hat{\Delta}_j^{\F,\tp} := \frakU_j^{\top} \cdot \hat{\Delta}_j^{\tp} \in \cM_j$ and $\hat{\Delta}^{\F,\tp} := \sumjd \hat{\Delta}_j^{\F,\tp} \in \cM$.
	Define a linear operator $\hat{T}: \cM_{\text{add}} \to \cM_{\text{add}}$ by $\hat{T} = (I-\hat{\pi}_d)\circ(I-\hat{\pi}_{d-1})\circ \cdots \circ(I-\hat{\pi}_1)$ where $I : \cM \to \cM$ is the identity map.
	Combining \eqref{eq:projection_restricted} and Lemma~\ref{lem:delta_j_k} yields
	\begin{align*}
		&\hspace{-0.5cm}\sup_{x_j\in[0,1]}\Big|\hat{\pi}_j\big(\hat{\Delta}^{\F,\tp}\big)(x_j) \Big| \\
		&= \sup_{x_j\in[0,1]} \Bigg| \, \frakU_j^{\top}	\cdot	\bigg(\hat{\Delta}_j^{\tp}(x_j)	+ \sumkneqj \int_0^1\hatbmM_{jj}(x_j)^{-1} \hatbmM_{jk}(x_j, x_k) \hat{\Delta}_k^{\tp}(x_k) \d x_k \bigg) \, \Bigg| \\
		&=O_p\bigg(\sum_{m=1}^d h_m^3 + \sqrt{\frac{\log n}{n}} \bigg).
	\end{align*}	
	Using the above result and the definition of the operator $\hat{\pi}_j$, we obtain
	\begin{equation}\label{eq:Delta_prop}
		\sup_{\bmx\in[0,1]^d} \Big| \hat{\Delta}^{\F,\tp}(\bmx) - \hat{T}(\hat{\Delta}^{\F,\tp})(\bmx) \Big| = O_p\bigg(\sum_{m=1}^d h_m^3 + \sqrt{\frac{\log n}{n}} \bigg).
	\end{equation}

	By (S.4) in \cite{jeon2022locally}, there exists a constant $0<\gamma<1$ such that
	\begin{equation*}
		\Prob\Big(\big\| \hat{T} \big\|_{\mathcal{L}(\cM_{\rmf{add}},\,\cM_{\rmf{add}})} < \gamma \Big) \to 1 \quad \text{as} \quad n\to\infty,
	\end{equation*}
	where
	\begin{align*}
		\big\| \hat{T} \big\|_{\mathcal{L}(\cM_{\text{add}}, \cM_{\text{add}})}	:=
		\sup \Big\{ \big\| \hat{T}(\bmf^{\F,\tp}) \big\|_2:
		\bmf^{\F,\tp} \in \cM_{\text{add}}, \big\| \bmf^{\F,\tp} \big\|_2 = 1
		\Big\}
	\end{align*}
	Applying this to \eqref{eq:Delta_prop}, we obtain
	\begin{align*}
		\big\| \hat{\Delta}^{\F,\tp} \big\|_2 =O_p\bigg(\sum_{m=1}^d h_m^3 + \sqrt{\frac{\log n}{n}} \bigg).
	\end{align*}
	Following the proof of Theorem 2 in \cite{jeon2022locally}, we have
	\begin{align*}
		\sup_{x_j\in[0,1]} \big| \hat{\Delta}_j^{\tp}(x_j) \big| = O_p\bigg(\sum_{m=1}^d h_m^3 + \sqrt{\frac{\log n}{n}} \bigg).
	\end{align*}
	This completes the proof of Proposition~\ref{prop:bias_expansion}.	
\end{proof}

\subsection{Proof of Theorem~\ref{thm:mhat_asymp}}\label{subsec:proof_mhat_asymp}

For $\ell = 1,\dots, L$, let $\hatbmM_{jj}(x_j;b_{j,\ell})$, $\tildebmm_j^{A, \tp}(x_j;b_{j,\ell})$, $\bmN_j(x_j;b_{j,\ell})$ and $\bmgamma_j(x_j;b_{j,\ell})$ denote the quantities obtained from $\hatbmM_{jj}(x_j)$, $\tildebmm_j^{A, \tp}(x_j)$, $\bmN_j(x_j)$ and $\bmgamma_j(x_j)$ by replacing $h_j$ with $b_{j,\ell}$.
Specifically,
\begin{equation}\label{h_lj_notation}
	\begin{split}
		\hatbmM_{jj}(x_j;b_{j,\ell}) :=& \,  \meanin \bmvec\bigg(\frac{X_{ij}-x_j}{b_{j,\ell}}\bigg)
		\bmvec\bigg(\frac{X_{ij}-x_j}{b_{j,\ell}}\bigg)^\top K_{b_{j,\ell}}(x_j,X_{ij}),   \\
		\tildebmm_j^{A, \tp}(x_j;b_{j,\ell}) :=& \,  \meanin \hatbmM_{jj}(x_j;b_{j,\ell})^{-1} \bmvec\bigg(\frac{X_{ij} - x_j}{b_{j,\ell}} \bigg) K_{b_{j,\ell}}(x_j, X_{ij})\cdot \varepsilon_i, \\
		\bmN_j(x_j;b_{j,\ell}) :=& \, \int_0^1 \bmvec\bigg(\frac{v-x_j}{b_{j,\ell}} \bigg)   \bmvec\bigg(\frac{v-x_j}{b_{j,\ell}} \bigg)^{\top} K_{b_{j,\ell}}(x_j, v)\d v, \\
		\bmgamma_j(x_j;b_{j,\ell}) :=& \, \int_0^1 \bigg( \frac{v-x_j}{b_{j,\ell}} \bigg)^2 \bmvec\bigg(\frac{v-x_j}{b_{j,\ell}} \bigg) K_{b_{j,\ell}}(x_j, v) \d v.
	\end{split}	
\end{equation}
Since $b_j = \max_{1 \le \ell \le L}b_{j,\ell}$, the interval $I_j = [2b_j, 1-2b_j]$ lies within the interior for every bandwidth $b_{j,\ell}$, i.e., $I_j=[2b_j, 1-2b_j] \subset [2b_{j,\ell}, 1-2b_{j,\ell}]$.
Hence, for all $x_j \in I_j$ and all $\ell = 1, \dots, L$, it holds that $\bmN_j(x_j;b_{j,\ell}) = \bmN$ and $\bmgamma_j(x_j;b_{j,\ell}) = \bmgamma$ for some $\bmN \in \bbR^{2\times 2}$ and $\bmgamma \in \bbR^2$.

Since the proofs for $\tildem_j^\BR$ and $\tildem_j^\WBR$ are similar, we only establish the convergence rates for $\tildem_j^\BR$.

For $\ell=1,\dots,L$ and $j=1,\dots,d$, let
\begin{equation*}
	d_{j,\ell} := \frac{L \cdot b_{j,\ell}^2 - \big( \sum_{\ell'=1}^L h_{\ell' j}^2 \big)}{L \cdot \sum_{\ell'=1}^L h_{\ell' j}^4 - \big( \sum_{\ell'=1}^L h_{\ell' j}^2 \big)^2}.
\end{equation*}
Then, it holds that
\begin{equation*}
	\sumellL d_{j,\ell} = 0, \quad \sumellL d_{j,\ell} \cdot b_{j,\ell}^2 = 1, \quad d_{j,\ell} = O\big(b_j^{-2}\big).
\end{equation*}
Recall the definitions of $\tildebmm_j^{A,\tp}(x_j;b_{j,\ell})$, $\bmN_j(x_j;b_{j,\ell})$ and $\bmgamma_j(x_j;b_{j,\ell})$ at \eqref{h_lj_notation}.
By Proposition~\ref{prop:bias_expansion}, we obtain the uniform asymptotic expansion for $\hatm_j^\BR$:
\begin{equation}\label{eq:mhat_expansion}
	\begin{split}
		\hatm_j^{\BR}(x_j) &= m_j(x_j) + \fraku_0^{\top}\tildebmm_j^{A, \tp}(x_j;h_j) + \frac{1}{2}h_j^2 \cdot \fraku_0^{\top}\bmN_j(x_j;h_j)^{-1}\bmgamma_j(x_j;h_j) m_j''(x_j) \\
		&\quad - h_j^2 \cdot \sumellL d_{j,\ell} \cdot \frac{1}{2}b_{j,\ell}^2 \cdot \fraku_0^\top \bmN_j(x_j;b_{j,\ell})^{-1} \bmgamma_j(x_j;b_{j,\ell}) m_j''(x_j) \\
		&\quad - h_j^2 \cdot \sumellL d_{j,\ell} \cdot \fraku_0^\top \tildebmm_j^{A, \tp}(x_j;b_{j,\ell}) \\
		&\quad + O_p\Bigg(h_j^3 + h_{-j}^3 + \sqrt{\frac{\log n}{n}}
		+ \frac{h_j^2}{b_j^2}\bigg(b_j^3 + b_{-j}^3 + \sqrt{\frac{\log n}{n}}\bigg)\Bigg),
	\end{split}
\end{equation}
where the $O_p$ term is uniform in $x_j\in[0,1]$.

For $1\le j\le d$ and $x_j\in[0,1]$, let
\begin{align*}
	\text{bias}_j(x_j) :=& \, \frac{1}{2}h_j^2 \cdot \fraku_0^{\top}\bmN_j(x_j;h_j)^{-1}\bmgamma_j(x_j;h_j) m_j''(x_j) \\
	&\quad - h_j^2 \cdot \sumellL d_{j,\ell} \cdot \frac{1}{2}b_{j,\ell}^2 \cdot \fraku_0^\top \bmN_j(x_j;b_{j,\ell})^{-1} \bmgamma_j(x_j;b_{j,\ell}) m_j''(x_j).
\end{align*}
Since $h_j = o(b_j)$, it holds that $I_j = [2b_j,1-2b_j]\subset [2h_j,1-2h_j]$ for sufficiently large $n$.
Then, similar to \eqref{BR_h_lj_prop}, we have
\begin{equation}\label{BR_h_j_prop}
	\sup_{x_j \in I_j} \big| \text{bias}_j(x_j) \big| = 0, \quad \sup_{x_j \in [0,1]} \big| \text{bias}_j(x_j) \big| = O(h_{j}^2).
\end{equation}
Furthermore, analogous to \eqref{mtilde_A_lj_prop}, it holds that
\begin{equation}\label{mtilde_A_j_prop}
	\begin{split}
		|\fraku_0^\top \tildebmm_j^{A, \tp}(x_j;h_{j})| &=  O_p(n^{-1/2}h_{j}^{-1/2} ), \quad x_j \in (0,1), \\
		\int_0^1 |\fraku_0^\top \tildebmm_j^{A, \tp}(x_j;h_{j})|^2 \d x_j &= O_p(n^{-1}h_{j}^{-1}), \\
		\sup_{x_j\in[0,1]} |\fraku_0^\top \tildebmm_j^{A, \tp}(x_j;h_{j})| &= O_p \Bigg(\sqrt{\frac{\log n}{nh_{j}}}  \Bigg).
	\end{split}
\end{equation}
Moreover, by (A6) and the assumption that $h_j = o(b_j)$, we have
\begin{equation}\label{BR_residual_rate}
	h_j^3 + h_{-j}^3 +\sqrt{\frac{\log n}{n}} + \frac{h_j^2}{b_j^2} \cdot \bigg(b_j^3 + b_{-j}^3 + \sqrt{\frac{\log n}{n}}\bigg) \lesssim h_j^2b_j + r_{nj} + \sqrt{\frac{\log n}{n}}.
\end{equation}
Substituting \eqref{mtilde_A_lj_prop} and \eqref{BR_h_j_prop}--\eqref{BR_residual_rate} into \eqref{eq:mhat_expansion}, we obtain the pointwise, $L_2$ and uniform convergence rates for $\hatm_j^\BR$.
\qed

\subsection{Proof of Theorem~\ref{thm:unif_CLT}}

To establish Theorem~\ref{thm:unif_CLT}, we first approximate the kernel-weighted sum associated with the centered estimator $\hat m_j^{\#}-m_j$ by a corresponding sum involving normal random variables.
In this proof, we denote $\#\in\{\BR,\WBR\}$.
We then approximate this kernel-weighted sum of normal random variables by a stationary Gaussian process with similar covariance structure.
The following two lemmas formalize these approximation steps.
Their proofs are deferred to
Sections~\ref{subsec:proof_normal_rv_approx} and \ref{subsec:proof_gaussian_process_approx},
respectively.
For any closed interval $I\subset(0,1)$, let
\begin{equation}\label{eq:TAB_nj_def}
	T_{nj}(I) := \frac{|I|}{h_j}, \quad A_{nj}(I) := \sqrt{2\log T_{nj}(I)}, \quad B_{nj}(I) := A_{nj}(I) + \frac{\log(\lambda_K/(4\pi^2))}{2A_{nj}(I)},
\end{equation}
where $\lambda_K := \nu_0^{-1} \cdot \int \{K'(t)\}^2\d t$.
Define
\begin{equation}\label{eq:b_nj_def}
	b_{nj}(u;I) := B_{nj}(I) + \frac{u}{A_{nj}(I)}, \quad u\in\bbR.
\end{equation}

\begin{lemma}\label{lem:normal_rv_approx}
	Assume that the conditions in Theorem~\ref{thm:unif_CLT} hold.
	Then, for each $1\le j\le d$, there exist random variables $Z_{1j},\ldots,Z_{nj}$ such that
	\begin{equation*}
		Z_{ij}\mid X_{ij} \sim N\{0,\sigma_j^2(X_{ij})\}, \qquad i=1,\ldots,n,
	\end{equation*}
	$Z_{1j},\ldots,Z_{nj}$ are independent conditionally on $X_{1j},\ldots,X_{nj}$, and for any closed interval $I\subset (0,1)$ and $\#\in\{\BR,\WBR\}$, it holds that
	\begin{align*}
		&\sup_{x_j\in I}	\Bigg|	\frac{\sqrt{n h_j p_j(x_j)}}{\sigma_j(x_j)}	\{\hat m_j^{\#}(x_j)-m_j(x_j)\}	-	\frac{1}{\sqrt{n h_j p_j(x_j)}\,\sigma_j(x_j)}		\sumin	K\bigg(\frac{x_j-X_{ij}}{h_j}\bigg)Z_{ij}	\Bigg| \\
		&\hspace{2cm} = o_p\big(A_{nj}(I)^{-1}\big).
	\end{align*}
\end{lemma}

Define the scaled sum of random variables
\begin{equation}\label{eq:Gn_def_hash}
	G_{nj}(x_j) := 	\frac{1}{\sqrt{n h_j p_j(x_j)}\,\sigma_j(x_j)}	\sum_{i=1}^n	K\bigg(\frac{x_j-X_{ij}}{h_j}\bigg)Z_{ij}, \qquad x_j\in [0,1],
\end{equation}
where $Z_{1j},\ldots,Z_{nj}$ are the normal random variables constructed in Lemma~\ref{lem:normal_rv_approx}.

\begin{lemma}\label{lem:gaussian_process_approx}
	Assume that the conditions in Theorem~\ref{thm:unif_CLT} hold.
	For any closed interval $I\subset(0,1)$, let $W_{nj}\equiv W_{nj}(\cdot;I)$ be a centered stationary Gaussian process on $I$ with
	covariance function
	\begin{equation*}
		\Cov\{W_{nj}(x_j),W_{nj}(x_j')\} = (K*K)\bigg(\frac{x_j-x_j'}{h_j}\bigg), \qquad x_j,x_j'\in I,		
	\end{equation*}
	where $(K*K)(t):=\int K(u)K(u+t)\,\d u$.
	Then, it holds that
	\begin{equation*}
		\sup_{t\in\bbR}	\Bigg| \Prob\bigg( \sup_{x_j\in I}|G_{nj}(x_j)|\le t
		\,\bigg| \, X_{1j},\ldots,X_{nj} \bigg) - \Prob\bigg( \sup_{x_j\in I}|W_{nj}(x_j)|\le t \bigg) \Bigg| = o_p(1).
	\end{equation*}
\end{lemma}

\begin{proof}[Proof of Theorem~{\upshape\ref{thm:unif_CLT}}]
	Since the proofs for $\hatm_j^\BR$ and $\hatm_j^\WBR$ are identical, we present only the proof for $\hatm_j^\BR$.	
	Given the closed interval $I = [a_I,b_I] \subset (0,1)$, the assumption of $T_{nj}$ at \eqref{eq:TAB_nj_def} implies that $T_{nj}(I) = (b_I - a_I) / h_j \to \infty$ as $n\to\infty$.
	Define the scaled error
	\begin{equation*}
		S_{nj}^{\BR}(x_j) := \frac{\sqrt{n h_j p_j(x_j)}}{\sigma_j(x_j)}	\{\hat m_j^{\BR}(x_j)-m_j(x_j)\},	\quad x_j\in I.
	\end{equation*}
	Then, we can write
	\begin{equation*}
		\hat\Delta_{nj}^{\BR}(I) =	\nu_0^{-1/2} \cdot \sup_{x_j\in I}|S_{nj}^{\BR}(x_j)|.
	\end{equation*}
	By Lemma~\ref{lem:normal_rv_approx}, it holds that
	\begin{equation*}
		\sup_{x_j\in I} \big|S_{nj}^{\BR}(x_j)-G_{nj}(x_j)\big| = o_p(A_{nj}(I)^{-1}),	
	\end{equation*}
	where $G_{nj}$ is defined at \eqref{eq:Gn_def_hash}.
	Hence, with
	\begin{equation*}
		M_{nj}(I) := \nu_0^{-1/2} \cdot \sup_{x_j \in I} |G_{nj}(x_j)|,
	\end{equation*}
	we have
	\begin{equation}\label{eq:Delta_Gn_close_hash}
		R_{nj}(I) := \big| \hat\Delta_{nj}^{\BR}(I)-M_{nj}(I) \big| \le	\nu_0^{-1/2} \cdot \sup_{x_j\in I} \big|S_{nj}^{\BR}(x)-G_{nj}(x_j)\big| = o_p(A_{nj}(I)^{-1}).
	\end{equation}
	Furthermore, by the definition of $b_{nj}$ at \eqref{eq:b_nj_def}, for every fixed $\delta>0$, the following inequalities hold uniformly over $t\in\bbR$:
	\begin{equation}\label{eq:sandwich_Delta_Gn_hash}
		\begin{split}
			\Prob\Big( M_{nj}(I)\le b_{nj}(t-\delta;I) \Big) &\le \Prob\Big(\hat\Delta_{nj}^\BR(I) \le b_{nj}(t;I)\Big) + \Prob\big(A_{nj}(I)R_{nj}(I) > \delta\big) \\
			&\le \Prob\Big( M_{nj}(I)\le b_{nj}(t+\delta;I) \Big) + 2\Prob\big(A_{nj}(I)R_{nj}(I) > \delta\big).		
		\end{split}
	\end{equation}
	This establishes a distributional approximation of $\hat\Delta_n^{\BR}(I)$ by $M_{nj}(I)$.

	Next, we approximate the distribution of $M_{nj}$ using the Gaussian process $W_{nj}$ defined in Lemma~\ref{lem:gaussian_process_approx}.
	Let
	\begin{equation*}
		U_{nj}(I) := \nu_0^{-1/2} \cdot \sup_{x_j\in I}|W_{nj}(x_j)|.
	\end{equation*}
	By combining Lemma~\ref{lem:gaussian_process_approx} with the fact that the conditional probability is bounded by $1$ almost surely, an application of the bounded convergence theorem gives
	\begin{equation*}
		\sup_{t\in\bbR}	\Bigg| \Prob\bigg( \sup_{x_j\in I}|G_{nj}(x_j)|\le t \bigg) - \Prob\bigg( \sup_{x_j\in I}|W_{nj}(x_j)|\le t \bigg) \Bigg| = o(1),
	\end{equation*}
	and equivalently,
	\begin{equation*}
		\sup_{t\in\bbR} \Big| \Prob\big(M_{nj}(I)\le t\big) - \Prob\big(U_{nj}(I)\le t\big) \Big| = o(1).
	\end{equation*}
	In particular, it holds that
	\begin{equation}\label{eq:Gn_Wn_at_bn_hash}
		\sup_{u\in\bbR} \bigg| \Prob\big(M_{nj}(I)\le b_{nj}(u;I) \big) - \Prob\big(U_{nj}(I)\le b_{nj}(u;I)\big) \bigg| = o(1).	
	\end{equation}

	It remains to obtain the extreme-value limit for $U_{nj}(I)$.
	Recall that $T_{nj}(I) = (b_I - a_I)/h_j$ for the interval $I=[a_I,b_I]$.
	Define the rescaled Gaussian process
	\begin{equation*}
		\bbW_{nj}(t) := \nu_0^{-1/2} \cdot W_{nj}(a_I + h_j t), \quad 0\le t\le T_{nj}(I).
	\end{equation*}
	Then, we can write
	\begin{equation*}
		U_{nj}(I) =	\sup_{0\le t\le T_{nj}(I)}|\bbW_{nj}(t)|.
	\end{equation*}
	Note that $\bbW_{nj}$ is a centered stationary Gaussian process with covariance function $\rho := \nu_0^{-1} \cdot (K * K)$; that is, $\E\big(\bbW_{nj}(s)\bbW_{nj}(t)\big) = \rho(s-t)$ for all $0\le s,t\le T_{nj}(I)$.
	By (A2), the function $\rho:\bbR\to\bbR$ satisfies $\rho(0)=1$, $0\le \rho(t) < 1$ for $t\neq 0$, $\rho(t)\log t \to 0$ as $t\to\infty$, and
	\begin{equation*}
		\rho(t) = 1 - \frac{\lambda_K t^2}{2}  + o(t^2), \quad t\to0
	\end{equation*}
	with $\lambda_K>0$ defined at \eqref{eq:b_nj_def}.
	Here, we have used the fact that
	\begin{equation*}
		-\rho''(0) = - \frac{(K*K)''(0)}{\nu_0} = \frac{\int \{K'(u)\}^2\d u}{\nu_0} = \lambda_K,
	\end{equation*}
	where the second equality follows from integration by parts.
	Thus, applying Theorem 8.2.7 of \cite{leadbetter1983extremes} for stationary Gaussian processes, we obtain
	\begin{equation*}
		\Prob\Bigg( A_{nj}(I) \bigg\{ \sup_{0\le t\le T_{nj}(I)}|\bbW_{nj}(t)|-	B_{nj}(I) \bigg\} \le u \Bigg) 	\to	\exp(-2e^{-u}), \quad u\in\bbR.
	\end{equation*}
	In other words,
	\begin{equation}\label{eq:MW_EV_limit_hash}
		\Prob\big(U_{nj}(I)\le b_{nj}(u;I) \big) \to \exp(-2e^{-u}), \quad u\in\bbR.
	\end{equation}

	Combining \eqref{eq:sandwich_Delta_Gn_hash}--\eqref{eq:MW_EV_limit_hash}, for every fixed $\delta>0$ and $u\in\bbR$, we have
	\begin{equation*}
		\exp\big(-2e^{-(u-\delta)}\big) \le \liminf_{n\to\infty}	\Prob\big(\hat\Delta_{nj}^{\BR}(I)\le b_{nj}(u;I)\big)
	\end{equation*}
	and
	\begin{equation*}
		\limsup_{n\to\infty} \Prob\big(\hat\Delta_{nj}^{\BR}(I)\le b_{nj}(u;I)\big) \le \exp\big(-2e^{-(u+\delta)}\big).
	\end{equation*}
	Letting $\delta\downarrow0$ and using the continuity of $v\mapsto\exp(-2e^{-v})$, we obtain
	\begin{equation*}
		\lim_{n \to \infty} \Prob\big(\hat\Delta_{nj}^{\BR}(I)\le b_{nj}(u;I)\big) = \exp\big(-2e^{-u}\big), \quad u\in\bbR.
	\end{equation*}
	By the definition of $b_{nj}$ at \eqref{eq:b_nj_def}, this concludes the proof.
\end{proof}

\subsection{Proof of Corollary~\ref{cor:unif_CLT}}

We first establish uniform consistency of the estimators of $p_j$ and
$\sigma_j^2$ on $I$.
By the condition (B1), the kernel density estimator satisfies
\begin{equation}\label{eq:p_unif}
	\sup_{x_j\in I}\abs{\hatp_j(x_j) - p_j(x_j)} = O_p\bigg(\sqrt{\frac{\log n}{nh_j}}\bigg) + o_p(h_j)
\end{equation}
Note that the estimate $\hat\sigma_j^2$ of $\sigma_j^2:=\E(\varepsilon^2|X_j=\cdot)$ is given by
\begin{align*}
	\hat\sigma_j^2(x_j) &= \frac{1}{n}\sumin w_{ij}(x_j) \cdot
	\Big(Y_i - \big(\bar{Y}+\sum_{j=1}^{d}\hatm_j(X_{ij})\big)\Big)^2 \\
	&= \frac{1}{n}\sumin w_{ij}(x_j) \bigg(\sum_{j=1}^{d}\big(m_j(X_{ij})-\hatm_j(X_{ij})\big) + \epsilon_i + m_0 - \bar{Y}\bigg)^2.
\end{align*}
where $\hatm_j$ denotes the local linear estimators of $m_j$, and $w_{ij}$ denotes the corresponding local linear weight function.
Using $m_0-\bar{Y} = O_p(n^{-1/2})$ and
\begin{equation*}
	\sup_{x_j\in I}\abs{\hatm_j(x_j)-m_j(x_j)}
	=O_p\bigg(\sqrt{\frac{\log n}{nh_j}}\bigg) + o_p(h_j),
\end{equation*}
we obtain
\begin{equation}\label{eq:sigma_unif}
	\sup_{x_j\in I}\abs{\hat\sigma_j^2(x_j)-\sigma_j^2(x_j)}
	=
	O_p\bigg(\sqrt{\frac{\log n}{nh_j}}\bigg)+o_p(h_j).
\end{equation}
By the definitions of $\widehat\Delta_{nj}^{\#}(I)$ and $\Delta_{nj}^{\#}(I)$,
\begin{equation*}
	\begin{split}
		\abs{\widehat\Delta_{nj}^{\#}(I)-\Delta_{nj}^{\#}(I)}
		\le&\, \nu_0^{-1/2}\sup_{x_j \in I}\abs{nh_j\sigma_j(x_j)^{-2}p_j(x_j)}^{1/2}
		\cdot \sup_{x_j \in I} \abs{1-\left\{\frac{\sigma_j^2(x_j)\hatp_j(x_j)}{\hat{\sigma}_j^2(x_j)p_j(x_j)}\right\}^{1/2}}\\
		&\cdot \sup_{x_j \in I}\abs{\hatm_j^{\#}(x_j) - m_j(x_j)}.
	\end{split}
\end{equation*}
Under the stronger assumption (B4), together with the bandwidth condition (B3), the uniform error of the bias-corrected estimator satisfies
\begin{equation*}
	\sup_{x_j \in I}\abs{\hatm_j^{\#}(x_j) - m_j(x_j)} = O_p\bigg(\sqrt{\frac{\log n}{nh_j}} + h_j^2b_j\bigg)
\end{equation*}
Since $p_j$ and $\sigma_j^2$ are bounded away from zero on $I$, the results \eqref{eq:p_unif} and \eqref{eq:sigma_unif} imply
\begin{equation*}
	\sup_{x_j\in I}
	\left|
	\left\{
	\frac{\hatp_j(x_j)\sigma_j^2(x_j)}
	{p_j(x_j)\hat\sigma_j^2(x_j)}
	\right\}^{1/2}-1
	\right|
	=O_p\bigg(\sqrt{\frac{\log n}{nh_j}}\bigg)
	+o_p(h_j).
\end{equation*}
Since $A_{nj}(I)=O(\sqrt{\log n})$ and $\sqrt{nh_j}A_{nj}(I)h_j^2b_j=o(1)$, the preceding bounds yield
\begin{equation}
	\begin{split}
		&A_{nj}(I)\big\{\widehat\Delta_{nj}^{\#}(I)-\Delta_{nj}^{\#}(I) \big\}  \\
		&\quad = O_p\bigg(\sqrt{\frac{(\log n)^3}{nh_j}} + \sqrt{h_j^4b_j^2(\log n)^2}\bigg) + o_p\Big(h_j\log n + \sqrt{nh_j^7 b_j^2 \log n}\Big).
	\end{split}
\end{equation}
Consequently, for
\begin{equation*}
	\delta_n := \sqrt{\frac{(\log n)^3}{nh_j}} + \sqrt{h_j^4b_j^2(\log n)^2} + h_j\log n + \sqrt{nh_j^7 b_j^2 \log n} \to 0 \quad \text{as} \quad n\to\infty,
\end{equation*}
we have
\[
\Prob\left(
A_{nj}(I)
\abs{\widehat\Delta_{nj}^{\#}(I)-\Delta_{nj}^{\#}(I)}
>\delta_n
\right)\longrightarrow0.
\]
Combining Theorem~\ref{thm:unif_CLT}, the continuity of $u\mapsto\exp(-2e^{-u})$ and the fact that
\begin{align*}
	&\Prob\left(
	A_{nj}(I)\{\Delta_{nj}^{\#}(I)-B_{nj}(I)\}\le u-\delta_n
	\right)
	-\Prob\left(
	A_{nj}(I)
	\abs{\widehat\Delta_{nj}^{\#}(I)-\Delta_{nj}^{\#}(I)}
	>\delta_n
	\right)\\
	&\quad\le \Prob\left(T_{nj}^{\#}(I)\le u\right)\\
	&\quad\le
	\Prob\left(
	A_{nj}(I)\{\Delta_{nj}^{\#}(I)-B_{nj}(I)\}\le u+\delta_n
	\right)
	+\Prob\left(
	A_{nj}(I)
	\abs{\widehat\Delta_{nj}^{\#}(I)-\Delta_{nj}^{\#}(I)}
	>\delta_n
	\right)
\end{align*}
completes the proof of Corollary~\ref{cor:unif_CLT}.
\qed

\section{Proofs of additional results}\label{suppsec:proof_lemmas}

\subsection{Proof of Theorem~\ref{thm:mtilde_asymp}}\label{subsec:proof_mtilde}
Similar to the proof of Theorem~\ref{thm:mhat_asymp}, we only establish the result for $\tildem_j^\BR$.

For $\ell=1,\dots,L$ and $j=1,\dots,d$, let
\begin{align*}
	c_{j,\ell} := \frac{\sum_{\ell'=1}^L h_{\ell' j}^4	-	\big(\sum_{\ell'=1}^L h_{\ell' j}^2\big) b_{j,\ell}^2}{	L \cdot \sum_{\ell'=1}^L h_{\ell' j}^4 -	\big(\sum_{\ell'=1}^L h_{\ell' j}^2\big)^2}.
\end{align*}
By the definition, it holds that
\begin{equation*}
	\sumellL c_{j,\ell} = 1, \quad \sumellL c_{j,\ell} \cdot b_{j,\ell}^2 = 0.
\end{equation*}
Also, since $b_{j,\ell} \asymp b_j$ for $\ell = 1,\dots, L$ by (A6), we have $c_{j,\ell} = O(1)$.
Using these weights, we can write $\tildem_j^\BR$ as a weighted average with respect to $\hatm_j(\cdot;\bmh_\ell)$:
\begin{align*}
	\tildem_{j}^{\BR}(x_j) = \sumellL c_{j,\ell} \cdot  \hatm_j(x_j;\bmh_\ell), \quad x_j\in[0,1].
\end{align*}
Applying the uniform asymptotic expansion given in Proposition~\ref{prop:bias_expansion} to each $\hatm_j(x_j;\bmh_\ell)$, we obtain
\begin{equation}\label{eq:mtilde_expansion}
	\begin{split}
		\tildem_j^{\BR}(x_j)
		&= m_j(x_j) + \sumellL c_{j,\ell} \cdot \frac{1}{2}b_{j,\ell}^2 \fraku_0^\top \bmN_j(x_j;b_{j,\ell})^{-1} \bmgamma_j(x_j;b_{j,\ell}) m_j''(x_j) \\
		&\quad + \sumellL c_{j,\ell} \cdot \fraku_0^\top \tildebmm_j^{A, \tp}(x_j;b_{j,\ell}) \\
		&\quad + O_p\bigg(b_j^3 + b_{-j}^3 + \sqrt{\frac{\log n}{n}}\bigg),
	\end{split}
\end{equation}
where the $O_p$ term is uniform in $x_j\in[0,1]$.

Recall that for $x_j \in I_j=[2b_j,\, 1-2b_j]$, both $\bmN_j(x_j;b_{j,\ell})$ and $\bmgamma_j(x_j;b_{j,\ell})$ coincide with the corresponding complete moments.
Hence, the second term in the right-hand side of \eqref{eq:mtilde_expansion} vanishes, as $\sumellL c_{j,\ell} \cdot b_{j,\ell}^2=0$.
In contrast, for $x_j \in [0,1]\setminus I_j = [0,\,2b_j) \cup (1-2b_j,\,1]$, the leading bias term of order $b_j^2$ in the decomposition \eqref{eq:mtilde_expansion} remains non-zero.
Hence, it holds that
\begin{equation}\label{BR_h_lj_prop}
	\begin{split}
		\sup_{x_j \in I_j} \bigg|\sumellL c_{j,\ell} \cdot \frac{1}{2}b_{j,\ell}^2 \fraku_0^\top \bmN_j(x_j;b_{j,\ell})^{-1} \bmgamma_j(x_j;b_{j,\ell}) m_j''(x_j)\bigg| &= 0, \\
		\sup_{x_j \in [0,1]} \bigg|\sumellL c_{j,\ell} \cdot \frac{1}{2}b_{j,\ell}^2 \fraku_0^\top \bmN_j(x_j;b_{j,\ell})^{-1} \bmgamma_j(x_j;b_{j,\ell}) m_j''(x_j)\bigg| &= O(b_j^2).
	\end{split}	
\end{equation}
Furthermore, by adapting the standard kernel smoothing arguments and using the assumption that $b_{j,\ell}\asymp  b_j$ for $\ell=1,\dots,L$, we can show that
\begin{equation}\label{mtilde_A_lj_prop}
	\begin{split}
		|\fraku_0^\top \tildebmm_j^{A, \tp}(x_j;b_{j,\ell})| &=  O_p(n^{-1/2}b_j^{-1/2} ), \quad x_j \in (0,1), \\
		\int_0^1 |\fraku_0^\top \tildebmm_j^{A, \tp}(x_j;b_{j,\ell})|^2 \d x_j &= O_p(n^{-1}b_j^{-1}), \\
		\sup_{x_j\in[0,1]} |\fraku_0^\top \tildebmm_j^{A, \tp}(x_j;b_{j,\ell})| &= O_p \Bigg(\sqrt{\frac{\log n}{nb_j}}  \Bigg).
	\end{split}
\end{equation}
Applying \eqref{BR_h_lj_prop} and \eqref{mtilde_A_lj_prop} to \eqref{eq:mtilde_expansion}, we obtain the desired convergence rates for $\tildem_j^{\BR}$.
\qed

\subsection{Proof of Lemma~\ref{lem:delta_j_k}}\label{subsec:proof_delta_j_k}

For simplicity, we present only the proof for the second part of Lemma~\ref{lem:delta_j_k}.	
Define the empirical weight
\begin{equation*}
	\hat{\boldsymbol\kappa}_{ij}(x_j) := \hatbmM_{jj}(x_j)^{-1} \frac{1}{n} \mathbf{vec}\bigg(\frac{X_{ij}-x_j}{h_j}\bigg) K_{h_j}(x_j,X_{ij}).
\end{equation*}
Note that
\begin{equation}\label{eq:kappa_sum}
	\sumin \hat{\bmkappa}_{ij}(x_j) = \Bigg[\text{first column of } \sumin \hat{\bmkappa}_{ij}(x_j)	\cdot	\bmvec\bigg( \frac{X_{ij}-x_j}{h_j}\bigg)^{\top} = I_2 \Bigg]	= \fraku_0,
\end{equation}
where $I_2$ denotes the $2\times2$ identity matrix.

For $1\le j\neq k\le d$, let
\begin{align*}
	\tildebmm_j^{B, \tp}(x_j) :=& \, \sumin  \hat{\bmkappa}_{ij}(x_j) \cdot \Bigg(m_j(X_{ij}) - \mbf{vec}\left( \frac{X_{ij}-x_j}{h_j}\right)^\top \bmm_j^{\tp}(x_j)  \Bigg)\\
	\tildebmm_{jk}^{C, \tp}(x_j) :=& \, \sumin \int_0^1  \hatbmM_{jj}(x_j)\hat{\bmkappa}_{ij}(x_j)K_{h_k}(x_k, X_{ik}) \\
	&\hspace{3cm} \cdot \Bigg(m_k(X_{ik}) - \mbf{vec}\bigg( \frac{X_{ik}-x_k}{h_k}\bigg)^\top \bmm_k^{\tp}(x_k)  \Bigg)\d x_k
\end{align*}
Then, the SBF equations at \eqref{eq:SBF_integral_eq} can be written as
\begin{align*}
	\hatbmm_j^{\tp}(x_j) &= \bmm_j^{\tp}(x_j) + (m_0 - \bar{Y})\cdot \fraku_0 + \tildebmm_j^{A, \tp}(x_j) + \tildebmm_j^{B, \tp}(x_j) + \sumkneqj \hatbmM_{jj}(x_j)^{-1}\tildebmm_{jk}^{C, \tp}(x_j) \\
	&\quad - \sumkneqj \int_0^1 \hatbmM_{jj}(x_j)^{-1}\hatbmM_{jk}(x_j, x_k)\left(\hatbmm_k^{\tp}(x_k) - \bmm_k^{\tp}(x_k)\right)\d x_k.
\end{align*}
Here, we used the normalization property that $\int_0^1 K_{h_k}(x_k,u)\d x_k = 1$ for all $u\in[0,1]$.
This also implies
\begin{equation}\label{eq:SBF_eq_decomp}
	\begin{split}
		&\hspace{-0.5cm} \hat{\Delta}_j^{\tp}(x_j) + \sumkneqj\int_0^1  \hatbmM_{jj}(x_j)^{-1}\hatbmM_{jk}(x_j, x_k)\hat{\Delta}_k^{\tp}(x_k)\d x_k \\
		&= (m_0 - \bar{Y}) \cdot \fraku_0 + \hatI_{1,j}(x_j) + \sumkneqj \Big(\hatI_{2,jk}(x_j) + \hatI_{3,jk}(x_j)\Big),
	\end{split}	
\end{equation}
where
\begin{align*}
	\hatI_{1,j}(x_j) :=& \, \tildebmm_j^{B, \tp}(x_j) - \frac{h_j^2}{2}\bmN_j(x_j)^{-1}\bmgamma_j(x_j)m_j''(x_j), \\
	\hatI_{2,jk}(x_j) :=& \, \hatbmM_{jj}(x_j)^{-1}\tildebmm_{jk}^{C, \tp}(x_j) \\
	&\hspace{1cm} - \frac{h_k^2}{2}\int_0^1 \hatbmM_{jj}(x_j)^{-1}\hatbmM_{jk}(x_j, x_k) \bmN_k(x_k)^{-1}\bmgamma_k(x_k)m_k''(x_k) \d x_k, \\	
	\hatI_{3,jk}(x_j) :=& \, \int_0^1 \hatbmM_{jj}(x_j)^{-1}\hatbmM_{jk}(x_j, x_k)\tildebmm_k^{A, \tp}( x_k) \d x_k.
\end{align*}
Since $m_0 - \bar{Y} = O_p(n^{-1/2})$ by (A3) and (A4), it suffices to show that
\begin{align}
	\sup_{x_j\in[0,1]} |\hatI_{1,j}(x_j)| &= O_p\bigg(h_j^3 + \sqrt{\frac{\log n}{nh_j}}h_j^2 \bigg), \label{tildem_claim1} \\
	\sup_{x_j\in[0,1]} |\hatI_{2,jk}(x_j)| &= O_p\bigg( h_k^3 + h_k^2h_j + \sqrt{\frac{\log n \cdot h_k^2}{nh_j}} \bigg), \label{tildem_claim2} \\
	\sup_{x_j\in[0,1]} |\hatI_{3,jk}(x_j)| &= O_p\bigg(\sqrt{\frac{\log n}{n}} \bigg). \label{tildem_claim3}
\end{align}

\underline{\bf Proof of \eqref{tildem_claim1}}. \quad
For $1\le j\le d$ and $x_j\in[0,1]$, let
\begin{equation*}
	\bmmu_j(x_j) := \int_0^1 \bmvec\bigg(\frac{u - x_j}{h_j}\bigg) \cdot K_{h_j}(x_j, u)\d u \in \bbR^2.
\end{equation*}
By standard kernel theory \citep{jeon2020additive}, under (A1), (A2) and (A5), we can show that
\begin{align}
	\sup_{x_j\in[0,1]} \Big| \hatbmM_{jj}(x_j) - \bmN_{j}(x_j)p_j(x_j) \Big| &= O_p\bigg(h_j + \sqrt{\frac{\log n}{n{h_j}}}\bigg)	\label{M_jj_approx}\\
	\sup_{(x_j, x_k) \in [0,1]^2} \Big| \hatbmM_{jk}(x_j, x_k) - \bmmu_j(x_j)\bmmu_k(x_k)^\top p_{jk}(x_j, x_k) \Big| &= O_p\bigg(h_j + h_k + \sqrt{\frac{\log n}{n{h_j}{h_k}}} \bigg)	\label{M_jk_approx}.
\end{align}
By (A1) and (A4), together with \eqref{M_jj_approx}, we obtain
\begin{align*}
	\tildebmm_j^{B, \tp}(x_j) &= \hatbmM_{jj}(x_j)^{-1} \cdot \meanin \bmvec\bigg( \frac{X_{ij}-x_j}{h_j}\bigg) K_{h_j}(x_j, X_{ij}) \cdot \bigg(\frac{(X_{ij}-x_j)^2}{2}m_j''(x_j) + o_p(h_j^2) \bigg) \\
	&= \Bigg(p_j(x_j)^{-1}\bmN_j(x_j)^{-1}+ O_p\bigg(h_j + \sqrt{\frac{\log n}{n{h_j}}}\bigg)\Bigg)  \\
	&\quad \times \Bigg( \int_0^1\left( \frac{v-x_j}{h_j}\right)^2\! \bmvec\!\bigg( \frac{v-x_j}{h_j}\bigg) K_{h_j}(x_j, v)p_j(x_j)\d v + O_p\bigg(h_j + \sqrt{\frac{\log n}{n{h_j}}}\bigg) \Bigg) \\
	&\quad \times  \bigg(\frac{h_j^2}{2}m_j''(x_j) + o_p(h_j^2) \bigg)\\
	&= \frac{h_j^2}{2}\bmN_j(x_j)^{-1}\bmgamma_j(x_j)m_j''(x_j) + O_p\bigg(h_j^3 + \sqrt{\frac{\log n}{nh_j}}h_j^2 \bigg) ,
\end{align*}
where the $O_p$ terms hold uniformly over $x_j \in [0,1]$.

\underline{\bf Proof of \eqref{tildem_claim2}}. \quad
By (A4), we have
\begin{equation}\label{eq:C_tp_app1}
	\begin{split}
		\hatbmM_{jj}(x_j)^{-1}\tildebmm_{jk}^{C, \tp}(x_j)	&=	\sumin \hat{\bmkappa}_{ij}(x_j) \cdot	\int_0^1 \frac{(X_{ik}-x_k)^2}{2}K_{h_k}(x_k, X_{ik})m_k''(x_k)\d x_k \\
		&\quad + O_p(h_k^3),
	\end{split}	
\end{equation}
where the $O_p$ term holds uniformly over $x_j \in [0,1]$.
Let
\begin{align*}
	\hat{\bmeta}_{ij}(x_j) := \hatbmM_{jj}(x_j)^{-1} \bmvec\bigg( \frac{X_{ij}-x_j}{h_j}\bigg)	\cdot	\int_0^1 \frac{(X_{ik}-x_k)^2}{2}K_{h_k}(x_k, X_{ik})m_k''(x_k)\d x_k.
\end{align*}
Since $|\hat\bmeta_{ij}(x_j) \bbI(|X_{ij}-x_j| \le h_j)| \le C h_k^2$ for some sufficiently large constant $C>0$, applying Proposition S.1 in \cite{jeon2020additive} yields
\begin{equation}\label{eq:C_tp_app2}
	\sup_{x_j\in [0,1]} \Big| \meanin K_{h_j}(x_j, X_{ij}) \big(\hat{\bmeta}_{ij}(x_j)- \E[\hat{\bmeta}_{ij}(x_j) \mid X_{ij}] \big) \Big| = O_p\bigg(\sqrt{\frac{\log n \cdot h_k^2}{nh_j}} \bigg).
\end{equation}
Note that
\begin{align*}
	\E[\hat{\bmeta}_{ij}(x_j) \mid X_{ij}] &= \hatbmM_{jj}(x_j)^{-1} \bmvec\bigg( \frac{X_{ij}-x_j}{h_j}\bigg) \\
	&\hspace{1cm} \cdot \int_0^1 \int_0^1 \frac{(v-x_k)^2}{2}K_{h_k}(x_k, v) \frac{p_{jk}(X_{ij}, v)}{p_j(X_{ij})} m_k''(x_k) \d v  \d x_k.
\end{align*}
Combining \eqref{eq:C_tp_app1} and \eqref{eq:C_tp_app2}, we have
\begin{equation}\label{eq:C_tp_app3}
	\begin{split}
		&\hspace{-0.5cm} \hatbmM_{jj}(x_j)^{-1}\tildebmm_{jk}^{C, \tp}(x_j)  \\
		&= \sumin \hat{\bmkappa}_{ij}(x_j)\int_0^1 \int_0^1 \frac{(v-x_k)^2}{2}K_{h_k}(x_k, v) \frac{p_{jk}(X_{ij}, v)}{p_j(X_{ij})} m_k''(x_k) \d v \d x_k \\
		&\quad + O_p\bigg(h_k^3 + \sqrt{\frac{\log n \cdot h_k^2}{nh_j}} \bigg)\\
		&= \sumin \hat{\bmkappa}_{ij}(x_j) \int_0^1 \frac{p_{jk}(X_{ij}, x_k)}{p_j(X_{ij})}
		\bigg[ \int_0^1 \frac{(v-x_k)^2}{2}K_{h_k}(x_k, v) \d v \bigg] m_k''(x_k)\d x_k\\
		&\quad + O_p\bigg(h_k^3 + \sqrt{\frac{\log n \cdot h_k^2}{nh_j}} \bigg)
	\end{split}	
\end{equation}
uniformly over $x_j\in[0,1]$.
Let $\mu_{\ell,j}(x) = \int_{-x/h_j}^{(1-x)/h_j} v^\ell K(v)\d v$ and $\mu_\ell = \int_{-1}^1v^\ell K(v)\d v$.
Since $K_h(x_, u) \le 2K_h(x-u)$, there exists a constant $C'>0$ such that
\begin{align*}
	&\Bigg|	\int_0^1 \frac{p_{jk}(X_{ij}, x_k)}{p_j(X_{ij})} \bigg[	\int_0^1 \frac{(v-x_k)^2}{2}K_{h_k}(x_k, v)\, \d v \bigg] m_k''(x_k)\, \d x_k	- h_k^2 \mu_2 \int_0^1 \frac{p_{jk}(X_{ij}, x_k)}{p_j(X_{ij})}	m_k''(x_k)\, \d x_k	\Bigg| \\
	&\quad \le	\int_{[0,2h_k]\cup[1-2h_k,1]} \frac{p_{jk}(X_{ij}, x_k)}{p_j(X_{ij})} \, 2h_k^2 \mu_{2,k}(x_k)\, m_k''(x_k)\, \d x_k \\
	&\quad \le C'h_k^3.
\end{align*}
for all $x_j\in[0,1]$.
Combining the above result with \eqref{eq:kappa_sum} and \eqref{eq:C_tp_app3}, we have
\begin{equation}\label{tildeC_approx_1}
	\begin{split}
		&\hspace{-0.5cm} \hatbmM_{jj}(x_j)^{-1}\tildebmm_{jk}^{C, \tp}(x_j)  \\
		&= \frac{\mu_2}{2}h_k^2 \cdot \sumin \hat{\bmkappa}_{ij}(x_j)\int_0^1 \frac{p_{jk}(X_{ij}, x_k)}{p_j(X_{ij})}m_k''(x_k)\d x_k + O_p\bigg(h_k^3 + \sqrt{\frac{\log n \cdot h_k^2}{nh_j}} \bigg) \\
		&= \frac{\mu_2}{2}h_k^2 \cdot \sumin \hat{\bmkappa}_{ij}(x_j)\int_0^1 \frac{p_{jk}(x_j, x_k)}{p_j(x_j)}m_k''(x_k)\d x_k	+ O_p\bigg(h_k^3 +h_jh_k^2+ \sqrt{\frac{\log n \cdot h_k^2}{nh_j}} \bigg) \\
		&= \frac{\mu_2}{2}h_k^2 \cdot \fraku_0 \cdot \int_0^1 \frac{p_{jk}(x_j, x_k)}{p_j(x_j)}m_k''(x_k)\d x_k + O_p\bigg(h_k^3 +h_jh_k^2+ \sqrt{\frac{\log n \cdot h_k^2}{nh_j}} \bigg).
	\end{split}	
\end{equation}

Moreover, by \eqref{M_jj_approx} and \eqref{M_jk_approx}, we get
\begin{equation}\label{tildeC_approx_2}
	\begin{split}
		&\hspace{-0.5cm} \int_0^1 \hatbmM_{jj}(x_j)^{-1}\hatbmM_{jk}(x_j, x_k) \bmN_k(x_k)^{-1}\bmgamma_k(x_k)m_k''(x_k) \d x_k \\
		&= \int_0^1 \bmN_j(x_j)^{-1}\bmmu_j(x_j) \bmmu_k(x_k)^{\top} \bmN_k(x_k)^{-1}\bmgamma_k(x_k) \frac{p_{jk}(x_j, x_k)}{p_j(x_j)}m_k''(x_k)\d x_k \\
		&\hspace{2cm} + O_p\bigg( h_j + h_k + \sqrt{\frac{\log n}{nh_jh_k}}\bigg) \\
		&= \int_0^1\fraku_0 \fraku_0^{\top} \bmgamma_k(x_k) \frac{p_{jk}(x_j, x_k)}{p_j(x_j)}m_k''(x_k)\d x_k + O_p\bigg( h_j + h_k + \sqrt{\frac{\log n}{nh_jh_k}}\bigg) \\
		&=\int_0^1\fraku_0 \cdot \mu_{2,k}(x_k) \frac{p_{jk}(x_j, x_k)}{p_j(x_j)}m_k''(x_k)\d x_k + O_p\bigg( h_j + h_k + \sqrt{\frac{\log n}{nh_jh_k}}\bigg) \\
		&= \mu_2 \cdot \fraku_0\int_0^1 \frac{p_{jk}(x_j, x_k)}{p_j(x_j)}m_k''(x_k)\d x_k + O_p\bigg( h_j + h_k + \sqrt{\frac{\log n}{nh_jh_k}}\bigg),
	\end{split}	
\end{equation}
where the second equality follows from the fact that $\bmN_j(x_j)^{-1}\bmmu_j(x_j) = \fraku_0$.
Also, for the third and fourth equalities, we have used $\mu_{2,k}(x_k)=\mu_2$ for $x_k\in[2h_k, 1-2h_k]$.
Combining \eqref{tildeC_approx_1} and \eqref{tildeC_approx_2} completes the proof of \eqref{tildem_claim2}.

\underline{\bf Proof of \eqref{tildem_claim3}}. \quad
For $1 \le i \le n$ and $1 \le j \neq k \le d$, we define two-dimensional weight functions $\bmw_{jk}^{(i)}(\cdot) = \big(w_{jk,0}^{(i)}(\cdot),w_{jk,1}^{(i)}(\cdot)\big)^\top : [0,1] \to \bbR^2$ by
\begin{equation*}
	\bmw_{jk}^{(i)}(x_j) := \int_0^1\hatbmM_{jj}(x_j)^{-1}\hatbmM_{jk}(x_j, x_k)\hatbmM_{kk}(x_k)^{-1} \bmvec\bigg(\frac{X_{ik}-x_k}{h_k} \bigg) \cdot K_{h_k}(x_k, X_{ik})\d x_k.
\end{equation*}
Then, we can write
\begin{align*}
	\int_0^1 \hatbmM_{jj}(x_j)^{-1}\hatbmM_{jk}(x_j, x_k)\tildebmm_k^{A, \tp}(x_k) \d x_k  =    \meanin\bmw_{jk}^{(i)}(x_j) \cdot \varepsilon_i.
\end{align*}
Consequently, it suffices to show that
\begin{equation}\label{wsum_sup}
	\sup_{x_j\in[0,1]}\bigg| \meanin w_{jk, \ell}^{(i)}(x_j) \cdot \varepsilon_i \bigg| = O_p\bigg(\sqrt{\frac{\log n}{n}} \bigg), \quad \ell=0,1.
\end{equation}

Since the proofs of \eqref{wsum_sup} for $\ell=0$ and $\ell=1$ are analogous, we focus only on the case where $\ell=0$.
We employ a truncation argument to prove \eqref{wsum_sup}.
Let $\tau_n := \sqrt{n / \log n}$.
For $1 \le j \neq k \le d$ and $x_j\in[0,1]$, define $Z_{jk, 0}(x_j) := \meanin \xi_{jk, 0}(x_j)$, where
\begin{equation*}
	\xi_{jk, 0}^{(i)}(x_j) := w_{jk, 0}^{(i)}(x_j) \cdot \Big[\varepsilon_i\cdot\bbI(|\varepsilon_i|\le \tau_n) - \E\big(\varepsilon_i\cdot\bbI(|\varepsilon_i|\le\tau_n) \, \big| \, \bmX_i\big) \Big].
\end{equation*}
We decompose
\begin{equation}\label{wsum_sup_decomp}
	\begin{split}
		\meanin w_{jk, \ell}^{(i)}(x_j) \cdot \varepsilon_i &= Z_{jk, 0}(x_j) + \meanin w_{jk, 0}^{(i)}(x_j) \cdot \varepsilon_i\cdot\bbI(|\varepsilon_i|> \tau_n) \\
		&\quad + \meanin w_{jk, 0}^{(i)}(x_j) \cdot \E\big(\varepsilon_i\cdot\bbI(|\varepsilon_i|> \tau_n) \, \big| \, \bmX_i \big).	
	\end{split}	
\end{equation}

First, there exists a constant $0<C<\infty$ such that
\begin{align*}
	\Prob\big(\max_{1 \le i \le n}|\varepsilon_i | \le \tau_n\big) \ge 1 - n\cdot \Prob\big(|\varepsilon_1|^\alpha > \tau_n^\alpha\big) \ge 1- C \frac{(\log n)^{\alpha/2}}{n^{\alpha/2 - 1}},
\end{align*}
where $\alpha$ is the constant that appears in (A5).
This implies
\begin{equation}\label{wsum_sup_decomp_prop1}
	\lim_{n \to \infty}\Prob\bigg(\sup_{x_j\in[0,1]}\Big|\meanin w_{jk, 0}^{(i)}(x_j) \varepsilon_i\cdot\bbI(|\varepsilon_i|>\tau_n)\Big|=0 \bigg) = 1.
\end{equation}
Moreover, using \eqref{M_jj_approx} and \eqref{M_jk_approx}, along with (A1), (A2) and (A5), we can show that there exists a constant $C'>0$ such that
\begin{align}\label{w_jk_bdd}
	\lim_{n \to \infty}\Prob\bigg( \sup_{x_j\in[0,1]} |w_{jk, 0}^{(i)}(x_j)| \le C' \bigg) = 1.
\end{align}
Then, by (A3), it holds that
\begin{align}\label{wsum_sup_decomp_prop2}
	\sup_{x_j\in[0,1]}\Big|\meanin w_{jk, 0}^{(i)}(x_j) \E\Big(\varepsilon_i\cdot\bbI(|\varepsilon_i|>\tau_n) \, \Big| \, \bmX_i\Big)\Big| = O_p(\tau_n^{-1}).
\end{align}
Substituting \eqref{wsum_sup_decomp_prop1} and \eqref{wsum_sup_decomp_prop2} into \eqref{wsum_sup_decomp}, we obtain
\begin{equation*}
	\sup_{x_j\in[0,1]}\bigg| \meanin w_{jk, 0}^{(i)}(x_j) \cdot \varepsilon_i \bigg| = \sup_{x_j\in[0,1]} \big| Z_{jk, 0}(x_j) \big| + O_p(\tau_n^{-1}).
\end{equation*}

For some fixed $\gamma>2$, let $\cN(n^{-\gamma})$ denote an $n^{-\gamma}$-net of $[0,1]$.
We further decompose
\begin{equation}\label{wsum_sup_decomp2}
	\sup_{x_j\in[0,1]}\big| Z_{jk, 0}(x_j) \big| \le \max_{x_j\in\cN(n^{-\gamma})}\big| Z_{jk, 0}(x_j) \big| + \sup_{|x_j - x_j'| \le n^{-\gamma}} \big| Z_{jk, 0}(x_j) - Z_{jk, 0}(x_j') \big|.
\end{equation}
We first bound the first term in the above decomposition \eqref{wsum_sup_decomp2}.
Since $\xi_{jk, 0}^{(i)}(x_j)$ for $i=1,\dots,n$ are i.i.d. random variables with $\E[\xi_{jk, 0}^{(i)}(x_j)] = 0$ for all $x_j\in[0,1]$, there exist some constants $0<C',C''<\infty$ such that
\begin{align*}
	\sup_{x_j\in[0,1]}\Prob\Big( Z_{jk, 0}(x_j) >C \cdot \tau_n^{-1}\Big) &= \sup_{x_j\in[0,1]}\Prob\bigg(\sqrt{\frac{\log n}{n}}\sumin \xi_{jk, 0}^{(i)}(x_j) > C\log n \bigg)\\
	& \le n^{-C} \sup_{x_j\in[0,1]} \bigg(\E \exp\bigg(\sqrt{\frac{\log n}{n}}\xi_{jk, 0}^{(1)}(x_j)\bigg) \bigg)^n\\
	& \le n^{-C} \cdot \bigg(1+e^{2C'}\cdot \frac{\log n}{n}\cdot \sup_{x_j\in[0,1]}\E(\xi_{jk, 0}^{(1)}(x_j)^2) \bigg)^n\\
	& \le n^{-C + C''},
\end{align*}
where for the last inequality we have used $\sup_{x_j\in[0,1]}\E\big[\xi_{jk, 0}^{(1)}(x_j)^2\big] \lesssim 1$, which follows from \eqref{w_jk_bdd} and (A3).
Similarly, we can show that
\begin{equation*}
	\sup_{x_j\in[0,1]}\Prob\Big( |Z_{jk, 0}(x_j)| >C \cdot \tau_n^{-1}\Big) \le 2 n^{-C+C''}.
\end{equation*}
Together with the fact that $|\cN(n^{-\gamma})| \lesssim n^c$, this implies
\begin{equation}\label{wsum_sup_decomp2_prop1}
	\sup_{x_j \in \cN(n^{-\gamma})} \lvert Z_{jk, 0}(x_j) \rvert = O_p(\tau_n^{-1}).	
\end{equation}

Next, we consider the supremum of the increment $| Z_{jk, 0}(x_j) - Z_{jk, 0}(x_j') |$ over all pairs $(x_j, x_j')$ with $\lvert x_j-x_j'\rvert \le n^{-\gamma}$.
Using the Lipschitz continuity of $K$, the fact that $n^{c_2} \min_{1 \le j \le d}h_j$ is bounded away from zero by (A3), and the inequality $K_h(x, u) \le 2K_h(x-u)$, we can show that
\begin{align*}
	\sup_{|x_j-x_j'| \le n^{-\gamma}} \big| \hatbmM_{jj}(x_j)  -  \hatbmM_{jj}(x_j') \big| &= O_p\big(n^{-\gamma} \cdot h_j^{-1}\big)\\
	\sup_{|x_j-x_j'| \le n^{-\gamma}} \big| \hatbmM_{jk}(x_j, x_k) - \hatbmM_{jk}(x_j', x_k) \big| &= O_p\big(n^{-\gamma} \cdot h_j^{-1}h_k^{-1}\big).
\end{align*}
This implies
\begin{equation}\label{w_jk_lipschitz}
	\sup_{|x_j-x_j'| \le n^{-\gamma}} \big| w_{jk, 0}^{(i)}(x_j) - w_{jk, 0}^{(i)}(x_j') \big| =    O_p\big( n^{-\gamma} \cdot h_j^{-1} h_k^{-1} \big).
\end{equation}
Then, we obtain
\begin{align}\label{wsum_sup_decomp2_prop2}
	\sup_{|x_j-x_j'| \le n^{-\gamma}} \big| Z_{jk, 0}(x_j) - Z_{jk, 0}(x_j') \big| = O_p\big(n^{-\gamma} \cdot h_j^{-1}h_k^{-1}\big) = O_p\big(n^{-\gamma + c_j+c_k}\big) =  o_p(\tau_n^{-1}),
\end{align}
where the last equality follows from the fact that $-\gamma + c_j + c_k \le -2 + c_j + c_k < -1$.
Applying \eqref{wsum_sup_decomp2_prop1} and \eqref{wsum_sup_decomp2_prop2} to \eqref{wsum_sup_decomp2}, we complete the proof of \eqref{tildem_claim3}.
\qed

\subsection{Proof of Lemma~\ref{lem:normal_rv_approx}}\label{subsec:proof_normal_rv_approx}

We present only the proof for $\hatm_j^\BR$.
Note that the same Gaussian variables $Z_{1j},\ldots,Z_{nj}$ suffice for both $\hatm_j^\BR$ and $\hatm_j^\WBR$, because the Gaussian approximation targets their common leading error term $\tildem_j^A := \fraku_0^\top \cdot \tildebmm_j^{A,\tp}$.

We first approximate the estimation error by a kernel-weighted error sum as follows.
\begin{equation}\label{eq:BR_weighted_epsilon_expansion}
	\hatm_j^{\BR}(x_j)-m_j(x_j)	= \frac{1}{n h_j p_j(x_j)}	\sumin 	K\bigg(\frac{x_j-X_{ij}}{h_j}\bigg) \varepsilon_i	+ R_{nj}^{\BR}(x_j), \quad x_j\in I.
\end{equation}
On the fixed interior interval $I$, \eqref{M_jj_approx} and the standard inequality for the kernel-weighted error sum give
\begin{equation}\label{eq:tilde_kernel_sum}
	\sup_{x_j\in I}\Bigg|
	\fraku_0^\top\tildebmm_j^{A,\tp}(x_j;h_j)
	-\frac{1}{nh_jp_j(x_j)}\sumin K\bigg(\frac{x_j-X_{ij}}{h_j}\bigg)\varepsilon_i
	\Bigg|
	=O_p\bigg(\sqrt{\frac{h_j\log n}{n}}+\frac{\log n}{nh_j}\bigg).
\end{equation}
Combining \eqref{eq:tilde_kernel_sum} with the uniform expansion in the proof of Theorem~\ref{thm:mhat_asymp}, we obtain
\begin{equation*}
	\begin{split}
		\sup_{x_j\in I}\big|R_{nj}^{\BR}(x_j)\big|
		=O_p\Bigg(&h_j^3+h_{-j}^3+\sqrt{\frac{\log n}{n}}
		+\frac{h_j^2}{b_j^2}\bigg\{b_j^3+b_{-j}^3+\sqrt{\frac{\log n}{n}}\bigg\}\\
		&+\frac{h_j^2}{b_j^2}\sqrt{\frac{\log n}{nb_j}}
		+\sqrt{\frac{h_j\log n}{n}}+\frac{\log n}{nh_j}\Bigg).
	\end{split}
\end{equation*}
It then follows from (B3) that
\begin{equation}\label{eq:BR_remainder_negligible}
	\sup_{x_j\in I}\big|R_{nj}^{\BR}(x_j)\big|
	=o_p\bigg(\frac{1}{\sqrt{nh_j}A_{nj}(I)}\bigg).
\end{equation}

We now approximate the kernel-weighted sum of the errors at \eqref{eq:BR_weighted_epsilon_expansion} by the corresponding kernel-weighted sum of normal random variables.
Since $1/2<(\alpha'-2)/\alpha'$, we can choose two fixed constants $q,\tau>0$ such that
\begin{equation}\label{eq:tau_choice_hash}
	\frac{2}{1-b} < q < \alpha' \quad \text{and} \quad \frac{1}{q} < \tau < \frac{1-b}{2},
\end{equation}
where $\alpha'$ and $b$ are the constants that appear in (B2) and (B3), respectively.
Let $X_{(1)j}\le\cdots\le X_{(n)j}$ denote the order statistics of $X_{1j},\ldots,X_{nj}$, and let $\varepsilon_{(i)}$ denote the random error corresponding to $X_{(i)j}$.
Notably, conditioned on $X_{1j},\ldots,X_{nj}$, the random variables $\varepsilon_{(1)},\ldots,\varepsilon_{(n)}$ are independent, have conditional mean zero, and satisfy $\Var(\varepsilon_{(i)} \mid X_{1j},\dots,X_{nj}) = \sigma_j^2(X_{(i)j})$.
Following the proof of Theorem 4.1 in \cite{hall2013simple} (see Equation (B.45) in that paper), we can show that there exist Gaussian random variables $Z_{1j},\ldots,Z_{nj}$ such that
\begin{itemize}
	\item[(i)] $Z_{1j},\dots,Z_{nj}$ are conditionally independent given $X_{1j},\dots,X_{nj}$;
	\item[(ii)] $Z_{ij}\mid X_{ij} \sim N\big(0,\sigma_j^2(X_{ij})\big)$ for $i=1,\dots,n$;
	\item[(iii)] it holds that
	\begin{equation}\label{eq:sakhanenko_hash}
		M_{nj} := \max_{1\le i\le n} \bigg| \sum_{k=1}^i \big(\varepsilon_{k}-Z_{kj}\big) \bigg| = O_p(n^\tau).
	\end{equation}	
\end{itemize}
Indeed, to establish the property (iii), Theorem 1 in \cite{gotze2008bounds} ensures that there exists a constant $C_q>0$ depending only on the constant $q>0$ such that
\begin{equation*}
	\Prob\big(M_{nj} > C_q n^{\tau} \big) \le C_q^{q}n^{-q\tau} \cdot \E\bigg(\sumin \E\big(|\varepsilon_i|^q \, \big| \, X_{ij}\big) \bigg) = O(n^{1-q\tau}) = o(1),
\end{equation*}
which implies \eqref{eq:sakhanenko_hash}.

Let $Z_{(i)j}$ denote the Gaussian variable associated with $X_{(i)j}$.
For $1\le i\le n$ and $x_j\in I$, define
\begin{equation*}
	a_{(i)j}(x_j) := K\bigg(\frac{x_j - X_{(i)j}}{h_j}\bigg) \quad \text{and} \quad S_{(i)j} := \sum_{k=1}^i \big(\varepsilon_{(k)} - Z_{(k)j}\big).
\end{equation*}
By the definition of the order statistics $X_{(i)j}$, the sequence $\big\{ (x_j - X_{(i)j}) /h_j : 1\le i\le n\big\}$ is monotonically increasing.
Since the total variation of $K$ is finite, there exists a constant $C_K>0$ depending only on the baseline kernel $K$ such that
\begin{equation}\label{eq:a_nj_prop}
	|a_{(n)j}(x_j)| + \sum_{i=1}^{n-1} |a_{(i+1)j}(x_j)-a_{(i)j}(x_j)| \le C_K, \quad x_j \in I.	
\end{equation}
Moreover, by Euler's summation formula,
\begin{equation}\label{eq:abel_hash}
	\sumin a_{(i)j}(x_j)\{\varepsilon_{(i)}-Z_{(i)j}\} = a_{(n)j}(x_j)S_{(n)j} - \sum_{i=1}^{n-1} S_{(i)j}\big\{a_{(i+1)j}(x_j)-a_{(i)j}(x_j)\big\}.	
\end{equation}
Combining \eqref{eq:sakhanenko_hash}--\eqref{eq:abel_hash}, we have
\begin{equation*}
	\sup_{x_j\in I} \Bigg| \sumin K\bigg(\frac{x_j-X_{ij}}{h_j}\bigg) \{\varepsilon_i-Z_{ij}\} \Bigg| = \sup_{x_j\in I} \Bigg| \sumin a_{(i)j}(x_j) \{\varepsilon_{(i)}-Z_{(i)j}\} \Bigg| = O_p(n^\tau).
\end{equation*}
Since $p_j$ is bounded away from zero on $[0,1]$, this implies
\begin{equation}\label{eq:epsilon_Z_negligible_hash}
	\sup_{x_j\in I} \Bigg| \frac{1}{n h_j p_j(x_j)} \sumin K\bigg(\frac{x_j-X_{ij}}{h_j}\bigg) \{\varepsilon_i-Z_{ij}\} \Bigg| = O_p\bigg(\frac{n^\tau}{n h_j}\bigg) = o_p\bigg(\frac{1}{\sqrt{n h_j}A_{nj}(I)}\bigg),
\end{equation}
where the last equality follows from (B3), \eqref{eq:tau_choice_hash} and the fact that $A_{nj}(I) = O(\sqrt{\log n})$.
Applying \eqref{eq:epsilon_Z_negligible_hash} to \eqref{eq:BR_weighted_epsilon_expansion}, we obtain
\begin{equation*}
	\sup_{x_j\in I} \Bigg| \hatm_j^\BR(x_j)-m_j(x_j) - \frac{1}{n h_j p_j(x_j)} \sumin	K\bigg(\frac{x_j-X_{ij}}{h_j}\bigg)Z_i \Bigg| = o_p\bigg(\frac{1}{\sqrt{n h_j}A_{nj}(I)} \bigg).	
\end{equation*}
This completes the proof of Lemma~\ref{lem:normal_rv_approx}.
\qed

\subsection{Proof of Lemma~\ref{lem:gaussian_process_approx}}\label{subsec:proof_gaussian_process_approx}

We adopt the notation of \cite{chernozhukov2014gaussian}, which is summarized below.
Let $P$ denote the distribution of $X_j$ and $P_n$ be its empirical version; that is,
\begin{equation*}
	P_n f := \frac{1}{n}\sumin f(X_{ij}), \qquad P f := \E f(X_j).
\end{equation*}
We also denote $\bbG_n :=\sqrt{n}(P_n-P)$.
For a probability measure $Q$ and a measurable function $g$, let $\norm{g}_{Q,2}:=\big( \int g^2 \, \d Q \big)^{1/2}$.
Given a function class $\cF$ with an envelope $F$, we write $\norm{\bbG_n}_{\cF}:= \sup_{f\in \cF} \abs{\bbG_n f}$.
Moreover, we define $ M(F) := \max_{1\le i\le n} F(X_{ij})$ and $\|M(F)\|_2 := \big(\E \{M(F)\}^2\big)^{1/2}$.

Let $\cX_j:=\sigma(X_{1j},\ldots,X_{nj})$ denote the $\sigma$-field generated by $X_{1j},\dots,X_{nj}$.
Note that, conditioned on $\cX_j$, the process $G_{nj}$ defined at \eqref{eq:Gn_def_hash} is a centered Gaussian process.
Specifically, its conditional covariance function $\Gamma_{nj}$ is given by
\begin{equation*}
	\begin{split}
		\Gamma_{nj}(x_j,x_j') :=& \, \Cov\big(G_{nj}(x_j),G_{nj}(x_j') \, \big| \, \cX_j\big) \\
		=& \, \frac{1}{n h_j\sqrt{p_j(x_j)p_j(x_j')}\sigma_j(x_j)\sigma_j(x_j')} \sumin K\bigg(\frac{x_j-X_{ij}}{h_j}\bigg)K\bigg(\frac{x_j'-X_{ij}}{h_j}\bigg) \sigma_j^2(X_{ij}).
	\end{split}	
\end{equation*}
Let $C_{nj}(x_j,x_j'):=(K*K)\big((x_j-x_j')/h_j\big)$.
The following lemma establishes the uniform approximation error between $\Gamma_{nj}$ and $C_{nj}$ over $I\times I$.
The proof of Lemma~\ref{lem:Cov_ftn_approx} is provided in Section~\ref{subsec:proof_Cov_ftn_approx}.

\begin{lemma}\label{lem:Cov_ftn_approx}
	Assume the conditions in Theorem~\ref{thm:unif_CLT}.
	Then, for any closed interval $I\subset (0,1)$, it holds that
	\begin{equation}\label{eq:Delta_n_rate}
		\sup_{x_j,x_j'\in I} \big|  \Gamma_{nj}(x_j,x_j')-C_{nj}(x_j,x_j') \big| = O_p\bigg( h_j+\sqrt{\frac{\log n}{n h_j}} \bigg).
	\end{equation}
\end{lemma}

\begin{proof}[Proof of Lemma~{\upshape\ref{lem:gaussian_process_approx}}]
	We employ a finite approximation argument for Gaussian suprema.
	First, by applying Dudley's inequality \citep{dudley2016suprema} and the Gaussian comparison inequality of Theorem 2 in \cite{chernozhukov2015comparison}, we obtain the desired comparison for the maxima of the Gaussian processes restricted to a finite grid.
	Then, the extension from the finite maxima to the continuous suprema is justified using the Gaussian anti-concentration inequality for suprema \citep{chernozhukov2014anti} to control the associated discretization error.

	Let $a_{nj}:=h_j+\sqrt{\frac{\log n}{n h_j}}$ and $\delta_{nj}:=\sqrt{a_{nj}}/\log n$.
	Since $h_j\asymp n^{-b}$, it follows that
	\begin{equation}\label{eq:delta_condition}
		\delta_{nj}(\log n)^2 =\sqrt{a_{nj}}\log n\to0, \qquad \frac{\delta_{nj}}{a_{nj}} =\frac{1}{\sqrt{a_{nj}}\log n}\to\infty \qquad \text{as} \quad n\to\infty.
	\end{equation}
	Define the event
	\begin{equation}\label{eq:event_nj}
		\cE_{nj}:=\big\{\Delta_{nj}\le \delta_{nj}\big\},
	\end{equation}
	where
	\begin{equation*}
		\Delta_{nj} := \sup_{x_j,x_j'\in I} \big|  \Gamma_{nj}(x_j,x_j')-C_{nj}(x_j,x_j') \big|.
	\end{equation*}
	By Lemma~\ref{lem:Cov_ftn_approx} and \eqref{eq:delta_condition}, the event $\cE_{nj}$ occurs with probability tending to one.
	Hence, it suffices to prove the desired assertion on the event $\cE_{nj}$.

	We first discretize the closed interval $I\subset(0,1)$.
	Let $\rho_{nj}:=(\log n)^{-(2+\gamma)}$ for some $\gamma>0$.
	Let $\cI_{nj}:=\{x_{1j},\ldots,x_{N_{nj}j}\}\subset I$ denote a grid such that every point of $I$ is within a distance of $h_j\rho_{nj}$ from some point in $\cI_{nj}$.
	Since $I$ is a fixed interval, we can choose $\cI_{nj}$ satisfying
	\begin{equation}\label{eq:net_condition}
		N_{nj}\le \frac{C_1}{h_j\rho_{nj}},	\qquad	\log N_{nj}\le C_2\log n,
	\end{equation}
	for some constants $C_1,C_2>0$.
	Let
	\begin{equation*}
		S_{nj}^W:=\sup_{x_j\in I}\abs{W_{nj}(x_j)}, \qquad S_{nj}^G:=\sup_{x_j\in I}\abs{G_{nj}(x_j)},
	\end{equation*}
	and define their discretized versions
	\begin{equation*}
		S_{nj}^{W,d}:=\max_{1\le k\le N_{nj}}\abs{W_{nj}(x_{kj})}, \qquad S_{nj}^{G,d}:=\max_{1\le k\le N_{nj}}\abs{G_{nj}(x_{kj})}.
	\end{equation*}
	For any $x_j\in I$, let $\pi_{nj}(x_j)\equiv \pi_{nj}(x_j;\cI_{nj})$ be a nearest grid point to $x_j$ within $\cI_{nj}$.
	Define
	\begin{equation*}
		\omega_{nj}^W := \sup_{x_j\in I} \abs{W_{nj}(x_j)-W_{nj}\big(\pi_{nj}(x_j)\big)}, \qquad \omega_{nj}^G := \sup_{x_j\in I} \abs{G_{nj}(x_j)-G_{nj}\big(\pi_{nj}(x_j)\big)}.
	\end{equation*}
	Then, we have
	\begin{equation*}
		S_{nj}^{W,d}\le S_{nj}^W\le S_{nj}^{W,d}+\omega_{nj}^W, \qquad S_{nj}^{G,d}\le S_{nj}^G\le S_{nj}^{G,d}+\omega_{nj}^G.
	\end{equation*}
	We claim that there exists a sufficiently large constant $C>0$ such that, with $u_{nj}:=C\sqrt{\big(\rho_{nj}+\delta_{nj}\big)\log n}$, it holds that
	\begin{equation}\label{eq:approx_claim}
		\Prob(\omega_{nj}^G>u_{nj}\mid \cX_j) + \Prob(\omega_{nj}^W>u_{nj}) = o(1)
	\end{equation}
	on the event $\cE_{nj}$ at \eqref{eq:event_nj}.
	This gives the discretization errors for $S_{nj}^W$ and $S_{nj}^G$.

	We now compare the maximum over the discretized set.
	Define the $2N_{nj}$-dimensional Gaussian vectors
	\begin{align*}
		\bbG_{nj} :=& \, \big(G_{nj}(x_{1j}),\ldots,G_{nj}(x_{N_{nj}j}),-G_{nj}(x_{1j}),\ldots,-G_{nj}(x_{N_{nj}j})\big)^\top, \\
		\bbW_{nj} :=& \, \big(W_{nj}(x_{1j}),\ldots,W_{nj}(x_{N_{nj}j}),-W_{nj}(x_{1j}),\ldots,-W_{nj}(x_{N_{nj}j})\big)^\top.
	\end{align*}
	Then, we can write
	\begin{equation*}
		S_{nj}^{W,d} = \max_{1\le k\le 2N_{nj}}\bbW_{nj,k}, \qquad S_{nj}^{G,d}	= \max_{1\le k\le 2N_{nj}}\bbG_{nj,k}.
	\end{equation*}
	On the event $\cE_{nj}$ at \eqref{eq:event_nj}, it holds that
	\begin{equation*}
		\max_{1\le k,l\le 2N_{nj}} \abs{\Cov(\bbG_{nj,k},\bbG_{nj,l}\mid \cX_j) - \Cov(\bbW_{nj,k},\bbW_{nj,l})} \le	\delta_{nj}.
	\end{equation*}
	Therefore, by Theorem 2 in \cite{chernozhukov2015comparison}, on the event $\cE_{nj}$ at \eqref{eq:event_nj}, we obtain
	\begin{equation}\label{eq:approx_2}
		D_{nj} := \sup_{t\in\bbR} \abs{\Prob(S_{nj}^{G,d}\le t\mid \cX_j) - \Prob(S_{nj}^{W,d}\le t)}  \lesssim \delta_{nj}^{1/3} \bigg\{ 1\vee \log\bigg(\frac{2N_{nj}}{\delta_{nj}}\bigg)\bigg\}^{2/3} = o(1).
	\end{equation}

	It remains to extend the result from the discretized set to the full interval $I$.
	Following the proof of the claim \eqref{eq:approx_claim}, we can show that
	\begin{equation*}
		\E[S_{nj}^W] = \E\bigg[\sup_{x_j\in I}\abs{W_{nj}(x_j)}\bigg] \lesssim \sqrt{\log n}.
	\end{equation*}
	Using this bound and Corollary 2.1 in \cite{chernozhukov2014anti}, we obtain
	\begin{equation}\label{eq:anti-concen}
		\sup_{t\in\bbR} \Prob\big(\abs{S_{nj}^W-t}\le 2u_{nj}\big) \le 8u_{nj}\{1+\E[S_{nj}^W]\} = o(1),
	\end{equation}
	where the last equality holds as $u_{nj}(\log n)^{1/2}\to0$, which follows from $\rho_{nj}=(\log n)^{-(2+\gamma)}$ and $\delta_{nj}(\log n)^2\to0$.

	Finally, for every $t\in\bbR$,
	\begin{equation*}
		\begin{split}
			\Prob(S_{nj}^G\le t\mid \cX_j) &\le \Prob(S_{nj}^{G,d}\le t\mid \cX_j)  \\
			&\le \Prob(S_{nj}^{W,d}\le t)+D_{nj}  \\
			&\le \Prob(S_{nj}^W\le t+u_{nj}) + \Prob(\omega_{nj}^W>u_{nj}) + D_{nj}
		\end{split}
	\end{equation*}
	and similarly,
	\begin{equation*}
		\begin{split}
			\Prob(S_{nj}^G\le t\mid \cX_j) &\ge \Prob(S_{nj}^{G,d}\le t-u_{nj}\mid \cX_j) - \Prob(\omega_{nj}^G>u_{nj}\mid \cX_j)  \\
			&\ge \Prob(S_{nj}^{W,d}\le t-u_{nj}) - D_{nj} - \Prob(\omega_{nj}^G>u_{nj}\mid \cX_j)  \\
			&\ge \Prob(S_{nj}^W\le t-u_{nj}) - D_{nj} - \Prob(\omega_{nj}^G>u_{nj}\mid \cX_j).
		\end{split}
	\end{equation*}
	By the above two results, on the event $\cE_{nj}$ at \eqref{eq:event_nj}, we have
	\begin{equation}\label{eq:G_W_prob_bdd}
		\begin{split}
			&\sup_{t\in\bbR} \abs{ \Prob(S_{nj}^G\le t\mid \cX_j) - \Prob(S_{nj}^W\le t)}  \\
			&\qquad\le D_{nj} + \Prob(\omega_{nj}^G>u_{nj}\mid \cX_j) + \Prob(\omega_{nj}^W>u_{nj})  + \sup_{t\in\bbR} \Prob\big(\abs{S_{nj}^W-t}\le 2u_{nj}\big).
		\end{split}
	\end{equation}
	Combining \eqref{eq:approx_claim}--\eqref{eq:G_W_prob_bdd}, we complete the proof of Lemma~\ref{lem:gaussian_process_approx}.

	\underline{\bf Proof of the claim \eqref{eq:approx_claim}}. \quad
	For $W_{nj}$, define the pseudo-metric
	\begin{equation}\label{eq:pseudometric_on_gaussian}
		d_W(x_j,x_j') := \Big[ \E\{W_{nj}(x_j)-W_{nj}(x_j')\}^2 \Big]^{1/2}.
	\end{equation}
	Since $K*K$ is Lipschitz continuous by (A2), there exists a constant $C_1>0$ such that, for all $\abs{x_j-x_j'}\le h_j\rho_{nj}$,
	\begin{equation}\label{eq:d_W}
		d_W^2(x_j,x_j') = 2(K*K)(0) - 2 (K*K)\bigg(\frac{x_j-x_j'}{h_j}\bigg) \le C_1\rho_{nj}.
	\end{equation}
	Similarly, define the conditional version of $d_W(x_j,x_j')$ as
	\begin{equation*}
		d_G(x_j,x_j') := \Big[ \E\big[ \{G_{nj}(x_j)-G_{nj}(x_j')\}^2 \,\big|\, \cX_j \big] \Big]^{1/2}.
	\end{equation*}
	Then, there exists a constant $C_2>0$ such that, on the event $\cE_{nj}$ at \eqref{eq:event_nj}, for all $\abs{x_j-x_j'}\le h_j\rho_{nj}$,
	\begin{equation*}
		\begin{split}
			d_G^2(x_j,x_j') &= \Gamma_{nj}(x_j,x_j)+\Gamma_{nj}(x_j',x_j') -2\Gamma_{nj}(x_j,x_j')  \\
			&\le C_{nj}(x_j,x_j)+C_{nj}(x_j',x_j') -2C_{nj}(x_j,x_j')+4\Delta_{nj}  \\
			&= 2(K*K)(0) -2 (K*K)\bigg(\frac{x_j-x_j'}{h_j}\bigg) + 4\Delta_{nj}  \\
			&\le C_2\big(\rho_{nj}+\delta_{nj}\big).
		\end{split}
	\end{equation*}

	Let $T_{nj}^{\pm}:=I\times\{-1,1\}$ denote the signed index set.
	For $(x_j,s)\in T_{nj}^{\pm}$, define
	\begin{align*}
		Y_{nj}^W(x_j,s):=&\, s \cdot \Big(W_{nj}(x_j)-W_{nj}\big(\pi_{nj}(x_j)\big)\Big), \\
		Y_{nj}^G(x_j,s):=&\, s \cdot \Big(G_{nj}(x_j)-G_{nj}\big(\pi_{nj}(x_j)\big)\Big).
	\end{align*}
	Then, $Y_{nj}^W$ and $Y_{nj}^G \, | \, \cX_j$ are centered Gaussian processes.
	Furthermore, we can write
	\begin{equation*}
		\omega_{nj}^W=\sup_{(x_j,s)\in T_{nj}^{\pm}}Y_{nj}^W(x_j,s), \qquad \omega_{nj}^G=\sup_{(x_j,s)\in T_{nj}^{\pm}}Y_{nj}^G(x_j,s).
	\end{equation*}
	Let $d_{Y,W}$ be the canonical pseudo-metric associated with $Y_{nj}^W$, defined as in \eqref{eq:pseudometric_on_gaussian}:
	\begin{equation*}
		d_{Y,W}\big((x_j,s),(x_j',s')\big):= \Big[ \E\{Y_{nj}^W(x_j,s) - Y_{nj}^W(x_j',s')\}^2 \Big]^{1/2}
	\end{equation*}
	Since $\abs{x_j-\pi_{nj}(x_j)}\le h_j\rho_{nj}$, \eqref{eq:d_W} implies $\sup_{x_j\in I}\Var\{Y_{nj}^W(x_j)\}=O(\rho_{nj})$, and hence,
	\begin{equation}\label{eq:diam_T_nj}
		\operatorname{diam}\big(T_{nj}^{\pm},d_{Y,W}\big) \lesssim \sqrt{\rho_{nj}}.
	\end{equation}

	We now bound the covering number of $T_{nj}^{\pm}$ under $d_{Y,W}$.
	For $1\le k\le N_{nj}$, let $I_{k,nj}:=\{x_j\in I:\pi_{nj}(x_j)=x_{kj}\}$ be the cell induced by the grid $x_{kj}\in\cI_{nj}$.
	Then $\abs{I_{k,nj}} = O(h_j\rho_{nj})$.
	For any $x_j,x_j'\in I_{k,nj}$ and a fixed $s\in\{-1,1\}$, we have
	\begin{equation*}
		Y_{nj}^W(x_j,s)-Y_{nj}^W(x_j',s) = s\{W_{nj}(x_j)-W_{nj}(x_j')\},
	\end{equation*}
	and therefore
	\begin{equation*}
		d_{Y,W}^2\big((x_j,s),(x_j',s)\big) = d_W^2(x_j,x_j') \le C_3\frac{\abs{x_j-x_j'}}{h_j}.
	\end{equation*}
	Thus, for each $u>0$, the number of $d_{Y,W}$-balls of radius $u$ to cover the set $I_{k,nj}\times\{s\}$ is $O\big(1+\rho_{nj}/ u^2\big)$.
	Summing over all $N_{nj}$ cells and the two signs, and using \eqref{eq:net_condition}, we get that, for any $0<u\le C_1'\sqrt{\rho_{nj}}$,
	\begin{equation}\label{eq:covering_number_bound}
		N\big(u,T_{nj}^{\pm},d_{Y,W}\big) \le C_2'N_{nj}\big(1+\frac{\rho_{nj}}{u^2}\big) \le \frac{C_3'}{h_j u^2}
	\end{equation}
	for some constants $C_1',C_2',C_3'>0$.
	Applying Dudley's entropy integral bound \citep{dudley2016suprema}, together with \eqref{eq:diam_T_nj} and \eqref{eq:covering_number_bound}, we obtain
	\begin{equation}\label{eq:omega^W bound}
		\E[\omega_{nj}^W] = \E\bigg[\sup_{(x_j,s)\in T_{nj}^{\pm}}Y_{nj}^W(x_j,s)\bigg] \le 24\int_0^{C_1'\sqrt{\rho_{nj}}} \sqrt{\log\bigg(\frac{C_3'}{h_j u^2}\bigg)}	\,\d u \lesssim \sqrt{\rho_{nj}\log n},
	\end{equation}
	where the last inequality holds as $\rho_{nj}=(\log n)^{-(2+\gamma)}$ and $h_j\asymp n^{-b}$.

	Analogously, we can show that, on the event $\cE_{nj}$ at \eqref{eq:event_nj},
	\begin{equation}\label{eq:omega_G bound}
		\E[\omega_{nj}^G\mid \cX_j] \lesssim \sqrt{\big(\rho_{nj}+\delta_{nj}\big)\log n}.
	\end{equation}
	Combining the preceding two bounds \eqref{eq:omega^W bound} and \eqref{eq:omega_G bound}, and applying Theorem 2.1.1 in \cite{adler2007random}, we obtain the claim \eqref{eq:approx_claim}.
\end{proof}

\subsection{Proof of Lemma~\ref{lem:Cov_ftn_approx}}\label{subsec:proof_Cov_ftn_approx}

Define the function class $\cF_{nj}:= \{f_{x_j,x_j'}: x_j,x_j'\in I\}$, where
\begin{equation*}
	f_{x_j,x_j'}(v_j) := K\bigg(\frac{x_j-v_j}{h_j}\bigg)K\bigg(\frac{x_j'-v_j}{h_j}\bigg) \sigma_j^2(v_j), \quad v_j \in I.
\end{equation*}
By (A2) and the fact that each $f_{x_j,x_j'}$ is supported on an interval of length $O(h_j)$, the class $\cF_{nj}$ admits a bounded envelope $F$ such that
\begin{equation}\label{eq:cF_prop}
	\sigma_{\cF_{nj}}^2 := \sup_{x_j,x_j' \in I} \E f_{x_j,x_j'}^2(X_j) \lesssim h_j, \quad \|F \|_{P,2} \lesssim 1, \quad \|M(F) \|_2 \lesssim 1.
\end{equation}
We next verify a covering number bound for $\cF_{nj}$.
By the Lipschitz continuity of $K$, there exists a constant $C>0$ such that, for every probability measure $Q$ and every $x_j,x_j',\tilde x_j,\tilde x_j'\in I$,
\begin{equation*}
	\big\| f_{x_j,x_j'}-f_{\tilde x_j,\tilde x_j'} \big\|_{Q,2}	\le	C h_j^{-1}\{|x_j-\tilde x_j| + |x_j'-\tilde x_j'|\} \cdot \norm{F}_{Q,2}.
\end{equation*}
Moreover, the parameter space $I^2$ can be covered by at most $C_1\eta^{-2}$ balls of radius $\eta$ under the metric $d\{(x_j,x_j'),(\tilde x_j,\tilde x_j')\}:=\abs{x_j-\tilde x_j}+\abs{x_j'-\tilde x_j'}$ for some constant $C_1>0$.
Therefore, by choosing $\eta\asymp u h_j$, we obtain
\begin{equation*}
	\sup_Q N\big(u \norm{F}_{Q,2},\cF_{nj},L_2(Q)\big) \le C_1\bigg(\frac{C}{u h_j}\bigg)^2 \le \frac{C_1C^2}{(u h_j)^2},	\quad 0<u \le 1,
\end{equation*}
where the supremum runs over all probability measure $Q$ on $[0,1]$.
By \eqref{eq:cF_prop}, applying Corollary 5.1 of \cite{chernozhukov2014gaussian} with $A=\sqrt{C_1}C/h_j$ and $\sigma^2=C'h_j$ satisfying $\sigma_{\cF_{nj}}^2 \le \sigma^2 \le \norm{F}_{P,2}^2$ yields
\begin{equation*}
	\E \norm{\bbG_n}_{\cF_{nj}}	\lesssim \sqrt{h_j\log(1/h_j)} + \frac{\log(1/h_j)}{\sqrt n}.
\end{equation*}
Hence, by Markov's inequality and the fact that $nh_j/\log(1/h_j)\to\infty$,
\begin{equation*}
	\sup_{x_j,x_j'\in I}\abs{(P_n-P)f_{x_j,x_j'}} = O_p\bigg(\sqrt{\frac{h_j\log(1/h_j)}{n}} + \frac{\log(1/h_j)}{n}\bigg) = O_p\bigg(\sqrt{\frac{h_j\log n}{n}}\bigg),
\end{equation*}
or equivalently,
\begin{equation}\label{eq:empirical_expectation_diff_bound}
	\sup_{x_j,x_j'\in I} \abs{ \meanin f_{x_j,x_j'}(X_{ij})	- \E f_{x_j,x_j'}(X_j) } = O_p\bigg( \sqrt{\frac{h_j\log n}{n}} \bigg).
\end{equation}

It remains to obtain a uniform bound of $\E f_{x_j,x_j'}(X_j)$ over $x_j,x_j'\in I$.
Observe that
\begin{equation*}
	h_j^{-1}\E f_{x_j,x_j'}(X_j)	= \int K(u)	K\bigg(u+\frac{x_j-x_j'}{h_j}\bigg) g_j(x_j-h_j u) \, \d u,
\end{equation*}
where $g_j(x_j):=p_j(x_j)\sigma_j^2(x_j)$.
Since $K$ is supported on $[-1,1]$, it holds that $\E f_{x_j,x_j'}(X_j) = 0$ whenever $|x_j-x_j'|>2h_j$.
Therefore, we suppose $\abs{x_j-x_j'}\le 2h_j$.
By (B1) and (B2), $g_j$ is continuously differentiable on an open neighborhood of $I$ and $\inf_{x_j\in I}g_j(x_j)>0$.
Thus,
\begin{equation}\label{eq:expectation_cov_approx_1}
	\begin{split}
		h_j^{-1}\E f_{x_j,x_j'}(X_j) &= g_j(x_j) \int K(u) K\bigg(u+\frac{x_j-x_j'}{h_j}\bigg)\, \d u + O(h_j) \\
		&=	g_j(x_j) (K*K)\bigg(\frac{x_j-x_j'}{h_j}\bigg) + O(h_j),
	\end{split}
\end{equation}
where the $O$ term holds uniformly over $x_j,x_j'\in I$.
We next approximate the normalizing denominator.
By (B1) and (B2), it holds that $\sqrt{p_j(x_j')}\,\sigma_j(x_j')=\sqrt{p_j(x_j)}\,\sigma_j(x_j)+O(h_j)$ uniformly over $x_j,x_j'\in I$ satisfying $\abs{x_j-x_j'}\le 2h_j$.
It then follows that
\begin{equation*}
	\sqrt{p_j(x_j)p_j(x_j')}\,\sigma_j(x_j)\sigma_j(x_j') =	p_j(x_j)\sigma_j^2(x_j)+O(h_j)	= g_j(x_j)+O(h_j).
\end{equation*}
Since $g_j$ is bounded away from zero on $I$, \eqref{eq:expectation_cov_approx_1} implies
\begin{equation}\label{eq:expectation_cov_approx_2}
	\frac{1}{h_j\sqrt{p_j(x_j)p_j(x_j')}\sigma_j(x_j)\sigma_j(x_j')} \E f_{x_j,x_j'}(X_j) = (K*K)\bigg(\frac{x_j-x_j'}{h_j}\bigg) + O(h_j),
\end{equation}
uniformly over $x_j,x_j'\in I$.
Combining \eqref{eq:empirical_expectation_diff_bound} and \eqref{eq:expectation_cov_approx_2} concludes the proof of Lemma~\ref{lem:Cov_ftn_approx}.
\qed

\end{document}